\documentclass[12pt]{amsart}
\usepackage[]{amscd,amsfonts,amsmath,amsxtra,amssymb,bm}
\usepackage{graphics}
\usepackage{extarrows}
\usepackage[all]{xy}
\usepackage{hyperref}
\usepackage{eucal}
\usepackage{rotating}
\usepackage[T1]{fontenc}
\newtheorem{sub}{}[section]
\newtheorem{subsub}{}[sub]
\usepackage{rotating}
\usepackage[active]{srcltx}
\def\ov#1{\overline{#1}}

\def\Hom{\mathop{\rm Hom}\nolimits}

\def\HHom{\mathop{\mathcal Hom}\nolimits}

\def\Ext{\mathop{\rm Ext}\nolimits}
\def\EExt{\mathop{\mathcal Ext}\nolimits}
\def\Tor{\mathop{\rm Tor}\nolimits}

\def\Pic{\mathop{\rm Pic}\nolimits}

\def\Aut{\mathop{\rm Aut}\nolimits}

\def\End{\mathop{\rm End}\nolimits}

\def\EEnd{\mathop{\mathcal End}\nolimits}
\def\AAd{\mathop{\mathcal Ad}\nolimits}
\def\GL{\mathop{\rm GL}\nolimits}

\def\imm{\mathop{\rm im}\nolimits}
\def\deg{\mathop{\rm deg}\nolimits}

\def\rk{\mathop{\rm rk}\nolimits}

\def\det{\mathop{\rm det}\nolimits}
\def\spec{\mathop{\rm spec}\nolimits}

\def\lra{\longrightarrow}

\def\sigg{\mathop{\hbox{$\displaystyle\sum$}}\limits}

\def\hfl#1#2{\smash{\mathop{\ \hbox to 12mm{\rightarrowfill}}
\limits^{\scriptstyle#1}_{\scriptstyle#2} \ }}
\def\hflb#1#2{\smash{\mathop{\hbox to 12mm{\leftarrowfill}}
\limits^{\scriptstyle#1}_{\scriptstyle#2}}}

\def\m#1{{\hbox{$#1$}}}

\def\ot{\otimes}
\def\og{\leavevmode\raise.3ex\hbox{$\scriptscriptstyle\langle\!\langle$}}
\def\fg{\leavevmode\raise.3ex\hbox{$\scriptscriptstyle\,\rangle\!\rangle$}}

\def\nsp{\lbrace 0\rbrace}

\def\dsp{\displaystyle}
\def\Ssect#1#2{\pagebreak[3]\begin{sub}\label{#2}{\sc\small\small
#1}\rm\medskip}

\def\sepsec{\vskip 1.4cm}
\def\sepsub{\vskip 0.7cm}
\def\sepsubsub{\vskip 0.5cm}
\def\sepprop{\vskip 0.5cm}

\def\xmat#1{\[\xymatrix{#1}\]}
\def\flinc{\ar@{^{(}->}}
\def\fleq{\ar@{=}}
\def\flon{\ar@{->>}}
\def\fmaps{\ar@{|-{>}}}
\def\fflat{\ar@{-}}
\def\fimpl{\ar@{=>}}
\def\Nligne{\hfil\break}

\def\ED{\vskip 1cm\end{document}}

\def\BS#1{{\boldsymbol{#1}}}

\newcommand{\A}{{\mathbb A}}
\newcommand{\B}{{\mathbb B}}
\newcommand{\Z}{{\mathbb Z}}

\newcommand{\C}{{\mathbb C}}

\renewcommand{\P}{{\mathbb P}}
\newcommand{\D}{{\mathbb D}}
\newcommand{\F}{{\mathbb F}}
\newcommand{\E}{{\mathbb E}}
\newcommand{\G}{{\mathbb G}}

\newcommand{\U}{{\mathbb U}}
\newcommand{\V}{{\mathbb V}}
\newcommand{\W}{{\mathbb W}}
\renewcommand{\L}{{\mathbb L}}
\renewcommand{\H}{{\mathbb H}}
\def\T{{\mathbb T}}
\newcommand{\ka}{{\mathcal A}}
\newcommand{\kb}{{\mathcal B}}

\newcommand{\ke}{{\mathcal E}}
\newcommand{\kf}{{\mathcal F}}
\newcommand{\kg}{{\mathcal G}}
\newcommand{\kh}{{\mathcal H}}
\newcommand{\ki}{{\mathcal I}}

\newcommand{\kl}{{\mathcal L}}
\newcommand{\km}{{\mathcal M}}
\newcommand{\kn}{{\mathcal N}}
\newcommand{\ko}{{\mathcal O}}
\newcommand{\kp}{{\mathcal P}}
\newcommand{\kq}{{\mathcal Q}}

\newcommand{\kv}{{\mathcal V}}

\begin{document}

\def\refname{References}
\def\contentsname{Summary}
\def\proofname{Proof}
\def\abstractname{Resume}

\author{Jean--Marc Dr\'{e}zet}
\address{
Institut de Math\'ematiques de Jussieu - Paris Rive Gauche\\
Case 247\\
4 place Jussieu\\
F-75252 Paris, France}
\email{jean-marc.drezet@imj-prg.fr}
\title[{Moduli of vector bundles}] {Moduli of vector bundles on primitive 
multiple schemes}

\begin{abstract}
A primitive multiple scheme is a Cohen-Macaulay scheme $Y$ such that the 
associated reduced scheme $X=Y_{red}$ is smooth, irreducible, and that $Y$ can 
be locally embedded in a smooth variety of dimension \m{\dim(X)+1}. If $n$ is 
the multiplicity of $Y$, there is a canonical filtration \m{ X=X_1\subset
X_2\subset\cdots\subset X_n=Y}, such that \m{X_i} is a primitive multiple 
scheme of multiplicity $i$. The simplest example is the trivial primitive 
multiple scheme of multiplicity $n$ associated to a line bundle $L$ on $X$: it 
is the $n$-th infinitesimal neighborhood of $X$, embedded in the line bundle 
$L^*$ by the zero section.

The main subject of this paper is the construction and properties of fine 
moduli spaces of vector bundles on primitive multiple schemes. Suppose that 
$Y=X_n$ is of multiplicity $n$, and can be extended to $X_{n+1}$ of 
multiplicity $n+1$, and let $M_n$ be a fine moduli space of vector bundles on 
$X_n$. With suitable hypotheses, we construct a fine moduli space $M_{n+1}$ for 
the vector bundles on $X_{n+1}$ whose restriction to $X_n$ belongs to $M_n$. It 
is an affine bundle over the subvariety $N_n\subset M_n$ of bundles that can be 
extended to $X_{n+1}$. In general this affine bundle is not banal. This applies 
in particular to Picard groups.

We give also many new examples of primitive multiple schemes $Y$ such that the 
dualizing sheaf $\omega_Y$ is trivial.
\end{abstract}

\maketitle
\tableofcontents


\vskip 1cm

\section{Introduction}\label{intro}

In this paper a {\em scheme} is a noetherian separated scheme over $\C$.
A {\em primitive multiple scheme} is a Cohen-Macaulay scheme $Y$ over $\C$ such 
that:
\begin{enumerate}
\item[--] $Y_{red}=X$ is a smooth connected variety, 
\item[--] for every closed point \m{x\in X}, there exists a neighbourhood $U$ 
of $x$ in $Y$, and a smooth variety $S$ of dimension \ \m{\dim(X)+1} such that 
$U$ is isomorphic to a closed subscheme of $S$.
\end{enumerate}
We call $X$ the {\em support} of $Y$. It may happen that $Y$ is 
quasi-projective, and in this case it is projective if $X$ is.

For every closed subscheme \m{Z\subset Y}, let \m{\ki_Z} (or \m{\ki_{Z,Y})} 
denote the ideal sheaf of $Z$ in $Y$. For every positive integer $i$, 
let \m{X_i} be the closed subscheme of $Y$ corresponding to the ideal sheaf 
\m{\ki_X^i}. The smallest integer $n$ such that \ \m{X_n=Y} \ is called the 
{\em multiplicity} of $Y$. For \m{1\leq i\leq n}, \m{X_i} is a primitive 
multiple scheme of multiplicity $i$, \m{L=\ki_X/\ki_{X_2}} is a line bundle on 
$X$, and we have \ \m{\ki_{X_{i}}/\ki_{X_{i+1}}=L^i}. 
We call $L$ the line bundle on $X$ {\em associated} to $Y$. The ideal sheaf 
\m{\ki_{X}} can be viewed as a line bundle on \m{X_{n-1}}. If \m{n=2}, $Y$ is 
called a {\em primitive double scheme}.

The simplest case is when $Y$ is contained in a smooth variety $S$ of dimension 
\m{\dim(X)+1}. Suppose that $Y$ has multiplicity $n$. Let \m{P\in X} and 
\m{f\in\ko_{S,P}}  a local equation of $X$. Then we have \ 
\m{\ki_{X_i,P}=(f^{i})} \ for \m{1<j\leq n} in $S$, in particular 
\m{\ki_{Y,P}=(f^n)}, and \ \m{L=\ko_X(-X)} .

For any \m{L\in\Pic(X)}, the {\em trivial primitive variety} of multiplicity 
$n$, with induced smooth variety $X$ and associated line bundle $L$ on $X$ is 
the $n$-th infinitesimal neighborhood of $X$, embedded by the zero section in 
the dual bundle $L^*$, seen as a smooth variety.

The primitive multiple curves where defined in \cite{fe}, \cite{ba_fo}. 
Primitive double curves were parametrized and studied in \cite{ba_ei} and 
\cite{ei_gr}. More results on primitive multiple curves can be found in 
\cite{dr2}, \cite{dr1}, \cite{dr4}, \cite{dr5}, \cite{dr6}, \cite{dr7}, 
\cite{dr8}, \cite{dr9}, \cite{ch-ka}, \cite{sa1}, \cite{sa2}, \cite{sa3}. Some 
primitive double schemes are studied in \cite{b_m_r}, \cite{ga-go-pu} and 
\cite{gonz1}. The case of varieties of any dimension is studied in \cite{dr10}, 
where the following subjects were treated:
\begin{enumerate}
\item[--] construction and parametrization of primitive multiple schemes,
\item[--] obstructions to the extension of a vector bundle on \m{X_m} to 
\m{X_{m+1}},
\item[--] obstructions to the extension of a primitive multiple scheme of 
multiplicity $n$ to one of multiplicity \m{n+1}.
\end{enumerate}
The main subject of this paper is the construction and properties of moduli 
spaces of vector bundles on \m{X_n}. If \m{X_{n+1}} is an extension of \m{X_n} 
to a primitive multiple scheme of multiplicity \m{n+1}, and $\bf M$ is a moduli 
space of vector bundles on \m{X_n}, we will see how to construct a moduli space 
\m{{\bf M}'} for the vector bundles on \m{X_{n+1}} whose restriction to \m{X_n} 
belongs to $\bf M$. It is an affine bundle over the closed subvariety \m{{\bf 
N}\subset{\bf M}} of bundles that can be extended to \m{X_{n+1}}. This applies 
in particular to Picard groups.

\sepsub

\Ssect{Fine moduli spaces}{F_M}

Let $\chi$ be a set of isomorphism classes of vector bundles on a scheme $Z$. 
Suppose that $\chi$ is {\em open}, i.e. for every scheme $V$ and vector bundle 
$\E$ on \m{Z\times V}, if \m{v\in V} is a closed point such that 
\m{\E_v\in\chi}, then there exists an open neighbourhood $U$ of $v$ such that 
\m{\E_u\in\chi} for every closed point \m{u\in U}.
A {\em fine moduli space for $\chi$ }is the data of a scheme $M$ and of
\begin{enumerate}
\item[--] a bijection
\[\xymatrix@R=5pt{M^0\ar[r] & \chi\\ m\fmaps[r] & E_m
}\]
(where \m{M^0} denotes the set of closed points of $M$),
\item[--] an open cover \m{(M_i)_{i\in I}} of $M$, and for every \m{i\in I}, a 
vector bundle \m{\ke_i} on \m{X\times M_i} such that for every \m{m\in M_i}, 
\m{\ke_{i,m}\simeq E_m},
\end{enumerate}
such that:
for any scheme $S$, any vector bundle $\kf$ on \m{Z\times S} such that for any 
closed point \m{s\in S}, \m{\kf_s\in\chi}, there is a morphism \ 
\m{f_\kf:S\to M} \ such that for every \m{s\in S}, if \m{m=f_\kf(s)} then 
\m{\kf_s\simeq E_m}, and if \m{m\in M_i}, then there exists an open 
neighbourhood \m{U} of $s$ such that \ \m{f_\kf(U)\subset M_i} and
\ \m{(I_X\times f_{\kf|U})^*(\ke_i)\simeq\kf_{|X\times U}}.

\sepsubsub

\begin{subsub}\label{111} Extensions of a vector bundle to higher multiplicity 
-- \rm Let $\E$ be a vector bundle on \m{X_n} and \ \m{E=\E_{|X}}. In 
\cite{dr10}, 7.1, a class \ \m{\Delta(\E)\in\Ext^2_{\ko_X}(E,E\ot L^n)} \ is 
defined, such that $\E$ can be extended to a vector bundle on \m{X_{n+1}} if 
and only if \ \m{\Delta(\E)=0}.

If $\ke$ is a vector bundle on \m{X_{n+1}} such that \ \m{\ke_{|X_n}=\E}, the 
kernel of the restriction morphism \ \m{\ke\to\E} \ is isomorphic to \m{E\ot 
L^n}, and we have an exact sequence on \m{X_{n+1}}
\begin{equation}\label{equ24}
0\lra E\ot L^n\lra\ke\lra\E\lra 0 \ . \end{equation}
Hence if we want to extend $\E$ to \m{X_{n+1}}, it is natural to study \ 
\m{\Ext^1_{\ko_{X_{n+1}}}(\E,E\ot L^n)}. Using the local-to-global Ext spectral 
sequence (cf. \cite{go}, Corollary of Theorem 7.3.3), and the fact that \ 
\m{\EExt^1_{\ko_{X_{n+1}}}(\E,E\ot L^n)\simeq\EEnd(E)}, we obtain a canonical 
exact sequence
\xmat{0\ar[r] & \Ext^1_{\ko_X}(E,E\ot L^n)\ar[r]^-\psi & 
\Ext^1_{\ko_{X_{n+1}}}(\E,
E\ot L^n)\ar[r]^-\phi & \End(E)\ar[r]^-{\delta_\E} & \Ext^2_{\ko_X}(E,E\ot L^n) 
\ , }
and we have, from \cite{dr10}, \m{\Delta(\E)=\delta_\E(I_E)}. In \cite{dr10} is 
also given the following description of \m{\Delta(\E)}: we have a canonical 
exact sequence of sheaves on \m{X_n}: \m{0\to 
L^n\to\Omega_{X_{n+1}|X_n}\to\Omega_{X_n}\to 0}, associated to \ 
\m{\sigma_{\Omega_{X_{n+1}}}\in\Ext^1_{\ko_{X_n}}(\Omega_{X_n},L^n)}, inducing
\ \m{\ov{\sigma}_{\E,X_{n+1}}\in
\Ext^1_{\ko_{X_n}}(\E\ot\Omega_{X_n},\E\ot L^n)}. We have then
\begin{equation}\label{equ25}\Delta(\E) \ = \ 
\ov{\sigma}_{\E,X_{n+1}}\nabla_0(\E) \ ,\end{equation}
where \ \m{\nabla_0(\E)\in\Ext^1_{\ko_{X_n}}(\E,\E\ot\Omega_{X_n})} \ is the 
{\em fundamental class} of $\E$. 
\end{subsub}

\sepsubsub

Let $\E$ be a vector bundle on \m{X_n}. Then we have \ 
\m{H^0(\ko_{X_n})\subset\End(\E)}. We say that $\E$ is {\em simple} if \ 
\m{\End(\E)=H^0(\ko_{X_n})}.

Let \m{X_n} be a primitive multiple scheme of support $X$ and multiplicity $n$, 
that can be extended to a primitive multiple scheme \m{X_{n+1}} of multiplicity 
\m{n+1}. Suppose that we have a nonempty open set \m{\chi_n} of isomorphism 
classes of simple vector bundles on \m{X_n} such that there is a fine moduli 
space $M_n$ for $\chi_n$. Then there exists a closed subvariety \ \m{N_n\subset 
M_n} \ such that for every closed point \m{t\in N_n}, the corresponding bundle 
can be extended to \m{X_{n+1}}, and that for every scheme $T$ and every vector 
bundle $\E$ on \m{X_{n+1}\times T} such that for every closed point \m{t\in T}, 
\m{\E_{t|X_n}\in\chi_n}, then if \ \m{\kf=\E_{|X_n\times T}} \ and \ 
\m{f_\kf:T\to M_n} \ is the associated morphism, we have \ \m{f_\kf(T)\subset 
N_n}.

Let \m{\chi_{n+1}} be the set of isomorphism classes of vector bundles $E$ on 
\m{X_{n+1}} such that \ \m{E_{|X_n}\in\chi_n}. We will show that with suitable 
hypotheses, there is also a fine moduli space \m{M_{n+1}} for \m{\chi_{n+1}}, 
and that the restriction morphism \ \m{M_{n+1}\to N_n} \ is an affine bundle.

\sepsubsub

\begin{subsub}\label{hyp_i}
 Construction of fine moduli spaces -- \rm Let \m{\chi_n} be an 
open set of isomorphism classes of vector bundles on \m{X_n}. Suppose that 
there is a fine moduli space for \m{\chi_n}, defined by a smooth irreducible 
variety \m{M_n}, and that the subvariety \m{N_n\subset M_n} of 
bundles that can be extended to \m{X_{n+1}} is smooth (this variety is defined 
in \ref{extens}). We suppose that for every \m{\ke\in\chi_n}, if \m{E=\ke_{|X}},
\begin{enumerate}
\item[(i)] $\ke$ and $E$ are simple,
\item[(ii)] $\dim(\Ext^1_{\ko_X}(E,E\ot L^n))$ and 
$\dim(\Ext^1_{\ko_{X_n}}(\ke,\ke\ot\Omega_{X_n}))$ are independent of $\ke$,
\item[(iii)] $\dim(\Ext^2_{\ko_X}(E,E\ot L^n))$ is independent of $\ke$.
\item[(iv)] $\Delta(\ke)=0$ .
\item[(v)] every vector bundle $\kf$ on \m{X_{n+1}} such that \ 
\m{\kf_{|X_n}\simeq\ke} \ is simple. 
\end{enumerate}
\end{subsub}
Let \m{\chi_{n+1}} be the set of isomorphism classes of vector bundles $\E$ on 
\m{X_{n+1}} such that \ \m{\E_{|X_n}\in\chi_n}. We have then (cf. \ref{ext_F})

\sepprop

\begin{subsub}\label{theo_0_1}{\bf Theorem: } The set \m{\chi_{n+1}} is open, 
and there is a fine moduli space \m{M_{n+1}} for \m{\chi_{n+1}}, which is 
an affine bundle over \m{N_n}. 
\end{subsub}

\sepprop

We have an exact sequence \ \m{0\to L^n\to\ko_{X_{n+1}}\to\ko_{X_n}\to 0}, 
whence a connecting homomorphism \ \m{\delta_0:H^0(X_n,\ko_{X_n})\to 
H^1(X,L^n)} . Let \m{\ke\in N_n}, \m{E=\ke_{|X}}. Since \Nligne \m{E^*\ot E\ot 
L^n\simeq L^n\oplus(\AAd(E)\ot L^n)} (where \m{\AAd(E)} is the subbundle of 
\m{E^*\ot E} of 
trace zero endomorphisms), we can see \m{\imm(\delta^0)} as a 
subspace of \m{\Ext^1_{\ko_X}(E,E\ot L^n)}. Let $\bf A$ be the vector bundle on 
\m{N_n} associated to the affine bundle \m{M_{n+1}\to N_n}. Then we have a 
canonical isomorphism
\[{\bf A}_\ke \ \simeq \ \Ext^1_{\ko_X}(E,E\ot L^n)/\imm(\delta_0) \ . \]
Suppose that for some open subset \m{S\in M_n}, there is a universal vector 
bundle $\F$ on \m{X_n\times S}. Let \m{S'} be the inverse image of $S$ in 
\m{M_{n+1}}. We have then an obstruction class \ \m{\Delta(\F)\in H^2(X_n\times 
S,\HHom(\F,E\ot p_X^*(L^n)))} such that \m{\Delta(\F)=0} if and only if $\F$ 
can be extended to a vector bundle on \m{X_{n+1}\times S}. Then we have
\[{\bf A}_{|S} \ = \ R^1p_{S*}(\HHom(\F,E\ot p_X^*(L^n))/(\ko_S\ot 
\imm(\delta_0)) \ . \]
Suppose that \ \m{p_{S*}(\HHom(\F,E\ot p_X^*(L^n)))=0} \ and \m{\delta^0=0}. We 
have an injective canonical map
\[\lambda:H^1(S,R^1p_{S*}(\HHom(\F,E\ot p_X^*(L^n))))\lra H^2(X_n\times 
S,\HHom(\F,E\ot p_X^*(L^n)))\]
coming from the Leray spectral sequence.
Since \ \m{\Delta(\F_s)=0} \ for every \m{s\in S}, we will see that \ 
\[\Delta(\F) \ \in \ H^1(S,R^1p_{S*}(\HHom(\F,E\ot p_X^*(L^n)))) \ = \ H^1({\bf 
A}_{|S}) \ . \]
Suppose that \m{\eta\in H^1({\bf A}_{|S})} corresponds to the affine bundle \ 
\m{S'\to S}. Then (cf. \ref{theo3})

\sepprop

\begin{subsub}\label{theo_0_10}{\bf Theorem: } We have \ 
\m{\C.\Delta(\F)=\C\eta}.
\end{subsub}

\sepprop

In other words there is a link between the possibility to extend a whole family 
in higher multiplicity and the structure of the affine bundle. The affine 
bundle is actually a vector bundle if and only if the family can be extended in 
multiplicity \m{n+1}.

The hypothesis (v) is satisfied if $L$ is a non trivial ideal sheaf and if the 
restrictions to $X$ of the bundles of \m{\chi_n} are simple (Proposition 
\ref{endo}).

The assumption on \m{N_n} is satisfied in the the following cases:
\begin{enumerate}
\item[--] The case of Picard groups (cf. \ref{intro2}), i.e. $\chi_n$ consists 
of line bundles.
\item[--] When $N_n=M_n$, in particular if $X$ is a curve, because in this 
case, all the obstructions to extend vector bundles on \m{X_n} to \m{X_{n+1}} 
vanish.
\end{enumerate}
\end{sub}

\sepsub

\Ssect{Picard groups}{intro2}

Let ${\bf P}$ be an irreducible component of \m{\Pic(X)} and \m{{\bf P}_0} the 
component that contains \m{\ko_X}. Let \m{\kl_{n,{\bf P}}} be the set of line 
bundles $\L$ on \m{X_n} such that \ \m{\L_{|X}\in{\bf P}}. 
We prove by induction on $n$ that if \m{\kl_{n,{\bf P}}} is not empty, then 
there is a fine moduli space \m{\Pic^{\bf P}(X_n)} for \m{\kl_{n,{\bf P}}}. 

Suppose that some \m{\kb\in\kl_{n,\bf P}} can be extended to \m{X_{n+1}}. Then
the product with $\kb$ defines an isomorphism \ \m{\eta_\kb:\Pic^{{\bf 
P}_0}(X_n)\to\Pic^{\bf P}(X_n)}. The 
set of line bundles in \m{\kl_{n,\bf P}} that can be extended to \m{X_{n+1}} is 
a smooth closed subvariety
\[\Gamma^{\bf P}(X_{n+1}) \ \subset \ \Pic^{\bf P}(X_n) \ .\]
We have \ \m{\Gamma^{\bf P}(X_{n+1})=\eta_\kb\big(\Gamma^{{\bf 
P}_0}(X_{n+1})\big)}, and \m{\Gamma^{{\bf P}_0}(X_{n+1})} is a subgroup of
\m{\Pic^{{\bf P}_0}(X_n)}: it is the kernel of the morphism of groups
\[\Delta:\Pic^{{\bf P}_0}(X_n)\to H^2(X,L^n) \ , \]
(the {\em obstruction morphism}). The case \m{n=1} is simpler: since the 
fundamental class of a line bundle on $X$ is a discrete invariant, we have 
either \ \m{\Gamma^{\bf P}(X_{2})=\emptyset} \ or \ \m{\Gamma^{\bf 
P}(X_{2})=\Pic^{\bf P}(X)}. We have (cf. \ref{theo8})

\sepprop

\begin{subsub}\label{theo_0_2}{\bf Theorem: } Suppose that \m{\kl_{n+1,{\bf 
P}}} is nonempty. Then \m{\Pic^{\bf P}(X_{n+1})} is an affine 
bundle over \m{\Gamma^{\bf P}(X_{n+1})} with associated vector bundle \ 
\m{\ko_{\Gamma^{\bf P}(X_{n+1})}\ot\big(H^1(X,L^n)/\imm(\delta^0)\big)} .
\end{subsub}
\end{sub}

\sepsub

\Ssect{The case of curves}{Cur}

Suppose that $X$ is a curve. In this case the obstructions to extend a 
primitive multiple curve or a vector bundle to higher multiplicity vanish (from 
\cite{dr10} they lie in the \m{H^2} of some vector bundles on $X$).

Suppose that \m{\deg(L)<0}. Let $r$, $d$ be integers such 
that \m{r>0} and $r$, $d$ are coprime. We take for \m{M_n} the set of 
isomorphism classes of vector bundles $\ke$ on \m{X_n} such that \m{\ke_{|X}} 
is stable of rank $r$ and degree $d$. The hypotheses of \ref{hyp_i} are 
satisfied, and we obtain by induction on $n$

\sepprop

\begin{subsub}\label{theo_0_3}{\bf Theorem: } There is a fine moduli space 
\m{M_{X_n}(r,d)} for \m{M_n}, which is a smooth irreducible variety. It is an 
affine bundle over \m{M_{X_{n-1}}(r,d)}.
\end{subsub}

Of course for \m{n=1}, \m{M_{X_n}(r,d)} is the well known moduli space of 
stable vector bundles of rank $r$ and degree $d$ on $X$. From \cite{dr5}, the 
vector bundles of \m{M_{X_n}(r,d)}, \m{n\geq 2}, are also stable. Hence 
\m{M_{X_n}(r,d)} is an open subset of a moduli space of stable sheaves on 
\m{X_n}. 

Let $\E$ be a Poincar\'e bundle on \ \m{X\times M_X(r,d)}. Let \ \m{p_X:X\times 
M_X(r,d)\to X},\Nligne \m{p_M:X\times M_X(r,d)\to M_X(r,d)} \ be the 
projections. Then if \m{n\geq 2},
\[{\bf A}_{n-1} \ = \ R^1p_{M*}(\E\ot\E^*\ot p_X^*(L^{n-1}))\]
is a vector bundle on \m{M_X(r,d)}. Let \ \m{p_{n-1}:M_{X_{n-1}}(r,d)\to 
M_X(r,d)} \ be the restriction morphism. Then \m{p_{n-1}^*({\bf A}_{n-1})} is 
the vector bundle on \m{M_{X_{n-1}}(r,d)} associated to the affine bundle \
\m{M_{X_n}(r,d)\to M_{X_{n-1}}(r,d)}.

\sepsubsub

\begin{subsub} The case of Picard groups -- \rm We don't require here that 
\m{\deg(L)<0}. The irreducible components of \m{\Pic(X)} are the smooth 
projective varieties \ \m{{\bf P}_d=\Pic^d(X)} \ of line bundles of degree $d$.
\end{subsub}
Then if follows from \ref{intro2} that

\sepprop

\begin{subsub}\label{theo_0_4}{\bf Theorem:} The irreducible components of 
\m{\Pic(X_n)} are the \m{\Pic^{{\bf P}_d}(X_n)}, \m{d\in\Z}, and\break 
\m{\Pic^{{\bf P}_d}(X_{n+1})} is an affine bundle over \m{\Pic^{{\bf 
P}_d}(X_n)}, with associated vector bundle \Nligne \m{\ko_{\Pic^{{\bf 
P}_d}(X_n)}\ot(H^1(X,L^n)/\imm(\delta_0))} .
\end{subsub}

\sepprop

These affine bundles need not to be banal: (cf. \ref{theo4})

\sepprop

\begin{subsub}{\bf Theorem :}\label{theo_0_5} 
\m{\Pic^{d}(X_2)} is not a vector bundle over \m{\Pic^d(X)} if $Y$ is not 
trivial and either $X$ is not hyperelliptic and \ \m{\deg(L)\leq 2-2g}, or 
\m{L=\omega_X}.
\end{subsub}

\sepprop

Using some parts of the proof of Theorem \ref{theo_0_5} we obtain

\sepprop

\begin{subsub}{\bf Theorem :}\label{theo_0_15} Suppose that $C$ is not 
hyperelliptic, \m{X_2} not trivial and \ \m{\deg(L)\leq 2-2g}.
Then \m{{\bf M}_{X_2}(r,-1)} is not a vector bundle over \m{M_X(r,-1)}.
\end{subsub}

\end{sub}

\sepsub

\Ssect{Primitive double schemes with trivial canonical line bundle}{triv_can}

Since \m{X_n} is a locally complete intersection, it has a {\em dualizing 
sheaf} \m{\omega_{X_n}} which is a line bundle on \m{X_n}. We call also 
\m{\omega_{X_n}} the canonical line bundle of \m{X_n}. We have
\[\omega_{X_n|X} \ \simeq \ \omega_X\ot L^{1-n} \ . \]
In particular we have \ \m{\omega_{X_2|X}\simeq\omega_X\ot L^{-1}}. We have an 
exact sequence
\xmat{0\ar[r] & L\ar[r] & \ko_{X_2}\ar[r]^-\zeta & \ko_X\ar[r] & 0 \ .}
We will prove (cf. \ref{theo7})

\sepprop

\begin{subsub}\label{theo_0_6}{\bf Theorem: } The map \ 
\m{H^1(\zeta):H^1(X_2,\ko_{X_2})\to H^1(X,\ko_X)} \ is surjective.
\end{subsub}

\sepprop

which will imply (cf. \ref{coro4})

\sepprop

\begin{subsub}\label{coro_04}{\bf Corollary:} We have \ 
\m{\omega_{X_2}\simeq\ko_{X_2}} \ if and only if \ \m{\omega_X\simeq L} .
\end{subsub}

\sepprop

If $X$ is a surface, \m{h^1(X_2,\ko_{X_2})=0} and \ 
\m{\omega_{X_2}\simeq\ko_{X_2}}, \m{X_2} is called a 
{\em K3-carpet} (cf. \cite{b_m_r}, \cite{ga-go-pu}). 

Suppose that \m{L\simeq\omega_X}. Then according to the exact sequence \ 
\m{0\to\omega_X\to\ko_{X_2}\to\ko_X\to 0}, we have \m{h^1(X_2,\ko_{X_2})=0} if 
and only if \m{h^1(X,\ko_X)=0} (cf. also \cite{ga-go-pu}, Proposition 1.6).

\end{sub}

\sepsub

\Ssect{Examples}{X_exa}

Let $C$ and $D$ be smooth projective irreducible curves, \m{X=C\times D} and
\ \m{\pi_C:X\to C}, \m{\pi_D:X\to D} \ the projections. Let \m{g_C} (resp. 
\m{g_D}) be the genus of $C$ (resp. $D$). Suppose that \m{g_C\geq 2}, 
\m{g_D\geq 2}. Let \m{L_C} (resp. \m{L_D}) be a line bundle on $C$ (resp. $D$) 
and \ \m{L=\pi_C^*(L_C)\ot\pi_D^*(L_D)}. The non trivial double schemes \m{X_2} 
with associated line bundle $L$ are parametrized by \m{\P(H^1(X,T_X\ot L))}, 
where \m{T_X} is the tangent bundle of $X$ (cf. \cite{ba_ei} and \cite{dr10}). 
We have
\[H^1(X,T_X\ot L) \ = \ \big(H^0(\omega_C^*\ot L_C)\ot 
H^1(L_D)\big)\oplus
\big(H^1(\omega_C^*\ot L_C)\ot H^0(L_D)\big)\oplus\]
\[ \qquad\qquad\qquad\qquad\quad
\big(H^0(L_C)\ot H^1(\omega_D^*\ot L_D)\big)\oplus
\big(H^1(L_C)\ot H^0(\omega_D^*\ot L_D)\big) \ , \]
and using this decomposition, the double schemes $X_2$ are defined by 
quadruplets \m{(\eta_1,\eta_2,\eta_3,\eta_4)}.

We denote by \m{\phi_{L_C}}, \m{\phi_{L_D}} the canonical maps \ 
\m{H^0(\omega_C^*\ot L_C)\ot H^1(\omega_C)\to H^1(L_C)} \ and \Nligne 
\m{H^0(\omega_D^*\ot L_D)\ot H^1(\omega_D)\to H^1(L_D)} \ respectively. 

Let \m{M_C} be a line bundle on $C$ (resp. $D$), and \ 
\m{M=\pi_C^*(M_C)\ot\pi_D^*(M_D)} . Using $(\ref{equ25})$ we will see that
\[\Delta(M) \ = \ (\phi_{L_C}\ot 
I_{H^1(L_D)})(\eta_1\ot\nabla_0(M_C))+(I_{H^1(L_C)}\ot
\phi_{L_D})(\eta_4\ot\nabla_0(M_D)) \ . \]
The canonical class of a line bundle on a smooth projective curve is in fact an 
integer. For every line bundle $\kl$ on $C$ we have \ 
\m{\nabla_0(\kl)=\deg(\kl)c}, where \ \m{c=\nabla_0(\ko_C(P))}, for any \m{P\in 
C}.

For example let $\Theta$ be a theta characteristic on $C$ such that 
\m{h^0(C,\Theta)>0}. We take \ \m{L_C=\Theta}, \m{L_D=\omega_D}. We have 
\m{\eta_1=0}. Suppose that \ \m{\eta_4=0}. Such non trivial double schemes are 
parametrized by \ \m{\P\big((H^1(\Theta^{-1})\ot H^0(\omega_D))\oplus
(H^0(\Theta)\ot H^1(\ko_D))\big)} (pairs \m{(\eta_2,\eta_3)}). Suppose that \ 
\m{\eta_3\not=0}. Then every line bundle on $X$ can be extended to \m{X_2}, 
which is projective, and (using Theorem \ref{theo_0_7} below) there exists an 
extension of \m{X_2} to a primitive multiple scheme \m{X_3} of multiplicity 3, 
and \m{X_3} is projective.

\sepsubsub

\begin{subsub} The case \m{L=\omega_X} -- \rm
This happens if and only if \m{L_C=\omega_C} and \m{L_D=\omega_D}. Let
\ \m{\eta=(\eta_1,\eta_2,\eta_3,\eta_4)\in H^1(X,T_X\ot\omega_X)}.
Here we have
\[\eta_1\in\C \ , \quad \eta_2\in H^1(C,\ko_C)\ot H^0(D,\omega_D) \ , \quad
\eta_3\in H^0(C,\omega_C)\ot H^1(D,\ko_D) \ , \quad \eta_4\in\C \ . \]
Let \m{M_C} (resp. \m{M_D}) be a line bundle on $C$ (resp. $D$). If \ 
\m{M=\pi_C^*(M_C)\ot\pi_D^*(M_D)}, then
\[\Delta(M) \ = \ \eta_1\nabla_0(M_C)+\eta_4\nabla_0(M_D) \ , \]
hence $M$ can be extended to a line bundle on \m{X_2} if and only if \ 
\m{\eta_1\deg(M_C)+\eta_4\deg(M_D)=0}.
\end{subsub}
It follows that
{\em the scheme \m{X_2} is projective if and only if \ \m{\eta_1=\eta_4=0} \ or 
\m{\eta_1\eta_4<0} \ and \m{\dsp\frac{\eta_1}{\eta_4}} is rational.}
We will see that \m{X_2} can be extended to a primitive multiple scheme of 
multiplicity 3 if and only \m{\omega_X} can be extended to a line bundle on 
\m{X_2}. This is the case if and only \ \m{(g_C-1)\eta_1+(g_D-1)\eta_4=0}, and 
then \m{X_2} is projective.

We will see also that there may exist components $\bf P$ of \m{\Pic(X)} such 
that the affine bundle \m{\Pic^{\bf P}(X_2)} is not empty, and not a vector 
bundle over $\bf P$.

\end{sub}

\sepsub

\Ssect{Improvements of results of \cite{dr10}}{Impr}

\begin{subsub}\label{impr01} Extensions of \m{X_n} to multiplicity \m{n+1} -- 
\rm If \m{X_n} can be extended to a primitive multiple scheme \m{X_{n+1}} of 
multiplicity \m{n+1}, then \m{\ki_{X,X_n}}, which is a line bundle on 
\m{X_{n-1}}, can be extended to a line bundle on \m{X_n}, namely 
\m{\ki_{X,X_{n+1}}}. Conversely, given an extension of \m{\ki_{X,X_n}} to a 
line bundle $\L$ on \m{X_n}, an invariant \ \m{\Delta''_n(X_n,\L)\in 
H^2(X,T_X\ot L^n)} \ is defined in \cite{dr10}, such that \m{X_n} can be 
extended to \m{X_{n+1}} with \ \m{\ki_{X,X_{n+1}}\simeq\L} \ if and only if \ 
\m{\Delta''_n(X_n,\L)=0}.

By \ref{111}, another extension $\L'$ of \m{\ki_{X,X_n}} to \m{X_{n+1}} comes 
from some \ \m{\eta\in H^1(X,L^{n-1})}. Let \ \m{\boldsymbol{\zeta}\in 
H^1(X,T_X\ot L)} \ correspond to the exact sequence \ \m{0\to 
L\to\Omega_{X_2|X}\to\Omega_X\to 0}.
We have a canonical product
\[H^1(X,T_X\ot L)\ot H^1(X,L^{n-1})\lra H^2(X,T_X\ot L^n)\]
and (cf. \ref{theo6})
\end{subsub}

\sepprop

\begin{subsub}\label{theo_0_7}{\bf Theorem: } We have \ 
\m{\Delta''_n(X_n,\L')-\Delta''_n(X_n,\L)=\boldsymbol{\zeta}\eta} .
\end{subsub}

\sepprop

This implies that if the product is a surjective map, then there exists an 
extension of \m{X_n} in multiplicity \m{n+1}.

\sepsubsub

\begin{subsub}\label{impr02} Extensions of a line bundle on \m{X_n} to 
\m{X_{n+1}} -- \rm
Let \ \m{\L=\ki_{X,X_{n+1}|X_2}}, which is a line bundle on \m{X_2}. Let 
\m{\delta^1_{\L^{n-1}}} be the connecting morphism \ \m{H^1(X,L^{n-1})\to 
H^2(X,L^n)} \ coming from the exact sequence \ \m{0\to L^n\to\L^{n-1}\to 
L^{n-1}\to 0}.
\end{subsub}
Suppose that \m{X_n} can be extended to \m{X_{n+1}}, and let $\D$ be a line 
bundle on \m{X_n}. By \ref{111}, $\D$ can be extended to a line bundle on 
\m{X_{n+1}} if and only if \ \m{\Delta(\D)=0}. Let \m{\D'\in\Pic(X_n)} be such 
that \ \m{\D'_{|X_{n-1}}\simeq\D_{|X_{n-1}}}. If \ \m{\E=\D_{|X_{n-1}}}, 
\m{E=\D_{|X}}, we have exact sequences
\[0\to E\ot L^{n-1}\to\D\to\E\to 0 \ , \qquad 0\to E\ot L^{n-1}\to\D'\to\E\to 0 
\ , \]
corresponding to \ \m{\eta,\eta'\in\Ext^1_{\ko_{X_n}}(\E,E\ot 
L^{n-1})=H^1(X,L^{n-1})} \ respectively. Then (cf. \ref{theo5})

\sepprop

\begin{subsub}\label{theo_0_8}{\bf Theorem:} We have \
\m{\Delta(\D')-\Delta(\D)=n\delta^1_{\L^{n-1}}(\eta'-\eta)} .
\end{subsub}

\sepprop

Theorem \ref{theo_0_8} implies for Picard groups (cf. \ref{intro2}, \ref{coro5})

\sepprop

\begin{subsub}\label{prop_0_1}{\bf Proposition:} Let $Z$ be the image of he 
restriction morphism \ \m{\Gamma^{\bf P}(X_{n+1})\to\Gamma^{\bf P}(X_n)}. Then 
$Z$ is smooth and \m{\Gamma^{\bf P}(X_{n+1})} is an affine bundle over $Z$ with 
associated vector bundle \ \m{\ko_Z\ot\ker(\delta^1_{\L^{n-1}})} .
\end{subsub}

\end{sub}

\sepsub

\Ssect{Outline of the paper}{Outl}

In Chapter 2, we give several ways to use \v Cech cohomology, and recall some 
basic definitions about fine moduli spaces and affine bundles.

In the very technical Chapter 3, we give a description of Leray's spectral 
sequence and the exact sequence of Ext's in terms on \v Cech cohomology (which 
is our main tool). 

The Chapter 4 is devoted to the definitions and properties of primitive 
multiple schemes, with the improvement \ref{impr01} of \cite{dr10}. It contains 
also supplementary results on vector bundles on primitive multiple schemes, 
that are used in the following chapters.

In Chapter 5 the definition and properties of the fundamental class of a vector 
bundles are recalled. We also give here the results of \ref{triv_can}.

Chapter 6 contains the proofs of the main results (the construction and 
properties of fine moduli spaces of vector bundles on primitive multiple 
schemes).

In Chapter 7 we treat the special case of Picard groups (cf. \ref{intro2}), and 
prove the improvement \ref{impr02} of \cite{dr10}.

In Chapter 8 we give the examples \ref{X_exa}, when $X$ is a product of curves.

Chapter 9 is devoted to moduli spaces of vector bundles on primitive 
multiple curves.

\end{sub}

\sepsub

{\bf Notations and terminology:} -- An {\em algebraic variety} is a 
quasi-projective scheme over $\C$. A {\em vector bundle} on a scheme is an 
algebraic vector bundle.

-- If $X$ is a scheme and \m{P\in X} is a closed point, we denote 
by \m{m_{X,P}} (or \m{m_P}) the maximal ideal of $P$ in \m{\ko_{X,P}}.

-- If $X$ is a scheme and \m{Z\subset X} is a closed subscheme, 
\m{\ki_Z} (or \m{\ki_{Z,X}}) denotes the ideal sheaf of $Z$ in $X$.

-- If \m{V} is a finite dimensional complex vector space, \m{\P(V)} denotes the 
projective space of lines in $V$, and we use a similar notation for projective 
bundles.

-- Let $X$ be a set and \m{(X_i)_{i\in I}} a family of subsets of $X$. If 
\m{i,j\in I}, \m{X_{ij}} denotes the intersection \m{X_i\cap X_j}, and 
similarly for more indices.

-- If $X$ and $Y$ are two schemes and \m{f:X\to Y} is a morphism, \m{{\bf 
d}(f):f^*(\Omega_Y)\to\Omega_X} \ is the corresponding canonical morphism.

-- If $R$ is a commutative ring, $n$ a positive integer, \m{{\bf a}\in 
R[t]/(t^n)} \ and $i$ an integer such that \m{0\leq i<n}, let \m{{\bf a}^{(i)}} 
denote the coefficient of \m{t^i} on $\bf a$, so that we have \ \m{{\bf 
a}=\sigg_{i=0}^{n-1}{\bf a}^{(i)}t^i}.

\sepsec

\section{Preliminaries}

\Ssect{\v Cech cohomology}{cech}

Let $X$ be a scheme over $\C$ and $E$, $V$ vector bundles on $X$ of rank 
$r$. Let \m{(U_i)_{i\in I}} be an open cover of $X$ such that 
we have isomorphisms:
\xmat{\alpha_i:E_{|U_i}\ar[r]^-\simeq & V_{|U_i} .}
Let
\xmat{\alpha_{ij}=\alpha_i\circ\alpha_j^{-1}:V_{|U_{ij}}\ar[r]^-\simeq &
V_{|U_{ij}} ,}
so that we have the relation \ \m{\alpha_{ij}\alpha_{jk}=\alpha_{ik}} \ on 
\m{U_{ijk}}.

\sepprop

\begin{subsub}\label{cech1}\rm 
Let $n$ be a positive integer, and for every sequence \m{(i_0,\ldots,i_n)} of 
distinct elements of $I$, \m{\sigma_{i_0\cdots i_n}\in H^0(U_{i_0\cdots 
i_n},V)}. Let
\[\theta_{i_0\cdots i_n} \ = \ \alpha_{i_0}^{-1}\sigma_{i_0\cdots i_n} \ \in \
H^0(U_{i_0\cdots i_n},E) \ . \]
The family \m{(\theta_{i_0\cdots i_n})} represents an element of \m{H^n(X,E)} 
if the cocycle relations are satisfied: for every sequence  
\m{(i_0,\ldots,i_{n+1})} of distinct elements of $I$, 
\[\sigg_{k=0}^n(-1)^k\theta_{i_0\cdots\widehat{i_k}\cdots i_{n+1}} \ = \ 0 \ , 
\]
which is equivalent to
\[\alpha_{i_0i_1}\sigma_{i_1\cdots i_{n+1}}+\sigg_{k=1}^{n+1}
(-1)^k\sigma_{i_0\cdots\widehat{i_k}\cdots i_{n+1}} \ = \ 0 \ . \]
For \m{n=1}, this gives that elements of \m{H^1(X,E)} are represented by 
families \m{(\sigma_{ij})},\Nligne \m{\sigma_{ij}\in H^0(U_{ij},V)}, 
such that
\[\alpha_{ij}\sigma_{jk}+\sigma_{ij}-\sigma_{ik} \ = \ 0\]
(instead of \m{(\alpha_i^{-1}\sigma_{ij})} with the usual cocycle relations).
In \v Cech cohomology, it is generally assumed that 
\m{\theta_{ji}=-\theta_{ij}}. This implies that \ 
\m{\sigma_{ji}=-\alpha_{ji}\sigma_{ij}}.
\end{subsub}

\sepprop

\begin{subsub}\label{cech2}\rm Let $F$ be a vector bundle on $X$.
Similarly, an element of \m{H^n(X,E\ot F)} is 
represented by a family \m{(\mu_{i_0\cdots i_n})}, with \ \m{\mu_{i_0\cdots 
i_n}\in H^0(U_{i_0\cdots i_n},V\ot F)}, satisfying the relations
\[(\alpha_{i_0i_1}\ot I_F)(\mu_{i_1\cdots i_{n+1}})+\sigg_{k=1}^{n+1}
(-1)^k\mu_{i_0\cdots\widehat{i_k}\cdots i_{n+1}} \ = \ 0 \ . \]
The corresponding element of \m{H^0(U_{i_0\cdots i_n},E\ot F)} is \ 
\m{\theta_{i_0\cdots i_n}=(\alpha_{i_0}\ot I_F)^{-1}(\mu_{i_0\cdots i_n})}.
\end{subsub}

\sepprop

\begin{subsub}\label{cech2f} Coboundaries -- \rm In the situation of 
\ref{cech2}, the family \m{(\theta_{i_0\cdots i_n})} represents $0$ in 
\m{H^n(X,E\ot F)} if and only if there exists a family \m{(\tau_{j_0\cdots 
j_{n-1}})}, with \ \m{\tau_{j_0\cdots j_{n-1}}\in H^0(U_{j_0\cdots 
j_{n-1}},E\ot F)}, such that
\[\theta_{i_0\cdots i_n} \ = \ 
\sigg_{k=0}^n(-1)^k\tau_{i_0\cdots\widehat{i_k}\cdots i_n} \ . \]
Similarly, the family \m{(\mu_{i_0\cdots i_n})} represents $0$ if and only 
there exists a family  \m{(\nu_{j_0\cdots j_{n-1}})}, with \ 
\m{\nu_{j_0\cdots j_{n-1}}\in H^0(U_{j_0\cdots j_{n-1}},V\ot F)}, such that
\[\mu_{i_0\cdots i_n} \ = \ 
(\alpha_{i_0i_1}\ot I_F)(\nu_{i_1\cdots i_n})+\sigg_{k=1}^{n}
(-1)^k\nu_{i_0\cdots\widehat{i_k}\cdots i_{n}} \ . \]
\end{subsub}

\sepprop

\begin{subsub}\label{cech2c}Representation of morphisms -- \rm Suppose that we 
have also local trivializations of $F$:
\xmat{\beta_i:F_{|U_i}\ar[r]^-\simeq & \ko_{U_i}\ot\C^s .}
We have then local trivializations of \m{\HHom(E,F)}
\[\xymatrix@R=5pt{\Delta_i:\HHom(E,F)_{|U_i}\ar[r]^-\simeq & 
\ko_{U_i}\ot L(\C^r,\C^s)\\ \phi\fmaps[r] & \beta_i\circ\phi\circ\alpha_i^{-1} 
}\]
such that, for every \ \m{\lambda\in H^0(\ko_{U_{ij}}\ot L(\C^r,\C^s))}, we have
\[\Delta_{ij}(\lambda) \ = \ \Delta_i\Delta_j^{-1}(\lambda) \ = \
\beta_{ij}\lambda\alpha_{ij}^{-1} \ . \]
\end{subsub}

\sepprop

\begin{subsub}\label{cech2e} Products -- \rm Let \ \m{\chi_{E,F}:H^1(X,E)\ot 
H^1(X,F)\to H^2(X,E\ot F)} \ be the canonical map. Let \ \m{\epsilon\in 
H^1(X,E)}, \m{\phi\in H^1(X,F)}, represented by cocycles \m{(\epsilon_{ij})}, 
\m{(\phi_{ij})}, \m{\epsilon_{ij}\in H^0(U_{ij},E)}, \m{\phi_{ij}\in 
H^0(U_{ij},F)}. Then \m{\chi_{E,F}(\epsilon\ot\phi)} is represented by the 
cocycle \m{(\epsilon_{ij}\ot\phi_{jk})}.

With respect to the canonical isomorphism \ \m{F\ot E\simeq E\ot F} \ we have
\Nligne \m{\chi_{F,E}(\phi\ot\epsilon)=-\chi_{E,F}(\epsilon\ot\phi)} .

Suppose that we have vector bundles $V$, $W$ on $X$ with isomorphisms
\[\alpha_i:E_{|U_i}\lra V_{|U_i} \ , \qquad \beta_i:F_{|U_i}\lra W_{|U_i} \ .  
\]
Let \ \m{e_{ij}=\alpha_i\epsilon_{ij}}, \m{f_{ij}=\beta_i\phi_{ij}} . We have 
then
\[\epsilon_{ij}\ot\phi_{jk} \ = \ (\alpha_i^{-1}\ot\beta_j^{-1})(e_{ij}\ot 
f_{jk}) \ , \]
and
\[(\alpha_i\ot\beta_i)(\epsilon_{ij}\ot\phi_{jk}) \ = \ (I\ot\beta_{ij})
(e_{ij}\ot f_{jk}) \ . \]
It follows that the family \m{(e_{ij}\ot\beta_{ij}f_{jk})} represents 
\m{\chi_{E,F}(\epsilon\ot\phi)} in the sense of \ref{cech1}, with respect to 
the local isomorphisms \ \m{\alpha_i\ot\beta_i:(E\ot F)_{|U_i}\lra (V\ot 
W)_{|U_i}}.
\end{subsub}

\sepsubsub

\begin{subsub}\label{cech2b}Construction of vector bundles via local 
isomorphisms -- \rm Let $Z$ be a scheme over $\C$, \m{(Z_i)_{i\in I}} an open 
cover of $Z$, and for every \m{i\in I}, a scheme \m{U_i}, with an isomorphism 
\m{\delta_i:Z_i\to U_i}. Let \ 
\m{\delta_{ij}=\delta_j\delta_i^{-1}:\delta_i(Z_{ij})\to\delta_j(Z_{ij})}, 
which is an isomorphism.
Then a vector bundle of $X$ can be constructed (in an obvious way) using vector 
bundles \m{E_i} on \m{U_i} and isomorphisms \ 
\[\Theta_{ij}:\delta_{ij}^*(E_{j|\delta_j(Z_{ij})})\lra 
E_{i|{\delta_i(Z_{ij})}}\]
such that \ \m{\Theta_{ij}\circ\delta_{ij}^*(\Theta_{jk})=\Theta_{ik}} \ on 
\m{\delta_i(Z_{ijk})}.
\end{subsub}

\sepsubsub

\begin{subsub}\label{cech3} Representation of extensions -- \rm Let $E$, $F$ be 
coherent sheaves on $X$, and\Nligne \m{\sigma\in 
H^1(X,\HHom(F,E))\subset\Ext^1_{\ko_X}(F,E)}. Let
\[0\lra E\lra\ke\lra F\lra 0\]
be the corresponding exact sequence. Suppose that $\sigma$ is represented by a 
cocycle \m{(\sigma_{ij})}, \m{\sigma_{ij}:F_{|U_{ij}}\to E_{|U_{ij}}}. Then 
$\ke$ can be constructed by gluing the sheaves \m{(E\oplus F)_{|U_i}} 
using the automorphisms of \m{(E\oplus F)_{|U_{ij}}} defined by the matrices
\ \m{\dsp\begin{pmatrix}I_E & \sigma_{ij}\\ 0 & I_F\end{pmatrix}}.
\end{subsub}

\end{sub}

\sepsub

\Ssect{Moduli spaces of sheaves}{Fine}

(cf. \cite{dr1c})

Let $X$ be a scheme over $\C$.

\sepprop

\begin{subsub} Families and sets of vector bundles -- \rm
Let $S$ be a scheme over $\C$.

A {\em family of vector bundles on $X$ parametrized by $S$} is a vector bundle 
on \m{X\times S}. If \m{f:T\to S} is a morphism of schemes and $\ke$ 
a family of vector bundles on $X$ parametrized by $S$, we note
\[f^\sharp(\ke) \ = \ (I_X\times f)^*(\ke) \ , \]
which is a family of vector bundles on $X$ parametrized by $T$.

Let $\chi$ be a nonempty set of isomorphism classes of vector bundles on $X$. 
A {\em family of vector bundles of $\chi$ parametrized by $\chi$} is a family 
of vector bundles $\ke$ on $X$ parametrized by $S$ such that for every closed 
point \m{s\in S}, \m{\ke_s} belongs to $\chi$. 

We say that $\chi$ is {\em open} if for every $S$ and every family of vector 
bundles $\ke$ on $X$ parametrized by $S$, if \m{s\in S} is a closed point such 
that \m{\ke_s\in\chi}, then there exists a neighbourhood $U$ of $s$ such that 
\m{\ke_u\in\chi} for every closed point \m{u\in U}. 
\end{subsub}

\sepprop

\begin{subsub}\label{Fine2} Fine moduli spaces -- \rm Let $\chi$ be a nonempty 
open set of isomorphism classes of vector bundles on $X$. A {\em fine moduli 
space for $\chi$} is the data of an integral scheme $M$ and of
\begin{enumerate}
\item[--] a bijection
\[\xymatrix@R=5pt{\phi:M^0\ar[r] & \chi\\ m\fmaps[r] & E_m
}\]
(where \m{M^0} denotes the set of closed points of $M$),
\item[--] An open cover \m{(M_i)_{i\in I}} of $M$, and for every \m{i\in I}, a 
vector bundle \m{\ke_i} on \m{X\times M_i} such that for every \m{m\in M_i}, 
\m{\ke_{i,m}\simeq E_m},
\end{enumerate}
such that for any scheme $S$, any family $\kf$ of vector bundles of $\chi$ 
parametrized by $S$, there exists a morphism \ \m{f_\kf:S\to M} such that: 
for every closed point \m{s\in S}, if \m{m=f_\kf(s)} and if \m{m\in M_i}, there 
exists an open neighbourhood $U$ of $s$ such that 
\ \m{f_\kf(U)\subset M_i} \ and \m{f^\sharp_{\kf|U}(\ke_i)\simeq\kf_{|X\times 
U}}. In particular we have \ \m{\kf_s\simeq E_m} .

If $\chi$ is an open set, then
for every closed point \m{m\in M}, if \m{m\in M_i}, then \m{\ke_i} is a 
semi-universal deformation of \m{\ke_{i,m}}, hence the Koda\"\i ra-Spencer map \
\m{T_m M_i\to\Ext^1_{\ko_X}(\ke_{i,m},\ke_{i,m})} \ is an isomorphism.

Fine moduli spaces are unique.
\end{subsub}
\end{sub}

\sepsub

\Ssect{Affine bundles}{aff}

Let \ \m{f:\ka\to S} \ be a morphism of schemes, and \m{r\geq 0} an integer. We 
say that $f$ (or $\ka$) is an {\em affine bundle} of rank $r$ over $S$ if there 
exists an open cover \m{(S_i)_{i\in I}} of $S$ such that for every \m{i\in I} 
there is an isomorphism \ \m{\tau_i:f^{-1}(S_i)\to S_i\times\C^r} \ over 
\m{S_i} such that for every distinct \m{i,j\in I}, \m{f_j\circ 
f_i^{-1}:S_{ij}\times\C^r\to S_{ij}\times\C^r} \ is of the form
\[(x,u)\longmapsto (x,A_{ij}(x)u+b_{ij}(x)) \ , \]
where \m{A_{ij}} is an \m{r\times r}-matrix of elements of \m{\ko_S(S_{ij})} 
and \m{b_{ij}} is a morphism from \m{S_{ij}} to \m{\C^r}. We have then the 
cocycle relations
\[A_{ij}A_{jk} \ = \ A_{ik} \ , \quad b_{ik} \ = \ A_{ij}b_{jk}+b_{ij} \ . \]
The first relation shows that the family \m{(A_{ij})} defines a vector bundle 
$\A$ on $S$, and the second that \m{(b_{ij})} defines \ \m{\lambda\in 
H^1(S,\A)} (according to \cite{dr10}, 2.1.1). 

The vector bundle $\A$ is uniquely defined, as well as \ 
\m{\eta(\ka)=\C\lambda\in\big(\P(H^1(S,\A))\cup\{0\}\big)/\Aut(\A)}.

We say that {\em $\ka$ is a vector bundle} if $f$ has a section. This is the 
case if and only \m{\lambda=0}, and then a section of $f$ induces an 
isomorphism \m{\ka\simeq\A} over $S$.

For every closed point \m{s\in S} there is a canonical action of the additive 
group \m{\A_s} on \m{\ka_s}
\[\xymatrix@R=5pt{\A_s\times\ka_s\ar[r] & \ka_s\\ (u,a)\fmaps[r] & a+u}\]
such that for every \m{a\in\ka_s}, $u\mapsto a+u$ is an isomorphism 
\m{\A_s\simeq\ka_s}.

\sepsubsub

\begin{subsub}\label{aff_quot} Quotients -- \rm Let \ \m{\B\subset\A} \ be a 
subbundle of $\A$, \m{q:\A\to\B/\A} \ the quotient morphism and \ 
\m{R:H^1(S,\A)\to H^1(S,\A/\B)} \ the induced morphism. Let \ \m{\lambda\in 
H^1(S,\ka)} \ and \ \m{\lambda''=R(\lambda)}. Let \ \m{f'':\ka''\to S} \ be the 
affine bundle corresponding to \m{\lambda''}. We will also use the notation \ 
\m{\ka''=\ka/\B}. Then there is an obvious 
$S$-morphism \ \m{p:\ka\to\ka''}, such that for every closed point \m{s\in S}, 
\m{a\in\ka_s} and \m{u\in\A_s}, we have \ \m{p_s(a+u)=p_s(a)+q_s(u)}. Moreover 
$p$ is an affine bundle with associated vector bundle \m{{f''}^*(\B)}.

\end{subsub}

\end{sub}

\sepsec

\section{Canonical exact sequences in \v Cech cohomology}\label{seb_cech}

The Ext spectral sequence allows to compute Ext groups in terms of cohomology 
groups and Leray's spectral sequence to compute the cohomology groups of a 
sheaf in terms of the cohomology groups of its direct images by a suitable 
morphism. We will describe simple examples of consequences of these spectral 
sequences in terms of \v Cech cohomology (this provides also elementary proofs 
of these consequences).

\sepsub

\Ssect{Leray's spectral sequence}{Ex_Ler}

Let \ \m{\pi:X\to S} \ be a proper flat morphism of sheaves, and $\ke$ a 
coherent sheaf on $X$ such that \ \m{\pi_*(\ke)=0}. Leray's spectral sequence 
implies that we have an exact sequence
\xmat{0\ar[r] & H^1(S,R^1\pi_*(\ke))\ar[r]^-\lambda & H^2(X,\ke)\ar[r]^-\eta & 
H^0(S,R^2\pi_*(\ke))\ar[r] & 0 \ . }
We will describe $\lambda$ and $\eta$ in terms of \v Cech cohomology.

Let \m{(S_m)_{m\in M}} (resp. \m{(X_i)_{i\in I}}) be an affine cover of $S$ 
(resp. $X$), such that for every \m{i\in I}, \m{m\in M}, \m{\pi^{-1}(S_m)\cap 
X_i} is affine.

\sepprop

\begin{subsub}\label{des_eta} Description of $\eta$ -- \rm
Let \m{\Sigma\in H^2(X,\ke)}, represented by the cocycle \m{(B_{ijk})_{i,j,k\in 
I}}, \Nligne \m{B_{ijk}\in H^0(X_{ijk},\ke)}. For every \m{m\in M}, 
\m{(B_{ijk|\pi^{-1}(S_m)\cap X_{ijk}})_{i,j,k\in I}} defines \ \m{\Sigma_m\in 
H^2(\pi^{-1}(S_m),\ke)}, and \m{(\Sigma_m)_{m\in M}} defines \ 
\m{\eta(\Sigma)\in H^0(S,R^2\pi_*(\ke))}.
\end{subsub}

\sepprop

\begin{subsub}\label{des_lam} Description of $\lambda$ -- \rm
Let \ \m{\theta\in H^1(S,R^1\pi_*(\ke))}, represented by the cocycle 
\m{(A_{mn})_{m,n\in M}}, \m{A_{mn}\in H^1(\pi^{-1}(S_{mn}),\ke)}. Suppose that 
\m{A_{mn}} is represented by the cocycle \m{(\alpha_{ij}^{mn})_{i,j\in 
I}},\Nligne \m{\alpha_{ij}^{mn}\in H^0(X_{ij}\cap\pi^{-1}(S_{mn}),\ke)}. The 
cocycle condition for \m{(A_{mn})} implies that \Nligne 
\m{(\alpha_{ij}^{mn}+\alpha_{ij}^{np}-\alpha_{ij}^{mp})_{i,j\in I}} represents 
0 in \m{H^1(\pi^{-1}(S_{mnp}),\ke)}. Hence we can write
\[\alpha_{ij}^{mn}+\alpha_{ij}^{np}-\alpha_{ij}^{mp} \ = \ 
\beta_i^{mnp}-\beta_j^{mnp} \ , \]
where \ \m{\beta_k^{mnp}\in H^0(\pi^{-1}(S_{mnp})\cap X_k,\ke)} \ for every 
\m{k\in I}. We have, for every distinct \Nligne \m{m,n,p,q\in M}, \m{i,j\in I}
\[\beta_i^{mnp}-\beta_i^{mnq}+\beta_i^{mpq}-\beta_i^{npq} \ = \
\beta_j^{mnp}-\beta_j^{mnq}+\beta_j^{mpq}-\beta_j^{npq} \ .\]
Hence there exists \ \m{\Delta^{mnpq}\in H^0(\pi^{-1}(S_{mnpq}),\ke)} \ such 
that
\[\beta_i^{mnp}-\beta_i^{mnq}+\beta_i^{mpq}-\beta_i^{npq} \ = \ 
\Delta^{mnpq}_{|\pi^{-1}(S_{mnpq})\cap V_i}\]
for every \m{i\in I}. Since \m{\pi_*(\ke)=0}, we have \ \m{\Delta^{mnpq}=0}. 
Hence for every \m{i\in I}, \m{(\beta_i^{mnp})_{m,n,p\in M}} satisfies the 
cocycle condition, and defines an element of \m{H^2(X_i,\ke)=0} (because 
\m{X_i} is affine). So we can write, for every distinct \m{m,n,p\in M}, \m{i\in 
I}, \[\beta_i^{mnp} \ = \ \Gamma_i^{mn}+\Gamma_i^{np}-\Gamma_i^{mp} \ , \]
with \ \m{\Gamma_i^{mn}\in H^0(\pi^{-1}(S_{mn})\cap X_i,\ke)}. Let
\[a^{mn}_{ij} \ = \ \alpha_{ij}^{mn}-\Gamma_i^{mn}+\Gamma_j^{mn} \ . \]
Then \m{A_{mn}} is also represented by \m{(a^{mn}_{ij})_{i,j\in I}}. For every 
distinct \m{i,j\in I}, we have \Nligne 
\m{a^{mn}_{ij}+a^{np}_{ij}-a^{mp}_{ij}=0}, 
hence \m{(a^{mn}_{ij})_{m,n\in M}} represents an element of \m{H^1(X_{ij},\ke)},
which is zero (because \m{X_{ij}} is affine). So we can write
\[a_{ij}^{mn} \ = \ \mu_{ij}^m-\mu_{ij}^n \ , \]
\m{\mu_{ij}^m\in H^0(\pi^{-1}(S_m)\cap X_{ij})}. We have, for every distinct 
\m{i,j,k\in I}, \m{m,n\in M},
\[\mu_{ij}^m+\mu_{jk}^m-\mu_{ik}^m \ = \ \mu_{ij}^n+\mu_{jk}^n-\mu_{ik}^n \]
on \m{\pi^{-1}(S_{mn})\cap X_{ijk}}. Hence there exists \ \m{B_{ijk}\in 
H^0(X_{ijk},\ke)} such that \Nligne \m{\mu_{ij}^m+\mu_{jk}^m-\mu_{ik}^m=B_{ijk|
\pi^{-1}(S_m)\cap X_{ijk}}}. Then \m{(B_{ijk})_{i,j,k\in I}} verifies the 
cocycle condition, and defines \ \m{\lambda(\theta)\in H^2(X,\ke)}.
\end{subsub}

\sepprop

With these definitions, it is easy to see that $\lambda$ is injective, $\eta$ 
surjective and \m{\ker(\eta)=\imm(\lambda)}.

\end{sub}

\sepsub

\Ssect{The exact sequence of Ext's}{Ex_Ext}

Let $X$ be a scheme and $E$, $F$ be coherent sheaves on $X$. From \cite{go}, 
7.3 (Ext spectral sequence), we have an exact sequence \m{\bm{\Gamma}_{E,F}}
\xmat{0\to H^1(X,\HHom(E,F))\ar[r]^-{\alpha_X} & 
\Ext^1_{\ko_X}(E,F)\ar[r]^-{\beta_X} & 
H^0(\EExt^1_{\ko_X}(E,F))\ar[r]^-{\delta_X} & H^2(X,\HHom(E,F))}
We will give a simple description of this exact sequence using \v Cech 
cohomology. 

\sepsubsub

\begin{subsub}\label{alpha_x} Description of \m{\alpha_X} -- \rm
Let \m{(X_i)_{i\in I}} be an open affine cover of $X$, and \Nligne 
\m{u\in H^1(X,\HHom(E,F))}, represented by the cocycle \m{(\phi_{ij})} with 
respect to the cover \m{(X_i)}, with \m{\phi_{ij}:E_{|X_{ij}}\to F_{|X_{ij}}}.
Let $\ke$ be the coherent sheaf on $X$ obtained by gluing the sheaves 
\m{(E\oplus F)_{|X_i}} using the automorphisms of \m{(E\oplus F)_{|X_{ij}}} 
defined by the matrices \m{\begin{pmatrix}I_E & 0\\ \phi_{ij} & 
I_F\end{pmatrix}}. Then we have an obvious exact sequence \ \m{0\to F\to\ke\to 
E\to 0}, whose associated element in \m{\Ext^1_{\ko_X}(E,F)} is \m{\alpha_X(u)}.

Let \m{\sigma\in\Ext^1_{\ko_X}(E,F)}, corresponding to the extension
\xmat{0\ar[r] & F\ar[r]^-\lambda & \ke\ar[r]^-\mu & E\ar[r] & 0.}
Let \m{u\in H^1(X,\HHom(E,F))}, represented by the cocycle 
\m{(\phi_{ij})}. Let \m{\sigma'=\sigma+\alpha_X(u)}, corresponding to the 
extension \ \m{0\to F\to\ke'\to E\to 0}. Then \m{\ke'} is obtained by gluing 
the sheaves \m{\ke_{|X_i}} using the automorphisms \ 
\m{I+\lambda\circ\phi_{ij}\circ\mu:\ke_{|X_{ij}}\to \ke_{|X_{ij}}}.
\end{subsub}

\sepsubsub

\begin{subsub}\label{beta_x} Description of \m{\beta_X} -- \rm
Let \m{\sigma\in\Ext^1_{\ko_X}(E,F)}, corresponding to an exact sequence \Nligne
\m{0\to F\to\ke\to E\to 0}, whose restriction to \m{X_i} defines \Nligne 
\m{\sigma_i\in\Ext^1_{\ko_{X_i}}(E_{|X_i},F_{|X_i})=\EExt^1_{\ko_X}(E,F)(X_i)}. 
then \m{\beta_X(\sigma)} is the section $s$ of \m{\EExt^1_{\ko_X}(E,F)} such 
that for every \m{i\in I}, \m{s_{|X_i}=\sigma_i}.
\end{subsub}

\sepsubsub

\begin{subsub}\label{delta_x} Description of \m{\delta_X} -- \rm
Let \ \m{s\in H^0(\EExt^1_{\ko_X}(E,F))}. For every \m{i\in I}, we have \Nligne 
\m{s_{|X_i}\in \Ext^1_{\ko_{X_i}}(E_{|X_i},F_{|X_i})}, which defines an exact 
sequence
\xmat{\Sigma_i:0\ar[r] & F_{|X_i}\ar[r]^-{\gamma_i} & \ke_i\ar[r]^-{p_i} & 
E_{|X_i}\ar[r] & 0 \ .}
The restrictions of \m{\Sigma_i} and \m{\Sigma_j} to \m{X_{ij}} define the same 
element of \m{\Ext^1_{\ko_{X_{ij}}}(E_{|X_{ij}},F_{|X_{ij}})}, i.e. we have a 
commutative diagram
\xmat{\Sigma_i: \ 0\ar[r] & F_{|X_{ij}}\ar[r]^-{\gamma_i}\fleq[d] & 
\ke_{i|X_{ij}}\ar[r]^-{p_i}\ar[d]^{\psi_{ij}} & E_{|X_{ij}}\ar[r]\fleq[d] & 0\\
\Sigma_j: \ 0\ar[r] & F_{|X_{ij}}\ar[r]^-{\gamma_j} & 
\ke_{j|X_{ij}}\ar[r]^-{p_j} & E_{|X_{ij}}\ar[r] & 0 }
We have then \ 
\m{\psi_{jk}\psi_{ij}-\psi_{ik}:\ke_{i|X_{ijk}}\to\ke_{k|X_{ijk}}}, and there 
exists \ \m{\sigma_{ijk}:E_{|X_{ijk}}\to F_{|X_{ijk}}} \ such that \ 
\m{\psi_{jk}\psi_{ij}-\psi_{ik}=\gamma_k\sigma_{ijk} p_i}. It is easily checked 
that the family \m{(\sigma_{ijk})} satisfies the cocycle relation, and then 
defines an element of \m{H^2(X,\HHom(E,F))}, which is \m{\delta_X(s)}.
\end{subsub}

\sepsubsub

\begin{subsub}\label{loc_ext} Local extensions -- \rm Let \ \m{s\in 
H^0(\EExt^1_{\ko_X}(E,F))}. Let \m{(U_m)_{m\in M}} be an open cover of $X$ such 
that for every \m{i\in I} and \m{m\in M}, the intersection \m{X_i\cap U_m} is 
affine. For every \m{i,j\in I} and \m{m,n\in M}, let \ \m{Z_i^m=X_i\cap U_m}, 
\m{Z_{ij}^m=X_{ij}\cap U_m}, etc.

Given two extensions
\xmat{\Sigma:0\ar[r] & F\ar[r]^-\lambda & \ke\ar[r]^-\mu & E\ar[r] & 0, \quad 
& \Sigma':0\ar[r] & F\ar[r]^-{\lambda'} &\ke'\ar[r]^-{\mu'} & E\ar[r] & 0 \ ,}
an isomorphism \m{\phi:\ke\to\ke'} is called a {\em 
I-isomorphism} (with respect to $\Sigma$, \m{\Sigma'}) if the following diagram
\xmat{0\ar[r] & F\ar[r]^-\lambda\fleq[d] & \ke\ar[r]^-\mu\ar[d]^\phi & 
E\ar[r]\fleq[d] & 0\\
0\ar[r] & F\ar[r]^-{\lambda'} & \ke'\ar[r]^-{\mu'} & E\ar[r] & 0}
is commutative. If \m{\phi'} is another I-isomorphism \m{\ke\to\ke'}, there 
exists a morphism \ \m{\tau:E\to F} \ such that \ 
\m{\phi-\phi'=\lambda'\tau\mu}.

Suppose that for every \m{m\in M} there exists \ 
\m{\sigma_m\in\Ext^1_{\ko_{U_m}}(E,F)} \ such that \ 
\m{\delta_{U_m}(\sigma_m)=s_{|U_m}}. Let
\xmat{0\ar[r] & F_{|U_m}\ar[r]^-{c_m} & \kf_m\ar[r]^-{\pi_m}\ar[r] & 
E_{|U_m}\ar[r] & 0}
be the extension corresponding to \m{\sigma_m}.
Then for every \m{i\in I}, \m{m\in M}, 
there exists an I-isomorphism \ \m{\alpha_i^m:\ke_{i|Z^m_i}\to\kf_{m|Z^m_i}}.
Let \m{j\in I}. Then \m{(\alpha_j^m)^{-1}\alpha_i^m} is an I-isomorphism \ 
\m{\ke_{i|Z^m_{ij}}\to\ke_{j|Z^m_{ij}}}. So we can write
\[\psi_{ij}-(\alpha_j^m)^{-1}\alpha_i^m \ = \ \gamma_j\mu^m_{ij}p_i \ , \]
with \ \m{\mu_{ij}^m:E_{|Z_{ij}^m}\to F_{|Z_{ij}^m}}. It follows that 
\[\big(\psi_{jk}\psi_{ij}-\psi_{ik}\big)_{|Z_{ijk}^m} \ = \ 
\gamma_k(\mu_{jk}^m+\mu_{ij}^m-\mu_{ik}^m)_{|Z_{ijk}^m}p_i \ , \]
so we have \ 
\[\sigma_{ijk|Z_{ijk}^m} \ = \ 
(\mu_{jk}^m+\mu_{ij}^m-\mu_{ik}^m)_{|Z_{ijk}^m} \ . \]
Let \ \m{a_{ij}^{mn}=\mu_{ij}^m-\mu_{ij}^n}. Then
for every \m{m,n\in M}, \m{(a^{mn}_{ij})_{i,j\in I}} is a cocycle and thus 
defines \ \m{A_{mn}\in H^1(U_{mn},\HHom(E,F))}.
\end{subsub}

\sepprop

\begin{subsub}\label{case_fam} The case of families -- \rm 
We suppose now that $X$ is a product : \m{X=Y\times S}, and that \ 
\m{p_{S*}(\HHom(E,F))=0}, where \ \m{p_S:X\to S} \ is the projection. We have 
from \ref{Ex_Ler} a canonical injection
\[\lambda : H^1(S,R^1p_{S*}(\HHom(E,F)))\lra H^2(X,\HHom(E,F)) \ . \]
Let \m{(S_m)_{m\in M}} be an affine open cover of $S$. We suppose that \ 
\m{U_m=p_S^{-1}(S_m)} \ for every \m{m\in M}. Then \m{(A_{mn})_{m,n\in M}} 
defines \ \m{\epsilon\in H^1(S,R^1p_{S*}(\HHom(E,F)))}. It follows from 
\ref{Ex_Ler} that 
\end{subsub}

\sepprop

\begin{subsub}\label{prop5}{\bf Proposition: } We have \ 
\m{\delta_X(s)=\lambda(\epsilon)}.
\end{subsub}

\sepprop

\begin{subsub} Generalization -- \rm If we don't assume that  \ 
\m{p_{S*}(\HHom(E,F))=0}, we have still an injective morphism
\[\lambda':H^1(S,R^1p_{S*}(\HHom(E,F)))\lra 
H^2(X,\HHom(E,F))/H^0(S,p_{S*}(\HHom(E,F))
\ , \]
and we have an analogous result (using a more complicated proof), i.e. 
\m{\lambda'(\epsilon)} is the class of \m{\delta_X(s)}.

\end{subsub}

\end{sub}

\sepsec

\section{Primitive multiple schemes}\label{PMS}

\Ssect{Definition and construction}{PMS_def}

Let $X$ be a smooth connected variety, and \ \m{d=\dim(X)}. A {\em multiple 
scheme with support $X$} is a Cohen-Macaulay scheme $Y$ such that 
\m{Y_{red}=X}. If $Y$ is quasi-projective we say that it is a {\em multiple 
variety with support $X$}. In this case $Y$ is projective if $X$ is.

Let $n$ be the smallest integer such that \m{Y=X^{(n-1)}}, \m{X^{(k-1)}}
being the $k$-th infinitesimal neighborhood of $X$, i.e. \
\m{\ki_{X^{(k-1)}}=\ki_X^{k}} . We have a filtration \ \m{X=X_1\subset
X_2\subset\cdots\subset X_{n}=Y} \ where $X_i$ is the biggest Cohen-Macaulay
subscheme contained in \m{Y\cap X^{(i-1)}}. We call $n$ the {\em multiplicity}
of $Y$.

We say that $Y$ is {\em primitive} if, for every closed point $x$ of $X$,
there exists a smooth variety $S$ of dimension \m{d+1}, containing a 
neighborhood of $x$ in $Y$ as a locally closed subvariety. In this case, 
\m{L=\ki_X/\ki_{X_2}} is a line bundle on $X$, \m{X_j} is a primitive multiple 
scheme of multiplicity $j$ and we have \ 
\m{\ki_{X_j}=\ki_X^j}, \m{\ki_{X_{j}}/\ki_{X_{j+1}}=L^j} \ for \m{1\leq j<n}. 
We call $L$ the line bundle on $X$ {\em associated} to $Y$. The ideal sheaf 
\m{\ki_{X,Y}} can be viewed as a line bundle on \m{X_{n-1}}.

Let \m{P\in X}. 
Then there exist elements \m{y_1,\ldots,y_d}, $t$ of \m{m_{S,P}} whose images 
in \m{m_{S,P}/m_{S,P}^2} form a basis, and such that for \m{1\leq i<n} we have 
\ \m{\ki_{X_i,P}=(t^{i})}. In this case the images of \m{y_1,\ldots,y_d} in 
\m{m_{X,P}/m_{X,P}^2} form a basis of this vector space.

A {\em multiple scheme with support $X$} is primitive if and only if 
\m{\ki_X/\ki_X^2} is zero or a line bundle on $X$ (cf. \cite{dr7}, Proposition 
2.3.1).  

Even if $X$ is projective, we do not assume that $Y$ is projective. In fact we 
will see examples of non quasi-projective $Y$.

The simplest case is when $Y$ is contained in a smooth variety $S$ of dimension 
\m{d+1}. Suppose that $Y$ has multiplicity $n$. Let \m{P\in X} and 
\m{f\in\ko_{S,P}}  a local equation of $X$. Then we have \ 
\m{\ki_{X_i,P}=(f^{i})} \ for \m{1<i\leq n} in $S$, in particular 
\m{\ki_{Y,P}=(f^n)}, and \ \m{L=\ko_X(-X)} .


For any \m{L\in\Pic(X)}, the {\em trivial primitive variety} of multiplicity 
$n$, with induced smooth variety $X$ and associated line bundle $L$ on $X$ is 
the $n$-th infinitesimal neighborhood of $X$, embedded by the zero section in 
the dual bundle $L^*$, seen as a smooth variety.

\sepsubsub

\begin{subsub}\label{PMS-1} Construction of primitive multiple schemes -- \rm
Let $Y$ be a primitive multiple scheme of multiplicity $n$, \m{X=Y_{red}}.
Let \ \m{{\bf Z}_n=\spec(\C[t]/(t^n))}.
Then for every closed point \m{P\in X}, there exists an open neighborhood $U$ 
of $P$ in $X$, such that if \m{U^{(n)}} is the corresponding neighborhood of 
$P$ in $Y$, there exists a commutative diagram
 \xmat{ & U\flinc[ld]\flinc[rd] \\
U^{(n)}\ar[rr]^-\simeq & & U\times {\mathbf Z}_n}
i.e. $Y$ is locally trivial (\cite{dr1}, Th\'eor\`eme 5.2.1, Corollaire 5.2.2).

It follows that we can construct a primitive multiple scheme of multiplicity 
$n$ by taking an open cover \m{(U_i)_{i\in I}} of $X$ and gluing the varieties 
\ \m{U_i\times{\bf Z}_n} (with automorphisms of the \ \m{U_{ij}\times{\bf Z}_n} 
\ leaving \m{U_{ij}} invariant).

Let \m{(U_i)_{i\in I}} be an affine open cover of $X$ such that we have 
trivializations
\xmat{\delta_i:U_i^{(n)}\ar[r]^-\simeq & U_i\times{\bf Z}_n , }
and \ \m{\delta_i^*:\ko_{U_i\times{\bf Z}_n}\to\ko_{U_i^{(n)}}} \ the 
corresponding isomorphism. Let
\xmat{\delta_{ij}=\delta_j\delta_i^{-1}:U_{ij}\times{\bf Z}_n\ar[r]^-\simeq & 
U_{ij}\times{\bf Z}_n \ . }
Then \ \m{\delta_{ij}^*=\delta_i^{*-1}\delta_j^*} \ is an automorphism of \ 
\m{\ko_{U_i\times Z_n}=\ko_X(U_{ij})[t]/(t^n)}, such that for every \ 
\m{\phi\in\ko_X(U_{ij})}, seen as a polynomial in $t$ with coefficients in 
\m{\ko_X(U_{ij})}, the term of degree zero of \m{\delta_{ij}^*(\phi)} is the 
same as the term of degree zero of $\phi$.
\end{subsub}

\sepsubsub

\begin{subsub}\label{I_X} The ideal sheaf of $X$ -- \rm
There exists \ \m{\alpha_{ij}\in\ko_X(U_{ij})\ot_\C\C[t]/(t^{n-1})} \ such that 
\ \m{\delta_{ij}^*(t)=\alpha_{ij}t}. Let \ 
\m{\alpha^{(0)}_{ij}=\alpha_{ij|X}\in\ko_X(U_i)}.
For every \m{i\in I}, \m{\delta_i^*(t)} is a generator of 
\m{\ki_{X,Y|{U^{(n)}}}}. So we have local trivializations
\[\xymatrix@R=5pt{\lambda_i:\ki_{X,Y|{U_i^{(n-1)}}}\ar[r] & 
\ko_{U_i^{(n-1)}}\\
\delta_i^*(t)\fmaps[r] & 1}\]
Hence \ \m{\lambda_{ij}=\lambda_i\lambda_j^{-1}: 
\ko_{U_{ij}^{(n-1)}}\to\ko_{U_{ij}^{(n-1)}}} \ is the multiplication by 
\m{\delta_j^*(\alpha_{ij})}. It follows that 
\m{(\delta_j^*(\alpha_{ij}))} (resp. 
\m{(\alpha^{(0)}_{ij})})  is a cocycle representing the line bundle 
\m{\ki_{X,Y}} (resp. \m{L}) on \m{X_{n-1}} (resp. $X$). 
\end{subsub}

\sepsubsub

\begin{subsub}\label{ext-I} The associated sheaves of non abelian groups and 
obstructions to extension in higher multiplicity -- \rm 
For every open subset $U$ of $X$, we have
\ \m{\ko_{X\times{\bf Z}_n}(U)=\ko_X(U)[t]/(t^n)}. Let  \m{\kg_n} be
the sheaf of (non abelian) groups on $X$ defined by: for every open subset $U$ 
of $X$, \m{\kg_n(U)} is the group of automorphisms $\theta$ of the $\C$-algebra 
\m{\ko_X(U)[t]/(t^n)} such that for every \m{\alpha\in\ko_X(U)[t]/(t^n)}, 
\m{\theta(\alpha)_{|U}=\alpha_{|U}}. We have \m{\kg_1=\ko_X^*}. Then 
\m{(\delta_{ij}^*)} is a cocycle of \m{\kg_n}, which describes completely 
\m{X_n}. 

In this way we see that there is a canonical bijection between the cohomology 
set \m{H^1(X,\kg_n)} and the set of isomorphism classes of primitive multiple 
schemes \m{X_n} such that \m{X=(X_n)_{red}}. There is an obvious surjective 
morphism \ \m{\rho_{n+1}:\kg_{n+1}\to\kg_n}, such that if \m{n\geq 2},
\[H^1(\rho_{n+1}):H^1(X,\kg_{n+1})\lra H^1(X,\kg_n)\]
sends a primitive multiple scheme of multiplicity \m{n+1} to the underlying 
scheme of multiplicity $n$, whereas 
\[H^1(\rho_2):H^1(X,\kg_2)\lra H^1(X,\ko_X^*)=\Pic(X)\]
sends \m{X_2} to $L$.

We have \m{\ker(\rho_2)\simeq T_X} \ and
\ \m{\ker(\rho_{n+1})\simeq T_X\oplus\ko_X} \ if \m{n\geq 2}.
The fact that they are sheaves of {\em abelian} groups allows to compute {\em 
obstructions}. Let \m{g_n\in H^1(X,\kg_n)}, corresponding to the primitive 
multiple scheme \m{X_n}, and for \m{1\leq i<n}, \m{g_i} the image of \m{g_n} in 
\m{H^1(X,\kg_i)}. Let  \m{\ker(\rho_{n+1})^{g_n}} be the associated sheaf of 
groups (cf. \cite{dr10}, 2.2, \cite{fr}). We have then
\[\ker(\rho_2)^{g_1} \ \simeq \ T_X\ot L  \] 
By cohomology theory, we find that there is a 
canonical surjective map
\[H^1(X,T_X\ot L)\lra H^1(\rho_2)^{-1}(L)\]
sending 0 to the trivial primitive scheme, and whose fibers are the orbits of 
the action of $\C^*$ on \m{H^1(X,T_X\ot L)} by multiplication. Hence there is a 
bijection between the set of non trivial double schemes with associated line 
bundle $L$, and \m{\P(H^1(X,T_X\ot L))}. For a non trivial double scheme 
\m{X_2}, we have an exact sequence 
\[0\lra L\lra\Omega_{X_2|X}\lra\Omega_X\lra 0 \ , \]
corresponding to \ \m{\sigma\in\Ext^1_{\ko_X}(\Omega_X,L)=H^1(T_X\ot L)}, and 
\m{\C\sigma} is the element of \m{\P(H^1(T_X\ot L))} corresponding to \m{X_2}.
We can write for \m{n=2}
\[(\delta_{ij}^*)_{|\ko_X(U_{ij})} \ = \ 
I_{|\ko_X(U_{ij})}+tD_{ij} \ , \]
where \m{D_{ij}} is a derivation of \m{\ko_X(U_{ij})}. Then the family 
\m{(D_{ij})} represents $\sigma$ in the sense of \ref{cech}.

If \m{n>2} we have
\[\ker(\rho_{n+1})^{g_n} \ \simeq \ (\Omega_{X_2|X})^*\ot L^n\]
We have then (by the theory of cohomology of 
sheaves of groups) an {\em obstruction map}
\[\Delta_n:H^1(X,\kg_n)\lra H^2((\Omega_{X_2|X})^*\ot L^n)\]
such that \ \m{g_n\in\imm(H^1(\rho_{n+1}))} \ if and only \ \m{\Delta_n(g_n)=0}.
\end{subsub}

\sepsubsub

\begin{subsub}\label{ext_id} Extensions of the ideal sheaf of $X$ -- \rm The 
ideal sheaf 
\m{\ki_{X,X_n}} is a line bundle on \m{X_{n-1}}. A necessary condition to 
extend \m{X_n} to a primitive multiple scheme \m{X_{n+1}} of multiplicity 
\m{n+1} is that \m{\ki_{X,X_n}} can be extended to a line bundle on \m{X_n} 
(namely \m{\ki_{X,X_{n+1}}}). This is why we can consider pairs \m{(X_n,\L)}, 
where $\L$ is a line bundle on \m{X_n} such that \ 
\m{\L_{|X_{n-1}}\simeq\ki_{X,X_n}}. 

The corresponding sheaf of groups on $X$ is defined as follows: for every open 
subset \m{U\subset X}, \m{\kh_n(U)} is the set of pairs \m{(\phi,u)}, where 
\m{\phi\in\kg_n(U)}, and \m{u\in\ko_X(U)[t]/(t^n)} is such that \ 
\m{\phi(t)=ut} (cf. \cite{dr10}, 4.5). The set of isomorphism classes of the 
above pairs \m{(X_n,\L)} 
can then be identified with the cohomology set \m{H^1(X,\kh_n)}.

There is an obvious morphism \ \m{\tau_n:\kg_{n+1}\to\kh_n}, such that
\[H^1(\tau_n):H^1(X,\kg_{n+1})\lra H^1(X,\kh_n)\]
sends \m{X_{n+1}} to \m{(X_n,\ki_{X,X_{n+1}})}.
Let \m{g\in H^1(X,\kh_n)}. Then \ 
\m{\ker(\tau_n)\simeq T_X} \ and
\[\ker(\tau_n)^g \ \simeq \ T_X\ot L^n \ . \]
Consequently there is again, by cohomology theory, an {\em obstruction map}
\[\Delta''_n:H^1(X,\kh_n)\lra H^2(T_X\ot L^n)\]
such that, if \m{(X_n,\L)} corresponds to $g$, there is an extension of \m{X_n} 
to a primitive multiple scheme \m{X_{n+1}} of multiplicity \m{n+1} with \ 
\m{\ki_{X,X_{n+1}}\simeq\L} \ if and only if \ \m{\Delta''_n(g)=0}.
\end{subsub}

\end{sub}

\sepsub

\Ssect{Descriptions using the open cover \m{(U_i)}}{const_sh}

{\bf (i)} \ {\em Construction of sheaves and morphisms -- } (cf. \cite{dr10}, 
2.1.3). Let $\ke$ be a 
coherent sheaf on \m{X_n}. We can define it in the usual way, starting with 
sheaves \m{\kf_i} on the open sets \m{U_i^{(n)}} and gluing them. We take 
these sheaves of the form \ \m{\kf_i=\delta_i^*(\ke_i)}, where \m{\ke_i} is a 
sheaf on \ \m{U_i\times {\bf Z}_n}. To glue the sheaves \m{\kf_i} on the 
intersection 
\m{U_{ij}^{(n)}} we use isomorphisms \ 
\m{\rho_{ij}:\kf_{j|U_{ij}^{(n)}}\to\kf_{i|U_{ij}^{(n)}}}, with the relations \ 
\m{\rho_{ik}=\rho_{jk}\rho_{ij}}. Let
\[\theta_{ij} = 
(\delta_i^*)^{-1}(\rho_{ij}):\delta_{ij}^*(\ke_{j|U_{ij}\times{\bf 
Z}_n})\lra\ke_{i|U_{ij}\times{\bf Z}_n} \ .
\]
We have then the relations \ 
\m{\theta_{ik}=\theta_{ij}\circ\delta_{ij}^*(\theta_{jk})}. Conversely, 
starting with sheaves \m{\ke_i} and isomorphisms \m{\theta_{ij}} satisfying the 
preceding relations, one obtains a coherent sheaf on \m{X_n}.

This applies to trivializations, i.e when \ \m{\ke_i=\ko_{U_i^{(n)}}\ot\C^r}. 
We have then \Nligne \m{\theta_{ij}:\ko_{U_{ij}\times{\bf 
Z}_n}\ot\C^r\to\ko_{U_{ij}\times{\bf Z}_n}\ot\C^r}. In particular, 
\m{\ki_{X,X_n}} is represented by \m{(\alpha_{ij})}.

Suppose that we have another sheaf \m{\ke'} on \m{X_n}, defined by sheaves 
\m{\ke'_i} on \ \m{U_i\times {\bf Z}_n} \ and isomorphisms \m{\theta'_{ij}}. 
One can see easily that a morphism \ \m{\Psi:\ke\to\ke'} \ is defined by 
morphisms \ \m{\Psi_i:\ke_i\to\ke'_i} \ such that \ 
\m{\theta'_{ij}\circ\delta_{ij}^*(\Psi_j)=\Psi_i\circ\theta_{ij}}.

\sepsubsub

{\bf (ii)} \ {\em Other construction of sheaves -- } Let $E$ be a vector bundle 
on $X$. For every open subset \m{i\in I}, let \ \m{p_i:U_i\times{\bf Z}_n\to 
U_i} \ be the projection. Then we have \ 
\m{\delta_i^*(p_i^*(E_{|U_i}))_{|U_i}=E_{|U_i}}.
We construct a vector bundle $\E$ on \m{X_n} by gluing the 
vector bundles \m{\delta_i^*(p_i^*(E_{|U_i}))} on \m{U_i^{(n)}}. For this we 
take isomorphisms
\[\epsilon_{ij}:\delta_j^*(p_j^*(E_{|U_{ij}}))_{U_{ij}^{(n)}}\lra
\delta_i^*(p_i^*(E_{|U_{ij}}))_{U_{ij}^{(n)}}\]
satisfying \ \m{\epsilon_{ij}\epsilon_{jk}=\epsilon_{ik}}. Let
\[\tau_{ij}={\delta_i^*}^{-1}(\epsilon_{ij}):\delta_{ij}^*(p_j^*(E_{|U_{ij}}))=
p_i^*(E_{|U_{ij}})\lra p_i^*(E_{|U_{ij}}) \ . \]
Then we have \ \m{\tau_{ij}\circ\delta_{ij}^*(\tau_{jk})=\tau_{ik}}. 
Conversely, given automorphisms \m{\tau_{ij}} satisfying the preceding 
relation, 
we can define \m{\epsilon_{ij}} and the vector bundle $\E$ on \m{X_n}.

Suppose that the \m{\tau_{ij}} are homotheties (multiplication by \m{\nu_{ij}
\in\ko(U_{ij}\times{\bf Z}_n)^*}). Then the family \m{(\nu_{ij})} defines a 
line bundle $\L$ on \m{X_n} (by (i)), and we have \ \m{\E_{|X}\simeq 
E\ot\L_{|X}}.

\sepsubsub

{\bf (iii)} \ {\em Cohomology of vector bundles -- } Suppose that a vector 
bundle $E$ on \m{X_n} is built using isomorphisms \ 
\m{\theta_{ij}:\ko_{U_{ij}\times{\bf 
Z}_n}\ot\C^r\to\ko_{U_{ij}\times{\bf Z}_n}\ot\C^r} \ as in (i).
Any \ \m{\beta\in H^1(X_n,E)} \ is represented by a family \m{(\beta_{ij})}, 
\m{\beta_{ij}\in H^0(\ko_{U_{ij}\times{\bf Z}_n}\ot\C^r)}, with the cocycle 
relations \ \m{\beta_{ik}=\beta_{ij}+\theta_{ij}\delta_{ij}^*(\beta_{jk})}.

 
\end{sub}

\sepsub

\Ssect{Connecting morphisms}{con_hom}

Recall that on \m{X_2} we can write \ 
\m{\delta_{ij|\ko_X(U_{ij})}^*=I+tD_{ij}}, where \m{D_{ij}} is a derivation 
of \m{\ko_X(U_{ij})} (cf. \ref{ext-I}).

Let $\E$ be a vector bundle on \m{X_2} of rank $r$, defined as in 
\ref{const_sh}, (i), using isomorphisms \Nligne 
\m{\theta_{ij}:\ko_{U_{ij}\times{\bf Z}_2}\ot\C^r\to\ko_{U_{ij}\times{\bf 
Z}_2}\ot\C^r}. We can write
\[\theta_{ij} \ = \ \theta_{ij}^{(0)}+\theta_{ij}^{(1)}t \ , \]
where \m{\theta_{ij}^{(0)}}, \m{\theta_{ij}^{(1)}} are matrices with 
coefficients in \m{\ko_X(U_{ij})}.

Let \ \m{E=\E_{|X}}. We have a canonical exact sequence
\[0\lra E\ot L\lra\E\lra E\lra 0 \ , \]
whence a connecting morphism
\[\delta^1_\E:H^1(X,E)\lra H^2(X,E\ot L) \ . \]
Let \ \m{\beta\in H^1(X,E)}, represented by a family \m{(\beta_{ij})}, 
\m{\beta_{ij}\in H^0(\ko_{U_{ij}}\ot\C^r)}, in the sense of 
\ref{cech}.

\sepprop

\begin{subsub}\label{lem10}{\bf Lemma: } \m{\delta^1_\E(\beta)} is represented, 
in the sense of \ref{const_sh}, (iii), by the family \m{(\nu_{ijk})}, with
\[\nu_{ijk} \ = \ 
\theta_{ij}^{(0)}D_{ij}(\beta_{jk})+\theta_{ij}^{(1)}\beta_{jk} \ . \]
\end{subsub}

\begin{proof}
We construct \m{\delta^1_\E(\beta)} as follows: we view \m{\beta_{ij}}, 
\m{\beta_{kj}}, \m{\beta_{ik}} as elements of \m{H^0(\ko_{U_{ij}\times{\bf 
Z}_n}\ot\C^r)}. Then we can take
\[\nu_{ijk} \ = \ \beta_{ij}-\beta_{ik}+\theta_{ij}\delta_{ij}^*(\beta_{jk})
\]
which gives Lemma \ref{lem10}.
\end{proof}

\sepprop

Suppose that \m{X_n} can be extended to a primitive multiple scheme \m{X_{n+1}} 
of multiplicity \m{n+1}. 

We have an exact sequence of sheaves on \m{X_{n+1}}
\[0\lra L^n\lra\ko_{X_{n+1}}\lra\ko_{X_n}\lra 0 \ , \]
whence a connecting morphism
\[\delta^0:H^0(\ko_{X_n})\lra H^1(X,L^n) \ . \]
We have \ \m{H^0(X,L^{n-1})\subset H^0(\ko_{X_n})}. Let \ \m{\eta\in 
H^0(X,L^{n-1})}, represented by a cocycle \m{(\eta_i)}, 
\m{\eta_i\in\ko_X(U_i)} (in the sense of \ref{cech1}, i.e. the cocycle relation 
is \ \m{\eta_i=(\alpha_{ij}^{(0)})^{n-1}\eta_j}). Recall that 
\m{\ki_{X,X_{n+1}}} is represented by the cocycle \m{(\alpha_{ij})}, where 
\m{\alpha_{ij}} is an invertible element of \m{\ko_X(U_{ij})[t]/(t^n)}.

The proof of the following lemma is similar to that of Lemma \ref{lem10}:

\sepprop

\begin{subsub}\label{lem17}{\bf Lemma: } \m{\delta^0(\eta)} is 
represented by the cocycle \m{(\gamma_{ij})}, with
\[\gamma_{ij} \ = \ (\alpha^{(0)}_{ij})^{n-1}D_{ij}(\eta_j)+ 
(n-1)(\alpha^{(0)}_{ij})^{n-2}\alpha^{(1)}_{ij}\eta_j \ . \]
\end{subsub}

\end{sub}

\sepsub

\Ssect{The sheaf of differentials}{PMS_Omeg}

We use the notations of \ref{PMS-1}. We will give a construction of 
\m{\Omega_{X_n}} using \ref{const_sh} : it is defined by the sheaves \ 
\m{\ke_i=\Omega_{U_i\times{\bf Z}_n}}, and
\[{\bf d}(\delta_{ij}):\delta_{ij}^*(\Omega_{U_i\times{\bf Z}_n})=
\Omega_{U_i\times{\bf Z}_n}\lra\Omega_{U_i\times{\bf Z}_n}\]
(more precisely \ \m{{\bf 
d}(\delta_{ij})(a.db)=\delta_{ij}^*(a)d(\delta_{ij}^*(b)))}).

Let \m{\sigma\in H^1(X_n,\Omega_{X_n})}, represented by the cocycle 
\m{(\sigma_{ij})}, \m{\sigma_{ij}\in\Omega_{U_{ij}^{(n)}}}. Then \Nligne 
\m{\sigma'_{ij}={\bf d}(\delta_i^{-1})(\sigma_{ij})\in\Gamma(\Omega_{U_{ij}
\times{\bf Z}_n})} , and we have the relations
\[\sigma'_{ik} \ = \ \sigma'_{ij}+{\bf d}(\delta_{ij})(\sigma'_{jk}) \ . \]
Conversely, every family \m{(\sigma'_{ij})} satisfying these relations defines 
an element of \m{H^1(X_n,\Omega_{X_n})}.

\sepsubsub

\begin{subsub}\label{dual_sh} The dualizing sheaf -- \rm Suppose that $X$ is 
projective. Then \m{X_n} is a proper Cohen-Macaulay scheme. There exists 
a dualizing sheaf \m{\omega_{X_n}} for \m{X_n}, which is a line bundle (because 
\m{X_n} is locally a complete intersection).
\end{subsub}
\end{sub}

\sepsub

\Ssect{Extensions of primitive multiple schemes}{Ext_sh}

Let \m{X_n} a primitive multiple scheme of multiplicity \m{n\geq 2}, with
underlying smooth variety $X$ projective, irreducible, and 
associated line bundle $L$ on $X$. For every open subset \m{U\subset X}, let 
\m{U^{(n)}} denote the corresponding open subset of \m{X_n}. Let \ 
\m{\boldsymbol{\zeta}\in H^1(X,T_X\ot L)}, corresponding to the exact 
sequence \ \m{0\to L\to\Omega_{X_2|X}\to\Omega_X\to 0}.

If \m{X_n} can be extended to a primitive multiple scheme \m{X_{n+1}} of 
multiplicity \m{n+1}, then \m{\ki_{X,X_n}}, which is a line bundle on 
\m{X_{n-1}}, can be extended to a line bundle on \m{X_n}, namely 
\m{\ki_{X,X_{n+1}}}. We have exact sequences
\begin{equation}\label{equ12}0\lra L^n\lra\ki_{X,X_{n+1}}\lra\ki_{X,X_n}\lra 0 
\ , 
\end{equation}
\xmat{0\ar[r] & H^1(X,L^{n-1})\ar[r] & 
\Ext^1_{\ko_{X_n}}(\ki_{X,X_n},L^n)\ar[r]^-\phi & \End(L)\ar[r]^-\delta & 
H^2(X,L^{n-1}) }
(cf \ref{ext_4}), and given any exact sequence
\[0\lra L^n\lra\ke\lra\ki_{X,X_n}\lra 0 \ , \]
corresponding to \ \m{\sigma\in\Ext^1_{\ko_{X_n}}(\ki_{X,X_n},L^n)},
the sheaf $\ke$ is a line bundle on \m{X_n} if and only if \m{\phi(\sigma)} 
is an automorphism. In this case we can suppose that \ \m{\phi(\sigma)=I_L} .
In this way we can see that, given $(\ref{equ12})$, the extensions of 
\m{\ki_{X,X_n}} to a line bundle on \m{X_n} are parametrized by 
\m{H^1(X,L^{n-1})}. 

More precisely, using the notations of \ref{PMS-1} and \ref{const_sh}, let 
\m{\L_0\in\Pic(X_n)} be an extension of \m{\ki_{X,X_n}}. Suppose that it is 
defined  by the family \m{(\alpha_{ij})}, with \ 
\m{\alpha_{ij}\in[\ko_X(U_{ij})[t]/(t^n)]^*}, satisfying the relations \ 
\m{\alpha_{ik}=\alpha_{ij}\delta_{ij}^*(\alpha_{jk})}, so that $L$ is defined 
by the family \m{\big(\alpha_{ij}^{(0)}\big)}. Let $\L$ be an extension 
of \m{\ki_{X,X_n}} to a line bundle on \m{X_n}, corresponding to \ \m{\eta\in 
H^1(X,L^{n-1})}. Then $\L$ is defined by a cocycle \m{(\theta_{ij})}, with
\[\theta_{ij} \ = \ \alpha_{ij}+\beta_{ij}t^{n-1} \ , \]
\m{\beta_{ij}\in\ko_X(U_{ij})}. The cocycle relation \ 
\m{\theta_{ik}=\theta_{ij}\delta_{ij}^*(\theta_{jk})} \ is equivalent to
\begin{equation}\label{equ14}\frac{\beta_{ik}}{\alpha_{ik}^{(0)}} \ = \ 
\frac{\beta_{ij}}{\alpha_{ij}^{(0)}}+\big(\alpha_{ij}^{(0)}\big)^{n-1}
\frac{\beta_{jk}}{\alpha_{jk}^{(0)}} \ . \end{equation}
According to \ref{cech}, this means that 
\m{\dsp\Big(\frac{\beta_{ij}}{\alpha_{ij}^{(0)}}\Big)} defines an element of 
\m{H^1(X,L^{n-1})}, which is $\eta$.

Now we use the results of \ref{ext_id}. The family 
\m{(\delta_{ij}^*,\alpha_{ij})} defines an element $\bf h$ of \m{H^1(X,\kh_n)}. 
The family \m{(\delta_{ij}^*,\theta_{ij})} defines an element $\bf k$ 
of \m{H^1(X,\kh_n)} (corresponding to \m{(X_n,\L)}), and we have \ 
\m{\Delta''_n({\bf k})=0} \ if and only there exists an extension \m{X'_{n+1}} 
of \m{X_n} in multiplicity \m{n+1} such that \ \m{\ki_{X,X'_{n+1}}\simeq\L}. We 
will compute \m{\Delta''_n({\bf k}-{\bf h})}.

We can write
\[(\delta_{ij}^*)_{|\ko_X(U_{ij})} \ = \ 
I_{|\ko_X(U_{ij})}+tD_{ij}+t^2\Phi_{ij} \ , \]
where \m{D_{ij}} is a derivation of \m{\ko_X(U_{ij})}. The family \m{(D_{ij})} 
represents the element $\boldsymbol{\zeta}$ of \m{H^1(X,T_X\ot L)} associated 
to \m{X_2} (cf. \ref{ext-I}) in the sense of \ref{cech} (i.e. the cocycle 
relations are \ \m{D_{ik}=D_{ij}+\alpha^{(0)}_{ij}D_{jk}}). We have a canonical
product
\[H^1(X,T_X\ot L)\ot H^1(X,L^{n-1})\lra H^2(X,T_X\ot L^n)\]
and

\sepprop

\begin{subsub}\label{theo6}{\bf Theorem: } We have \ \m{\Delta''_n({\bf 
k}-{\bf h})=\boldsymbol{\zeta}\eta} .
\end{subsub}

\begin{proof}
To construct \m{\Delta''_n({\bf k}-{\bf h})} we start from the exact sequence
\xmat{0\ar[r] & T_X\ar[r] & \kg_{n+1}\ar[r]^-{\tau_n} & \kh_n\ar[r] & 0}
(cf. \cite{dr10}, 2.2, \cite{fr}). We take \ 
\m{\psi_{ij}\in\Aut(\ko_X(U_{ij})[t]/(t^n))} \ over 
\m{(\delta_{ij}^*,\theta_{ij})} (such that \ \m{\psi_{ji}=\psi_{ij}^{-1}}). 
Then \m{\psi_{ij}\psi_{jk}\psi_{ki}} is of the form
\[\psi_{ij}\psi_{jk}\psi_{ki} \ = \ I+\nu_{ijk}t^n \ , \]
where \m{\nu_{ijk}} is a derivation of \m{\ko_X(U_{ij})}. The family 
\m{(\nu_{ijk})} represents \m{\Delta''_n({\bf k})} in the sense of \ref{cech}, 
i.e. we have the cocycle relations \ \m{\nu_{ijk}-\nu_{ijl}+\nu_{ikl}-
\big(\alpha^{(0)}_{ij}\big)^n\nu_{jkl}=0} .

Similarly we take \ 
\m{\phi_{ij}\in\Aut(\ko_X(U_{ij})[t]/(t^n))} \ over 
\m{(\delta_{ij}^*,\alpha_{ij})} (such that \ \m{\phi_{ji}=\phi_{ij}^{-1}}). 
Then \m{\phi_{ij}\phi_{jk}\phi_{ki}} is of the form
\[\phi_{ij}\phi_{jk}\phi_{ki} \ = \ I+\tau_{ijk}t^n \ , \]
where \m{\tau_{ijk}} is a derivation of \m{\ko_X(U_{ij})}. The family 
\m{(\tau_{ijk})} represents \m{\Delta''_n({\bf h})} in the sense of \ref{cech}.

We can take
\[\psi_{ij}-\phi_{ij} \ = \ t^n\Delta_{ij} \ , \]
where, if we want that \m{\psi_{ij}} is a morphism of rings, 
\m{\Delta_{ij}:\ko_X(U_{ij})[t]/(t^{n+1})\to\ko_X(U_{ij})} \ is such that there 
exist a derivation \m{E_{ij}} of \m{\ko_X(U_{ij})} and \ 
\m{\mu_{ij}\in\ko_X(U_{ij})} \ such that for every \ 
\m{u\in\ko_X(U_{ij})[t]/(t^{n+1})} \ we have
\[\Delta_{ij}(u) \ = \ E_{ij}(u^{(0)})+u^{(1)}\mu_{ij} \ . \]
Since \m{\psi_{ij}} is over \m{(\delta_{ij}^*,\theta_{ij})} we have \ 
\m{\mu_{ij}=\beta_{ij}}. It is easy to see that if \ 
\m{\phi_{ji}=\phi_{ij}^{-1}}, the condition \ 
\m{\psi_{ji}=\psi_{ij}^{-1}} \ is equivalent to 
\begin{equation}\label{equ16}\beta_{ji} \ = \ 
-\frac{\beta_{ij}}{(\alpha_{ij}^{(0)})^{n+1}}\end{equation}
which is already fulfilled by $(\ref{equ14})$ and
\begin{equation}\label{equ15} E_{ji} \ = \ 
\frac{\beta_{ij}}{(\alpha_{ij}^{(0)})^{n+1}}D_{ij}-
\frac{1}{(\alpha_{ij}^{(0)})^n}E_{ij} \ . \end{equation}
We have then
\[\nu_{ijk}-\tau_{ijk} \ = \ 
t^n\Delta_{ij}\phi_{jk}\phi_{ki}+\phi_{ij}(t^n\Delta_{jk})\phi_{ki}+
\phi_{ij}\phi_{jk}(t^n\Delta_{ki}) \ . \]
We have, for every \ \m{u\in\ko_X(U_{ij})[t]/(t^{n+1})}
\begin{eqnarray*}
t^n\Delta_{ij}\phi_{jk}\phi_{ki}(u) & = & 
\big(E_{ij}(u^{(0)})+(D_{ji}(u^{(0)})+\alpha^{(0)}_{ji}u^{(1)})\beta_{ij}
\big)t^n \ ,\\
\phi_{ij}(t^n\Delta_{jk})\phi_{ki}(u) & = & 
\big(E_{jk}(u^{(0)})+(D_{ki}(u^{(0)})+
\alpha^{(0)}_{ki}u^{(1)})\beta_{jk}\big)\big(\alpha^{(0)}_{ij}\big)^nt^n \ , \\
\phi_{ij}\phi_{jk}(t^n\Delta_{ki})(u) & = & \big(\alpha^{(0)}_{ik}\big)^n\big(
E_{ki}(u^{(0)})+\beta_{ki}u^{(1)}\big)t^n \ . \\
\end{eqnarray*}
Hence we have
\begin{eqnarray*}
\nu_{ijk}(u)-\tau_{ijk}(u) & = & 
\big(E_{ij}+\big(\alpha^{(0)}_{ij}\big)^nE_{jk}+
\big(\alpha^{(0)}_{ik}\big)^nE_{ki}\big)(u^{(0)})+\\
& & \big(\beta_{ij}D_{ji}+\big(\alpha^{(0)}_{ij}\big)^n\beta_{jk}D_{ki}
\big)(u^{(0)})+\\
& & \big(\alpha^{(0)}_{ji}\beta_{ij}+
\big(\alpha^{(0)}_{ij}\big)^n\alpha^{(0)}_{ki}\beta_{jk}+
\big(\alpha^{(0)}_{ik}\big)^n\beta_{ki}\big)u^{(1)} \ . \\
\end{eqnarray*}
By $(\ref{equ14})$ and $(\ref{equ16})$ we have \ \m{\alpha^{(0)}_{ji}\beta_{ij}+
\big(\alpha^{(0)}_{ij}\big)^n\alpha^{(0)}_{ki}\beta_{jk}+
\big(\alpha^{(0)}_{ik}\big)^n\beta_{ki}=0}. On the other hand, by \ref{cech2f} 
the family \m{\big(\big(\alpha^{(0)}_{ij}\big)^nE_{jk}-E_{ik}+E_{ij}\big)} is a 
coboundary, and by $(\ref{equ15})$ we have \Nligne \m{\dsp E_{ki}= 
\frac{\beta_{ik}}{(\alpha_{ik}^{(0)})^{n+1}}D_{ik}-
\frac{1}{(\alpha_{ik}^{(0)})^n}E_{ik}}. Hence we can assume that
\begin{eqnarray*}\nu_{ijk} & = & \frac{\beta_{ik}}{\alpha^{(0)}_{ik}}D_{ik}
+\beta_{ij}D_{ji}+\big(\alpha^{(0)}_{ij}\big)^n\beta_{jk}D_{ki}\\
& = & \Big(\frac{\beta_{ik}}{\alpha^{(0)}_{ik}}-
\frac{\big(\alpha^{(0)}_{ij}\big)^n}{\alpha^{(0)}_{ik}}\beta_{jk}\Big)D_{ik}
+\beta_{ij}D_{ji} \ . \\ \end{eqnarray*}
But \ \m{\dsp\frac{\beta_{ik}}{\alpha^{(0)}_{ik}}-
\frac{\big(\alpha^{(0)}_{ij}\big)^n}{\alpha^{(0)}_{ik}}\beta_{jk}=
\frac{\beta_{ij}}{\alpha^{(0)}_{ij}}} \ by $(\ref{equ14})$, and \ \m{D_{ik}=
D_{ij}+\alpha^{(0)}_{ij}D_{jk}}, \m{\dsp D_{ji}=-\frac{1}{\alpha^{(0)}_{ij}}
D_{ij}}. Hence
\[\nu_{ijk} \ = \ \beta_{ij}D_{jk} \ = \ \alpha^{(0)}_{ij}\Big(
\frac{\beta_{ij}}{\alpha^{(0)}_{ij}}\Big)D_{jk} \ . \]
Since \m{\dsp\Big(\frac{\beta_{ij}}{\alpha^{(0)}_{ij}}\Big)} represents $\eta$ 
and \m{(D_{ij})} represents $\boldsymbol{\zeta}$, the result follow from 
\ref{cech2e}.
\end{proof}

\end{sub}

\sepsub

\Ssect{Extensions of vector bundles to higher multiplicity}{ext_B}

Suppose that \m{X_n} can be extended to a primitive multiple scheme \m{X_{n+1}} 
of multiplicity \m{n+1}. Let $\E$ be a vector bundle on \m{X_n}, and \ 
\m{E=\E_{|X}}. If $\E$ can be extended to a vector bundle \m{\E_{n+1}} on 
\m{X_{n+1}}, then we have an exact sequence
\[0\lra E\ot L^n\lra\E_{n+1}\lra\E=\E_{n+1|X_n}\lra 0 \ . \]

\sepprop

\begin{subsub}\label{ext_B1} Obstruction to the extension of a vector bundle in 
higher multiplicity -- \rm In \cite{dr10}, 7.1, a class \ 
\m{\Delta(\E)\in\Ext^2_{\ko_X}(E,E\ot L^n)} \ is defined, such that $\E$ can be 
extended to a vector bundle on \m{X_{n+1}} if and only \ \m{\Delta(\E)=0}.

The canonical exact sequence \ \m{0\to  
L^n\to\Omega_{X_{n+1}|X_n}\to\Omega_{X_n}\to 0} \ induces \Nligne 
\m{\ov{\sigma}_{\E,X_{n+1}}\in\Ext^1_{\ko_{X_n}}(\E\ot\Omega_{X_n},\E\ot L^n)}. 
We have a canonical product
\[\Ext^1_{\ko_{X_n}}(\E\ot\Omega_{X_n},\E\ot L^n)\times
\Ext^1_{\ko_{X_n}}(\E,\E\ot\Omega_{X_n})\lra
\Ext^2_{\ko_{X_n}}(\E,\E\ot L^n)= \Ext^2_{\ko_X}(E,E\ot L^n) \ , \]
and
\begin{equation}\label{equ3}
\Delta(\E)=\ov{\sigma}_{\E,X_{n+1}}\nabla_0(\E) \ ,
\end{equation}

where \m{\nabla_0(\E)} is the {\em canonical class} of $\E$ (cf. \ref{CCVB}, 
\cite{dr10}, Theorem 7.1.2).
\end{subsub}

\sepprop

\begin{subsub}\label{lem4}{\bf Lemma: } We have \
\m{\EExt^1_{\ko_{X_{n+1}}}(\E,E\ot L^n)\simeq\EEnd(E)}.
\end{subsub}
\begin{proof} Let \m{U\subset X} be a nonempty subset and \m{U_n}, \m{U_{n+1}} 
the corresponding open subsets of \m{X_n}, \m{X_{n+1}} respectively. 
Suppose that $U$ is affine. It follows from \ref{ext_B} that \m{\E_{|U_{n}}} 
(resp. \m{\ki_{U,U_{n+1}}}) can be extended to a vector bundle $\F$ (resp. 
$\L$) on \m{U_{n+1}}. We have a canonical morphism \m{\L\to\ko_{U_{n+1}}}:
\xmat{\L\ar[r] & \L_{|U_n}=\ki_{U,U_{n+1}}\flinc[r] & \ko_{U_{n+1}} \ ,}
inducing \ \m{\alpha:\L^{n+1}\to\L^n} \ and \ \m{\beta:\L^n\to\ko_{n+1}}. We 
have then an obvious locally free resolution of \m{\E_{|U_{n+1}}}:
\xmat{\cdots\ar[r] & \F\ot\L^{n+1} \ar[rr]^-{I_F\ot\alpha} & & \F\ot\L^n 
\ar[rr]^-{I_F\ot\beta} & & \ar[r] \F\ar[r] & \E_{|U_{n}}\ar[r] & 0 \ . }
Using this resolution it follows that \ 
\m{\EExt^1_{\ko_{U_{n+1}}}(\E_{U_{n+1}},E_{|U}\ot L^n)\simeq\EEnd(E_{|U})}. It 
is easy to see that these isomorphisms can be glued together to define the 
isomorphism of Lemma \ref{lem4}.
\end{proof}

\sepprop

\begin{subsub}\label{ext_4} Other description of \m{\Delta(\E)} and of the 
extensions -- \rm
From \ref{Ex_Ext} and Lemma \ref{lem4}, we have the exact sequence 
\m{\bm{\Gamma}_{\E,E\ot L^n}} on \m{X_{n+1}}
\xmat{0\ar[r] & \Ext^1_{\ko_X}(E,E\ot L^n)\ar[r]^-\psi & 
\Ext^1_{\ko_{X_{n+1}}}(\E,
E\ot L^n)\ar[r]^-\phi & \End(E)\ar[r]^-{\delta_\E} & \Ext^2_{\ko_X}(E,E\ot L^n) 
\ , }
and we have \ \m{\Delta(\E)=\delta_\E(I_E)} (cf. \cite{dr10}).
\end{subsub}

\sepsubsub

\begin{subsub}\label{ext_5} Other description of the extensions -- \rm Let
\[0\lra E\ot L^n\lra\ke\to\E\lra 0\]
be an extension on \m{X_{n+1}}, with $\ke$ locally free, and 
\m{\beta\in\Ext^1_{\ko_X}(E,E\ot L^n)}. There is an open cover \m{(U_i)_{i\in 
I}} of $X$ such that there are trivializations
\[\theta_i:\ke_{|U_i^{(n+1)}}\lra\ko_{U_i^{(n+1)}}\ot\C^r \ \]
so that $\ke$ is represented by the family \m{(\theta_{ij})}, 
\m{\theta_{ij}\in\Aut(\ko_{U_{ij}^{(n+1)}}\ot\C^r)}.
Recall that \ \m{\ki_{X_n,X_{n+1}}=L^n}, viewed as a sheaf on \m{X_{n+1}}. Let 
\ \m{\beta\in\Ext^1_{\ko_X}(E,E\ot L^n)} . Then (for a suitable cover 
\m{(U_i)}) $\beta$ is represented by \m{(\beta_{ij})}, with
\[\beta_{ij} \ \in \ 
H^0(\ko_{U_{ij}}\ot\End(\C^r)\ot L^n) \ \subset \ 
H^0(\ko_{U_{ij}^{(n+1)}}\ot\End(\C^r)) \ , \]
in the sense of \ref{cech1}, i.e. 
the cocycle relation is \ 
\m{\beta_{ik}=\beta_{ij}+\theta_{ij}\beta_{jk}\theta_{ij}^{-1}}. Let
\[0\lra E\ot L^n\lra\ke'\lra\E\lra 0\]
be the extension, corresponding to \m{\sigma+\Psi(\beta)}. Then \m{\ke'} 
is represented by the family \m{((1+\beta_{ij})\theta_{ij})}. This description 
is similar to that of \ref{cech3}.
\end{subsub}

\sepsubsub

If \ \m{\Delta(\E)=0} then $\phi$ is surjective (because 
\m{\imm(\phi)} contains all the automorphisms of $E$).

This is also a consequence of

\sepprop

\begin{subsub}\label{lem5}{\bf Lemma: } Let
\ \m{0\to E\ot L^n\to\ke\to\E\lra 0} \ be an extension, corresponding to \ 
\m{\sigma\in\Ext^1_{\ko_{X_{n+1}}}(\E,E\ot L^n)}. Then the sheaf $\ke$ is 
locally free if and only if \m{\phi(\sigma)} is an isomorphism.
\end{subsub}
\begin{proof} We need only to prove the analogous following result: let $U$ be 
a nonempty open\break subset of $X$, \m{U_n}, \m{U_{n+1}} the corresponding 
open  subsets of \m{X_n}, \m{X_{n+1}} respectively, such that\break 
\m{\ki_{U,U_{n+1}}} can 
be extended to a line bundle $\L$ on \m{U_{n+1}}. Let $r$ be a positive integer 
and \break \m{0\to L^n_{|U}\ot\C^r\to\ke\to\ko_{U_n}\ot\C^r\to 0} \ an exact 
sequence of sheaves on \m{U_{n+1}}, and \Nligne 
\m{\phi\in\End(L^n_{|U}\ot\C^r)} \ 
the associated morphism. Then we have \ \m{\ko\simeq\ko_{U_{n+1}}\ot\C^r} \ if 
and only $\phi$ is an isomorphism. 

The result follows easily from the following fact, using the locally free 
resolution of Lemma \ref{lem4}: $\ke$ is isomorphic to the cokernel of the 
morphism
\xmat{L^n_{|U}\ot\C^r\ar[rr]^-{(\phi,i)} & & 
(L^n_{|U}\ot\C^r)\oplus(\ko_{U_{n+1}}\ot\C^r) \ ,}
where $i$ is the canonical inclusion (cf. \cite{dr1b}, 4.2).
\end{proof}

\sepprop

It follows that when $\ke$ is locally free we can always assume that \ 
\m{\phi(\sigma)=I_E}.

\sepsubsub

\begin{subsub}\label{isom_cl} Isomorphism classes -- \rm We have an exact 
sequence \ \m{0\to L^n\to\ko_{X_{n+1}}\to\ko_{X_n}\to 0}, whence a connecting 
homomorphism
\[\delta^0:H^0(X_n,\ko_{X_n})\lra H^1(X,L^n) \ . \]
\end{subsub}
Since \ \m{E^*\ot E\ot L^n\simeq L^n\oplus(\AAd(E)\ot L^n)}, we can see 
\m{\imm(\delta^0)} as a subspace of\Nligne \m{\Ext^1_{\ko_X}(E,E\ot L^n)}.

Recall that $\E$ is called {\em simple} if the only endomorphisms of $\E$ are 
the homotheties, i.e. if the canonical morphism \ \m{H^0(\ko_{X_n})\to\End(\E)} 
\ is an isomorphism. 

\sepprop

\begin{subsub}\label{prop15}{\bf Proposition: } Suppose that $\E$ is simple. 
Let
\[0\lra E\ot L^n\lra\ke\to\E\lra 0\]
be an extension on \m{X_{n+1}}, with $\ke$ locally free, corresponding to \ 
\m{\sigma\in\Ext^1_{\ko_{X_{n+1}}}(\E,E\ot L^n)}, and
\[0\lra E\ot L^n\lra\ke'\lra\E\lra 0\]
another extension, corresponding to \m{\sigma+\Psi(\beta)}, 
\m{\beta\in\Ext^1_{\ko_X}(E,E\ot L^n)}. Then \m{\ke'\simeq\ke} if and only if \ 
\m{\beta\in\imm(\delta^0)}.
\end{subsub}

\begin{proof}
We use the notations of \ref{ext_4}. Let \ 
\m{\theta'_i:\ke'_{|U_i^{(n+1)}}\to\ko_{U_i^{(n+1)}}\ot\C^r} \ be 
trivializations, with \ \m{\theta'_{ij}=(1+\beta_{ij})\theta_{ij}}.

There is an automorphism \ \m{f:\ke\to\ke'} \ if and only for every \m{i\in I} 
there exists an automorphism \m{f_i} of \m{\ko_{U_i^{(n+1)}}\ot\C^r} such that 
the following square is commutative
\xmat{\ke_{|U_i^{(n+1)}}\ar[r]^-{\theta_i}\ar[d]^f & 
\ko_{U_i^{(n+1)}}\ot\C^r\ar[d]^{f_i} \\
\ke'_{|U_i^{(n+1)}}\ar[r]^-{\theta'_i} & \ko_{U_i^{(n+1)}}\ot\C^r}
Suppose that $f$ exists. Then since $\E$ is simple the morphism \ 
\m{\ke_{|X_n}=\E\to\ke'_{|X_n}=\E} \ induced by $f$ is the multiplication by \ 
\m{\mu\in H^0(\ko_{X_n})^*}. Hence we can write \ \m{f_i=\mu_i.I+\eta_i}, where
\ \m{\mu_i\in H^0(\ko_{X_{n+1}}(U_i^{(n+1)}))} \ is such that \ 
\m{\mu_{i|U_i^{(n)}}=\mu} \ and \ 
\m{\eta_i:\ko_{U_i^{(n+1)}}\ot\C^r\to\ki_{X_n,X_{n+1}|U_i^{(n+1)}}}. Let \ 
\m{\mu_0=\mu_{|X}\in\C^*}. We have on \m{U_{ij}^{(n+1)}}, 
\m{{\theta'}_i^{-1}f_i\theta_i={\theta'}_j^{-1}f_j\theta_j=f}, i.e. 
\m{f_i\theta_{ij}=\theta'_{ij}f_j}, which gives
\[\beta_{ij} \ = \ 
\frac{\mu_i-\mu_j}{\mu_j}+\frac{\eta_i-\theta_{ij}\eta_j\theta_{ij}^{-1}}{\mu_0}
 \ .\]
By \ref{cech2f}, the family \m{(\eta_i-\theta_{ij}\eta_j\theta_{ij}^{-1})} is a 
coboundary, hence we can suppose that \ \m{\dsp\beta_{ij}=
\frac{\mu_i-\mu_j}{\mu_j}}. Then, since \ \m{\mu_i-\mu_j\in 
H^0(U_{ij}^{(n+1)},\ki_{X_n,X_{n+1}})}, we have \ 
\m{\dsp\frac{\mu_i-\mu_j}{\mu_j}=\frac{\mu_i}{\mu_0}-\frac{\mu_j}{\mu_0}}. Hence
\ \m{\dsp\beta=\delta^0\Big(\frac{\mu}{\mu_0}\Big)}. The converse is proved in 
the same way.
\end{proof}

\end{sub}

\sepsub

\Ssect{Canonical filtrations, quasi locally free sheaves and Serre duality}{QF}

Let \m{P\in X} be a closed point, \m{z\in\ko_{Y,P}} an equation of $X$ and $M$ 
a \m{\ko_{X_n,P}}-module of finite type. Let $\ke$ be a coherent sheaf on 
\m{X_n}.

The two {\em canonical filtrations} are useful tools to study the coherent
sheaves on primitive multiple curves.

\sepprop

\begin{subsub}\label{QLL-def1} First canonical filtration -- \rm
The first canonical filtration of $M$ is
\[M_n=\nsp\subset M_{n-1}\subset\cdots\subset M_{1}\subset M_0=M\]
where for \m{0\leq i< n}, \m{M_{i+1}} is the kernel of the surjective
canonical morphism\Nligne \m{M_{i}\to M_{i}\ot_{\ko_{n,P}}\ko_{X,P}} . So we
have
\[M_{i}/M_{i+1} \ = \ M_{i}\ot_{\ko_{n,P}}\ko_{X,P}, \ \ \ \
M/M_i \ \simeq \ M\ot_{\ko_{n,P}}\ko_{X_i,P}, \ \ \ \
M_i \ = \ z^iM .\]

One defines similarly the {\em first canonical filtration of $\ke$}: it is the
filtration
\[\ke_n=0\subset \ke_{n-1}\subset\cdots\subset \ke_{1}\subset \ke_0=\ke\]
such that for \m{0\leq i< n}, \m{\ke_{i+1}} is the kernel of the canonical
surjective morphism \ \m{\ke_i\to\ke_{i\mid X}}. So we have \
\m{\ke_{i}/\ke_{i+1}=\ke_{i\mid X}} ,
\m{\ke/\ke_i=\ke_{\mid X_i}} . 
\end{subsub}

\sepsubsub

\begin{subsub}\label{2-fc} Second canonical filtration -- \rm
One defines similarly the {\em second canonical filtration of $M$}: it is the
filtration
\[M^{(0)}=\nsp\subset M^{(1)}\subset\cdots\subset M^{(n-1)}\subset M^{(n)}=M\]
with \ \m{M^{(i)} \ = \ \big\lbrace u\in M ; z^iu=0\big\rbrace} . If \
\m{M_n=\nsp\subset M_{n-1}\subset\cdots\subset M_1\subset M_0=M} \ is
the first canonical filtration of $M$, we have \ \m{M_i\subset M^{(n-i)}} \
for \m{0\leq i\leq n}. 

One defines in the same way the {\em second canonical filtration of $\ke$}:
\[\ke^{(0)}=\nsp\subset \ke^{(1)}\subset\cdots\subset
\ke^{(n-1)}\subset \ke^{(n)}=\ke  .\]
\end{subsub}

Let $P$ be a closed point of $X$.
Let $M$ be a \m{\ko_{X_n,P}}-module of finite type. Then $M$ is called {\em
quasi free} if there exist non negative integers \m{m_1,\ldots,m_n} and an
isomorphism \ \m{M\simeq\oplus_{i=1}^nm_i\ko_{X_i,P}}~. The integers
\m{m_1,\ldots,m_n} are uniquely determined: it is easy to recover them from 
the first canonical filtration of $M$. We say that \m{(m_1,\ldots,m_n)} is the 
{\em type} of $M$.

Let $\ke$ be a coherent sheaf on \m{X_n}. We say that $\ke$ is {\em quasi free
at $P$} if \m{\ke_P} is quasi free, and that $\ke$ is {\em quasi locally free} 
on a nonempty open subset \m{U\subset X} if it is quasi free at every point of 
$U$. If \m{U=X} we say that $\ke$ is quasi locally free. In this case the types 
of the modules \m{\ke_x}, \m{x\in X}, are the same, and for every \m{x\in X} 
there exists a neighbourhood \m{V\subset X_n} of $x$ and an isomorphism
\begin{equation}\label{equ20}\ke_{|V} \ \simeq \ \sigg_{n=1}^n m_i\ko_{X_i\cap 
V} \ .
\end{equation}

The proof of the following theorem is the same as in the case \m{\dim(X)=1} 
(th. 5.1.3 of \cite{dr2}).

\sepprop

\begin{subsub}\label{theo1}{\bf Theorem: } The following two assertions are 
equivalent:
\begin{enumerate}
\item[(i)] The $\ko_{X_n,P}$-module $M$ is quasi free.
\item[(ii)] All the $M_i/M_{i+1}$ are free $\ko_{X,P}$-modules.
\end{enumerate}
\end{subsub}

\sepprop

It follows that 

\sepprop

\begin{subsub}\label{theo2}{\bf Theorem: } The following two assertions are 
equivalent:
\begin{enumerate}
\item[(i)] $\ke$ is quasi locally free.
\item[(ii)] All the \m{\ke_i/\ke_{i+1}} are locally free on $X$.
\end{enumerate}
\end{subsub}

\sepprop

For every coherent sheaf $\ke$ on \m{X_n} there exists a nonempty open subset 
\m{U\subset X} such that \m{\ke} is quasi locally free on $U$.

\sepprop

\begin{subsub}\label{prop10}{\bf Proposition: } Let $\ke$ be a quasi locally 
free sheaf on \m{X_n}. Then 

{\rm (i) } $\ke$ is reflexive, 

{\rm (ii)} for every positive integer 
$i$ we have \ \m{\EExt^i_{\ko_{X_n}}(\ke,\ko_{X_n})=0}.

{\rm (iii)} for every vector bundle $\F$ on \m{X_n} and every integer $i$ we 
have \Nligne \m{\Ext^i_{\ko_{X_n}}(\ke,\F)\simeq H^i(X_n,\ke^\vee\ot\F)} .
\end{subsub}

\begin{proof}
The assertion (iii) follows easily from (i), with the Ext spectral sequence. 
Assertions (i) and (ii) are local, so we can replace \m{X_n} with \ \m{U\times 
{\bf Z}_n}, where $U$ is an open subset of \m{X_n}. Let \m{U_i=U\times{\bf 
Z}_i}, for \m{1\leq i\leq n}. The results follow from $(\ref{equ20})$ and the 
free resolutions
\xmat{....\ar[r] & \ko_{U_n}\ar[r]^-{\times t^{i}}\ar[r] & 
\ko_{U_n}\ar[r]^-{\times t^{n-i}}\ar[r] & \ko_{U_n}\ar[r]^-{\times t^{i}}
\ar[r] & \ko_{U_i}}
\end{proof}

\sepsubsub

\begin{subsub}\label{ser_dual} Serre duality for quasi locally free sheaves -- 
\rm Suppose that $X$ is projective of dimension $d$. It follows from 
Proposition \ref{prop10} that
\end{subsub}

\sepprop

\begin{subsub}\label{prop11}{\bf Proposition: } For every \m{i\geq 0} and every 
quasi locally free sheaf $\ke$ on \m{X_n}, there is a functorial isomorphism
\[H^i(X_n,\ke) \ \simeq \ H^{d-i}(X_n,\ke^\vee\ot\omega_{X_n})^* \ . \]
\end{subsub}

\sepprop

\begin{subsub}\label{coro3}{\bf Corollary: } We have \ 
\m{\omega_{X_n|X}\simeq\omega_X\ot L^{1-n}} .
\end{subsub}

The proof is similar to that of Proposition 7.2, III, of \cite{ha}.

\end{sub}

\sepsub

\Ssect{Simple vector bundles}{vb}

We use the notations of \ref{PMS_def}.

\sepprop

\begin{subsub}\label{endo0}{\bf Lemma:} Suppose that the restriction map \ 
\m{H^0(\ko_{X_n})\to H^0(\ko_{X_{n-1}})} \ is surjective, that 
\m{\E_{|X_{n-1}}} is simple on \m{X_{n-1}} and \ \m{\Hom(E,E\ot 
L^{n-1})=H^0(X,L^{n-1})}. Then $\E$ is simple.
\end{subsub}
\begin{proof} Let \m{f\in\End(\E)}. Then the induced morphism \ 
\m{\E_{|X_{n-1}}\to\E_{|X_{n-1}}} \ is the multiplication by some \ 
\m{\lambda\in H^0(\ko_{X_{n-1}})}. Let \ \m{\mu\in H^0(\ko_{X_n})} \ be such 
that \ \m{\mu_{|X_{n-1}}=\lambda}. Then the restriction of \ \m{f-\mu I_{X_n}} 
\ to \m{X_{n-1}} vanishes. It follows that \ \m{f-\mu I_{X_n}} \ can be 
factorized in the following way
\xmat{\E\flon[r] & \E_{|X}=E\ar[r] & E\ot L^{n-1}=\E^{(1)}\flinc[r]& \E}
Since \ \m{\Hom(E,E\ot L^{n-1})=H^0(X,L^{n-1})}, \m{f-\mu I_{X_n}} is the 
multiplication by some \Nligne \m{\eta\in H^0(X,L^{n-1})\subset H^0(\ko_{X_n})},
and \m{f=(\mu+\eta)I_{X_n}}.
\end{proof}

\sepprop

It follows easily by induction on $n$ that

\sepprop

\begin{subsub}\label{endo}{\bf Corollary:} Suppose that $L$ is a non trivial 
ideal sheaf, or \m{\deg(L)<0} if $X$ is a curve, and that \m{E=\E_{|X}} is 
simple. Then $\E$ is simple.
\end{subsub}

\sepprop

\end{sub}

\newpage

\section{Canonical class associated to a vector bundle}\label{CCVB}

\Ssect{Definition}{def_can}

Let $Z$ be a scheme over $\C$ 
and $E$ a vector bundle on $Z$ 
of rank $r$. To $E$ one associates an element \m{\nabla_0(E)} of \m{H^1(Z,E\ot 
E^*\ot\Omega_Z)}, called the {\em canonical class of $E$} (cf. \cite{dr10}, 
3-). If $L$ is a line bundle on $Z$, then \ \m{\nabla_0(L)\in H^1(Z,\Omega_Z)}.
If $Z$ is smooth and projective, \m{L=\ko_Z(Y)}, where \m{Y\subset Z} is a 
smooth hypersurface, then \m{\nabla_0(L)} is the cohomology class of $Y$.

Let \m{(Z_i)_{i\in I}} be an open cover of $Z$ such that $E$ is defined by a 
cocycle \m{(\theta_{ij})}, \m{\theta_{ij}\in\GL(r,\ko_{Z_{ij}})}. Then 
\m{((d\theta_{ij})\theta_{ij}^{-1})} is a cocycle (in the sense of 
\ref{cech2c}) which represents \m{\nabla_0(E)}.
 
We have \ \m{E\ot E^*\simeq\AAd(E)\oplus\ko_Z}, so \ \m{H^1(E\ot 
E^*\ot\Omega_Z)\simeq H^1(Z,\AAd(E)\ot\Omega_Z)\oplus H^1(Z,\Omega_Z)}. The 
component of \m{\nabla_0(E)} in \m{H^1(Z,\Omega_Z)} is 
\m{\dsp\frac{1}{r}\nabla_0(\det(E))}.

If \m{L_1}, \m{L_2} are line bundles on $Z$, then
\begin{equation}\label{equ21}\nabla_0(L_1\ot L_2) \ = \ 
\nabla_0(L_1)+\nabla_0(L_2) \ .
\end{equation}

Let $S$ be a scheme over $\C$ and \m{s\in S} a closed point. Let $Z$ be a 
scheme and \m{p_Z:Z\times S\to Z}, \m{p_S:Z\times S\to S} the projections. Let 
$\E$ be a vector bundle over \m{Z\times S}. We have
\[\nabla_0(\E) \ \in \ H^1(Z\times S,\E\ot\E^*\ot\Omega_{Z\times S})=
H^1(Z\times S,\E\ot\E^*\ot p_Z^*(\Omega_Z))\oplus H^1(Z\times S,\E\ot\E^*\ot 
p_S^*(\Omega_S)) \ . \]
By restricting to \m{Z\times\{s\}} we get a canonical map
\[r_s:H^1(Z\times S,\E\ot\E^*\ot\Omega_{Z\times S}) \lra
H^1(Z,\E_s\ot\E^*_s\ot\Omega_Z) \ . \]
From Lemma 3.2.1 of \cite{dr10} we have
\begin{equation}\label{equ2}
r_s(\nabla_0(\E)) \ = \ \nabla_0(\E_s) \ .
\end{equation}
\end{sub}

\sepsub

\Ssect{The case of smooth projective varieties}{can_smo}

Let $Z$ be a scheme over $\C$, $T$ a projective integral variety and 
$\kl$ a line bundle on \m{Z\times T}, viewed as a family of line bundles on $Z$ 
parametrized by $T$. For any closed point \m{t\in T}, let \ 
\m{\kl_t=\kl_{|Z\times\{t\}}\in\Pic(Z)}.

\sepprop

\begin{subsub}\label{prop7}{\bf Proposition: } The map
\[\xymatrix@R=5pt{\nabla_\kl:T\ar[r] & H^1(Z,\Omega_Z)\\ t\fmaps[r] & 
\nabla_0(\kl_t)}\] 
is constant.
\end{subsub}

\begin{proof}
Let \m{t\in T} and \m{S\subset T} an affine open neighbourhood of $t$. Then we 
have \Nligne \m{H^1(Z\times S,\Omega_{Z\times S})\simeq(\ko_S(S)\ot
H^1(Z,\Omega_Z)\oplus(H^0(S,\Omega_S)\ot H^1(T,\ko_T))}. Let
\[\nabla_0(\kl_{|Z\times S}) \ = \ \sigg_{k=1}^p f_k\ot\nabla_k+\eta \ , \]
with \ \m{f_1,\ldots,f_p\in\ko_S(S)}, \m{\nabla_1,\ldots,\nabla_p\in 
H^1(Z,\Omega_Z)} and \ \m{\eta\in H^0(S,\Omega_S)\ot H^1(T,\ko_T))}. Then by 
\ref{def_can} we have \ \m{\nabla_\kl(s)=\sigg_{k=1}^p 
f_k(s)\nabla_k} \ for every \m{s\in S}. Hence \m{\nabla_\kl} is a regular map 
to a finite dimensional vector space. Since $T$ is projective, \m{\nabla_\kl} 
must be constant.
\end{proof}

\sepprop

\begin{subsub}\label{coro2}{\bf Corollary: } If $Z$ is smooth and projective,
and \ \m{L\in\Pic(Z)}, then \m{\nabla_0(L)} is invariant by deformation of $L$.
\end{subsub}

\begin{proof} This follows from the fact that the connected component of 
\m{\Pic(Z)} containing $L$ parametrizes a projective family of line bundles 
containing $L$. 
\end{proof}

\sepprop

\begin{subsub}\label{difmo} The differential morphism -- \rm Let \ \m{\eta\in 
H^1(Z,\ko_Z)}, represented by a cocycle \m{(\eta_{ij})} (with respect to an 
open cover \m{(Z_i)_{i\in I}} of $Z$), \m{\eta_{ij}\in\ko_Z(Z_{ij})}. Then 
\m{(d\eta_{ij})} is a cocycle which represents \ \m{d_Z(\eta)\in 
H^1(\Z,\Omega_Z)} (which depends only on $\eta$). This defines a linear map
\[d_Z:H^1(Z,\ko_Z)\lra H^1(Z,\Omega_Z) \ . \]
\end{subsub}

\sepprop

The following result and its proof are classical:

\sepprop

\begin{subsub}\label{prop8}{\bf Proposition:} If $Z$ is smooth and projective, 
then \ \m{d_Z=0} .
\end{subsub}

\begin{proof}
Let $\kl$ be a Poincaré bundle on \ \m{\Pic^0(Z)\times Z}. We can compute the 
canonical classes in analytic cohomology. In particular, there exist an open 
cover \m{(Z_i)} of $Z$ and a neighbourhood $S$ of \ \m{\ko_Z\in\Pic^0(Z)} (for 
the usual topology) such that \m{\kl_{|Z\times S}} is represented by a cocycle 
\m{(\theta_{ij})}, where \m{\theta_{ij}} is an invertible analytic function on 
\ \m{Z_{ij}\times S}. Let \ \m{\BS{u}\in H^1(Z,\ko_Z)=T_{\ko_Z}(\Pic^0(Z))}. 
There exists an analytic map \ \m{\iota:U\to S} (where $U$ is a neighbourhood 
of 0 in $\C$), which is an embedding such that \ \m{\phi'(0):\C\to 
T_{\ko_Z}(\Pic^0(Z))} \ sends 1 to $\BS{u}$. Let 
\[\alpha_{ij}^U \ = \ \alpha_{ij|Z_{ij}\times U} \ , \qquad
\alpha_{ij}^0 \ = \ \alpha_{ij|Z_{ij}\times\{0\}} \ . \]
we can write
\[\alpha_{ij}^U(x,t) \ = \ \alpha_{ij}^0(x)+t\alpha_{ij}^0(x)\beta_{ij}(x,t) \ ,
\]
where \m{\beta_{ij}} is an analytic map. Let \ \m{\beta_{ij}^0 \ = \ 
\beta_{ij|Z_{ij}\times\{0\}}}. We have \ \m{\beta_{ij}^0 + \beta_{jk}^0 
=\beta_{ik}^0} \ on \m{Z_{ijk}}, and the family \m{(\beta_{ij}^0)} is a cocycle 
which represents $\BS{u}$.

A simple computation shows that we can write
\[\frac{d\alpha_{ij}^U}{\alpha_{ij}^U} \ = \ 
\frac{d\alpha_{ij}^0}{\alpha_{ij}^0}+td\beta_{ij}+\rho_{ij}dt+t^2\mu_{ij} \ , \]
where \m{\rho_{ij}} is an analytic function on \ \m{Z_{ij}\times U} \ and \ 
\m{\mu_{ij}\in H^0(Z_{ij}\times U,\Omega_{Z_{ij}\times U})}. It follows that 
\[d(\BS{u}) \ = \ (\nabla_\kl\circ\iota)'(1) \ . \]
Since \m{\nabla_\kl\circ\iota} is constant by Proposition \ref{prop7}, we have 
\m{d_Z(\BS{u})=0}.
\end{proof}

\sepprop

Similarly, in the general case, let $Y$ be a scheme and $\kl$ a line bundle on 
\m{Z\times Y}. Let \m{y\in Y} be a closed point and \ \m{\omega_y:T_yY\to 
H^1(Z,\ko_Z)} \ the Koda\"\i ra-Spencer map. Then we have

\sepprop

\begin{subsub}\label{prop8b}{\bf Proposition:} We have \ 
\m{T\nabla_{\kl,y}=d_Z\circ\omega_y}.
\end{subsub}

\end{sub}

\sepsub

\Ssect{The dualizing sheaf of a primitive double scheme}{dua_X2}

Let $X$ be a smooth, projective and irreducible variety. Let \m{X_2} be a 
primitive double scheme, with underlying smooth variety $X$, and associated 
line bundle $L$ on $X$. We have a canonical exact sequence
\xmat{0\ar[r] & L\ar[r] & \ko_{X_2}\ar[r]^-\zeta & \ko_X\ar[r] & 0}
associated with \ \m{\eta\in H^1(X,T_X\ot L)} (cf. \ref{ext-I}).
We will prove

\sepprop

\begin{subsub}\label{theo7}{\bf Theorem: } The map \ 
\m{H^1(\zeta):H^1(X_2,\ko_{X_2})\to H^1(X,\ko_X)} \ is surjective.
\end{subsub}

\sepprop

The following result generalizes Proposition 1.5 of \cite{ga-go-pu}:

\sepprop

\begin{subsub}\label{coro4}{\bf Corollary:} We have \ 
\m{\omega_{X_2}\simeq\ko_{X_2}} \ if and only if \ \m{\omega_X\simeq L} .
\end{subsub}

\begin{proof}
If \ \m{\omega_{X_2}\simeq\ko_{X_2}} \ then \ \m{\omega_X\simeq L} \ by 
Corollary \ref{coro3}. Conversely, suppose that \ \m{\omega_X\simeq L}. Then \ 
\m{\omega_{X_2|X}\simeq\ko_X}. We have an exact sequence
\xmat{0\ar[r] & \omega_X=L\lra\omega_{X_2}\ar[r]^-\rho & \ko_X\ar[r] & 0 \ . }
and the associated long exact sequence
\xmat{0\ar[r] & H^0(X,\omega_X)\ar[r] & 
H^0(X_2,\omega_{X_2})\ar[rr]^-{H^0(\rho)} & &
H^0(X,\ko_X)\ar[r]^-\delta & H^1(X,L) \ .}
On the other hand, from the exact sequence
\[0\lra L\lra\ko_{X_2}\lra\ko_X\lra 0\]
we deduce the following
\xmat{H^1(X_2,\ko_{X_2})\ar[rr]^-{H^1(\xi)} & & H^1(X,\ko_X)\ar[r]^-{\delta'} & 
H^2(X,L) \ . }
From Serre duality on \m{X_2} we get the commutative square
\xmat{H^0(X,\ko_X)\ar[rr]^-\delta\ar[d]^\simeq & & H^1(X,L)\ar[d]^\simeq\\
H^2(X,L)^*\ar[rr]^-{^t\delta'} & & H^1(X,\ko_X)^*}
From Theorem \ref{theo7} we have \m{\delta'=0}. Hence \m{\delta=0} and 
\m{H^0(\rho)} is surjective. If \ \m{\sigma\in H^0(X_2,\omega_{X_2})} \ is such 
that \ \m{H^0(\rho)(\sigma)=1}, the corresponding morphism \ 
\m{\ko_{X_2}\to\omega_{X_2}} \ is an isomorphism, so \ 
\m{\omega_{X_2}\simeq\ko_{X_2}}.
\end{proof}

\sepprop

We now prove Theorem \ref{theo7}, i.e. that the connecting morphism \ 
\m{\delta:H^1(X,\ko_X)\to H^2(X,L)} \ is null. We use the notations of 
\ref{PMS_def}. Let \ \m{\beta\in H^1(X,\ko_X)}, defined by a cocycle 
\m{(\beta_{ij})}, \m{\beta_{ij}\in\ko_X(U_{ij})}. We can view \m{\beta_{ij}} as 
an element of \m{H^0(\ko_{U_{ij}\times Z_2})}. Then \ \m{\delta(\beta)} is 
represented by the cocycle \m{(\sigma_{ijk})}, with 
\[\sigma_{ijk} \ = \ \beta_{ij}+\delta_{ij}^*(\beta_{jk})-\beta_{ik}\]
(cf. \ref{con_hom}). We have then
\begin{equation}\label{equ22}\sigma_{ijk} \ = \ D_{ij}(\beta_{jk}) \ .
\end{equation}
Let
\[\Theta:H^1(X,T_X\ot L)\ot H^1(X,\Omega_X)\lra H^2(X,L)\]
be the canonical map. Then it follows from $(\ref{equ22})$ that \ 
\m{\delta(\beta)=\Theta(\eta\ot d_X(\beta))}. Since \ \m{d_X(\beta)=0} \ by 
Proposition \ref{prop8}, we have \ \m{\delta(\beta)=0}.

\end{sub}

\sepsec

\section{Extensions of families of vector bundles to higher 
multiplicity}\label{ext_F}

Let \m{Y=X_n} be a primitive multiple scheme of multiplicity $n$, with
underlying smooth projective irreducible variety $X$, and associated line 
bundle $L$ on $X$. 

If \m{X_n} can be extended to a primitive multiple scheme \m{X_{n+1}} 
of multiplicity \m{n+1}, we have a canonical exact sequence of sheaves on 
\m{X_n}
\begin{equation}\label{equ1}\Sigma_{X_{n+1}}:\qquad
0\lra L^n\lra\Omega_{X_{n+1}|X_n}\lra\Omega_{X_n}\lra 0 \ , \end{equation}
corresponding to \ 
\m{\sigma_{\Omega_{X_{n+1}}}\in\Ext^1_{\ko_{X_n}}(\Omega_{X_n},L^n)}.

Let $S$ be a smooth algebraic variety over $\C$, and 
\m{p_{X_n}:X_n\times S\to X_n}, \m{p_S:X_n\times S\to S} the projections. We 
can see \m{X_n\times S} as a primitive multiple scheme of multiplicity $n$, 
with induced smooth variety \m{X\times S} and associated line bundle 
\m{p_{X_n}^*(L)_{|X\times S}}. If \m{X_n} is extended to \m{X_{n+1}},
\m{X_{n+1}\times S} as a primitive multiple scheme of multiplicity \m{n+1} 
extending \m{X_n\times S}. We have \m{\Omega_{X_n\times 
S}=p_{X_n}^*(\Omega_{X_n})\oplus p_S^*(\Omega_S)} and
\[\sigma_{X_{n+1}\times S} \ \in \ \Ext^1_{\ko_{X_n\times 
S}}(p_{X_n}^*(\Omega_{X_n}),p_{X_n}^*(L^n))\oplus\Ext^1_{\ko_{X_n\times 
S}}(p_{S}^*(\Omega_S),p_{X_n}^*(L^n)) \ . \]
Then \m{\Sigma_{X_{n+1}\times S}} is \m{p_{X_n}^*(\Sigma_{X_n})\oplus 
p_S^*(\Sigma)}, where $\Sigma$ is the exact sequence
\xmat{0\ar[r] & 0\ar[r] & \Omega_S\fleq[r] & \Omega_S\ar[r] & 0 \ .}
It follows that the component of \m{\sigma_{X_{n+1}\times S}} in 
\m{\Ext^1_{\ko_{X_n\times S}}(p_{S}^*(\Omega_S),p_{X_n}^*(L^n))} vanishes. 

\sepsub

\Ssect{The family of extensions to \m{X_{n+1}} of a vector bundle on 
\m{X_n}}{fam_ext}

Let $\E$ be a vector bundle on \m{X_n}, and \ \m{E=\E_{|X}}. We keep the 
notations of \ref{ext_B}.

Let \m{V=\Ext^1_{\ko_{X_{n+1}}}(\E,E\ot L^n)}, and \m{\P=\P(V)}. Let 
\m{p:X_{n+1}\times\P\to X_{n+1}}, \m{q:X_{n+1}\times\P\to\P} be the 
projections. If
\[0\lra p^*(E\ot L^n)\ot q^*(\ko_\P(1))\lra\ke\lra p^*(\E)\lra 0\]
is an exact sequence of sheaves on \m{X_{n+1}\times\P}, the induced sequence on 
\m{X_{n+1}}
\xmat{0\ar[r] & E\ot L^n\ot V^*\fleq[d]\ar[r] & p_*(\ke)\ar[r] & 
\E\fleq[d]\ar[r] & 0\\
& p_*(p^*(E\ot L^n)\ot q^*(\ko_\P(1))) & & p_*(p^*(\E))}
is also exact. This defines a linear map
\[{\bf F}:\Ext^1_{\ko_{X_{n+1}\times\P}}(p^*(\E),p^*(E\ot L^n)\ot 
q^*(\ko_\P(1)))\lra\Ext^1_{\ko_{X_{n+1}}}(\E,E\ot L^n)\ot V^*=V\ot V^* \ . \]

\sepprop

\begin{subsub}\label{lem6}{\bf Proposition: } ${\bf F}$ is an isomorphism.
\end{subsub}
\begin{proof} We prove first the surjectivity of $\bf F$. Let \m{\phi\in V^*}, 
inducing \m{\ov{\phi}\in H^0(q^*(\ko_\P(1)))}, \m{\sigma\in V}, corresponding 
to the exact sequence
\[0\lra E\ot L^n\lra\ke\lra\E\lra 0 \ . \]
We have a commutative diagram
\xmat{0\ar[r] & E\ot L^n\ar[r]\ar[d]^{I\ot\phi} & \ke\ar[r]\ar[d] & 
\E\fleq[d]\ar[r] & 0\\
0\ar[r] & E\ot\L^n\ot V^*\ar[r] & \kf\ar[r] & \E\ar[r] & 0 \ , }
where the exact sequence in the bottom corresponds to \m{\sigma\ot\phi\in V\ot 
V^*} (cf. \cite{dr1b}, prop. 4.3.2). This diagram is the image (by \m{p_*}) of 
the following commutative diagram on \m{X_{n+1}\times\P}
\xmat{0\ar[r] & p^*(E\ot L^n)\ar[r]\ar[d]^{I\ot\ov{\phi}} & 
p^*(\ke)\ar[r]\ar[d] & p^*(\E)\fleq[d]\ar[r] & 0\\
0\ar[r] & p^*(E\ot L^n)\ot q^*(\ko_\P(1))\ar[r] & p^*(\kf)\ar[r] & 
p^*(\E)\ar[r] & 0 \ , }
where the exact sequence in the bottom is associated to \Nligne
\m{\sigma'\in\Ext^1_{\ko_{X_{n+1}\times\P}}(p^*(\E),p^*(E\ot L^n)\ot 
q^*(\ko_\P(1)))} \ such that \m{{\bf F}(\sigma')=\sigma\ot\phi}. It follows 
that \ \m{\sigma\ot\phi\in\imm({\bf F})}, and that $\bf F$ is surjective.

Now we prove that $\bf F$ is injective. For any coherent sheaves $\ke$, $\kf$ 
on a scheme $Y$, let \m{\Gamma_{\ke,\kf}} be the natural exact sequence
\[0\lra H^1(Y,\HHom(\ke,\kf))\lra\Ext^1_{\ko_Y}(\ke,\kf)\lra 
H^0(Y,\EExt^1_{\ko_Y}(\ke,\kf))\]
which is a part of \m{\bm{\Gamma}_{\ke,\kf}} (cf. \ref{Ex_Ext}). For 
\m{\Gamma_{\E,E\ot L^n\ot V^*}} on \m{X_{n+1}} we have
\[H^1(\HHom(\E,E\ot L^n\ot V^*)) \ = \ \Ext^1_{\ko_X}(E,E\ot L^n)\ot V^* \ , \]
\[H^0(\EExt^1_{\ko_{X_{n+1}}}(\E,E\ot L^n\ot V^*)) \ = \ \End(E)\ot V^* \ . \]
For \m{\Gamma_{p^*(\E),p^*(E\ot L^n)\ot q^*(\ko_\P(1))}} on \m{X_{n+1}\times\P} 
we have
\[\HHom(p^*(\E),p^*(E\ot L^n)\ot q^*(\ko_\P(1))) \ = \ p^*(\HHom(E,E\ot 
L^n))\ot q^*(\ko_\P(1)) \ , \]
whence
\[H^1(\HHom(p^*(\E),p^*(E\ot L^n)\ot q^*(\ko_\P(1)))) \ = \ 
\Ext^1_{\ko_X}(E,E\ot L^n)\ot V^* \ . \]
We have
\[\EExt^1_{\ko_{X_{n+1}\times\P}}(p^*(\E),p^*(E\ot L^n)\ot q^*(\ko_\P(1))) \ = 
\ \EEnd(p^*(E))\ot q^*(\ko_\P(1)) \ , \]
whence
\[H^0(\EExt^1_{\ko_{X_{n+1}\times\P}}(p^*(\E),p^*(E\ot L^n)\ot q^*(\ko_\P(1)))) 
\ = \ \End(E)\ot V^* \ . \]
Let \ \m{W=\Ext^1_{\ko_{X_{n+1}\times\P}}(p^*(\E),p^*(E\ot L^n)\ot 
q^*(\ko_\P(1)))}. We have then a commutative diagram (a morphism from 
\m{\Gamma_{p^*(\E),p^*(E\ot 
L^n)\ot q^*(\ko_\P(1))}} to \m{\Gamma_{\E,E\ot L^n\ot V^*}})
\xmat{0\ar[r] & \Ext^1_{\ko_X}(E,E\ot L^n)\ot V^*\ar[r]\fleq[d] &
W\ar[r]\ar[d]^{\bf F} & \End(E)\ot V^*\fleq[d]\\
0\ar[r] & \Ext^1_{\ko_X}(E,E\ot L^n)\ot V^*\ar[r] & 
\Ext^1_{\ko_{X_{n+1}}}(\E,E\ot L^n)\ot V^*\ar[r] & \End(E)\ot V^*
}
which implies easily the injectivity of $\bf F$.
\end{proof}

\sepprop

\begin{subsub}\label{ext_B3} The universal family of extensions of $\E$ -- \rm
Let
\[\bm{\sigma}\in\Ext^1_{\ko_{X_{n+1}\times\P}}(p^*(\E),p^*(E\ot L^n)\ot 
q^*(\ko_\P(1))) \ \simeq \ V\ot V^* \ , \]
corresponding to \m{I_V}, and
\[0\lra p^*(E\ot L^n)\ot q^*(\ko_\P(1))\lra\ke\lra p^*(\E)\lra 0\]
the corresponding extension on \ \m{X_{n+1}\times\P}. The sheaf $\ke$ is flat 
on $\P$, and for every \m{y\in\P}, the restriction of this exact sequence to \ 
\m{X_{n+1}\times\{y\}}
\[0\lra E\ot L^n\ot y^*\lra\ke_y\lra \E\lra 0\]
is exact (by Lemma \ref{lem3}), associated to the inclusion 
\m{i:y\hookrightarrow V}, \Nligne \m{i\in y^*\ot 
V=\Ext^1_{\ko_{X_{n+1}}}(\E,E\ot L^n\ot y^*)}.
\end{subsub}

Let \ \m{\delta:\Ext^1_{\ko_{X_{n+1}}}(\E,E\ot L^n)\to\End(E)} \ be the 
canonical map, and \ \m{V_\E=\delta^{-1}(I_E)}. Then \ 
\m{\ke_\E=\ke_{|X_{n+1}\times V_\E}} \ is a family of locally free sheaves on 
\m{X_{n+1}} parametrized by \m{V_\E}, flat on \m{V_\E}. Hence it is a vector 
bundle on \ \m{X_{n+1}\times V_\E}.

Let \m{p_X:X_{n+1}\times S\to X_{n+1}}, \m{p_S:X_{n+1}\times S\to S} be the 
projections. Let
\begin{eqnarray*}A & = & H^1(X_{n+1}\times S,\HHom(p_X^*(\E),p_X^*(E\ot 
L^n))), \\
B & = & H^0(X_{n+1}\times S,\EExt^1(p_X^*(\E),p_X^*(E\ot L^n))),\\
A' & = & H^0(\ko_S)\ot\Ext^1_{\ko_X}(E,E\ot L^n),\\
B' & = & H^0(\ko_S)\ot\End(E).
\end{eqnarray*}
We have canonical isomorphisms \m{\alpha:A\to A'}, \m{\beta:B\to B'}.

\sepprop

\begin{subsub}\label{lem8}{\bf Lemma: } If $S$ is affine, then we have a 
commutative diagram
\xmat{0\ar[r] & A\ar[r]\ar[d]^\alpha & \Ext^1_{\ko_{X_{n+1}\times 
S}}(p_X^*(\E), p_X^*(E\ot L^n))\ar[r]\ar[d]^{\gamma_S} & B\ar[d]^\beta\\
0\ar[r] & A'\ar[r] &  H^0(\ko_S)\ot\Ext^1_{\ko_{X_{n+1}}}(\E, E\ot L^n)\ar[r] & 
B',}
where the bottom exact sequence is \m{H^0(\ko_S)\ot\Gamma_{\E,E\ot L^n}}, and 
\m{\gamma_S} is an isomorphism.
\end{subsub}
\begin{proof}Analogous to the proof of Proposition \ref{lem6}.
\end{proof}

\sepprop

\begin{subsub}\label{prop3}{\bf Proposition: } Let $\ke$ be a vector bundle on 
\m{X_{n+1}\times S} such that for every closed 
point \m{s\in S}, \m{\ke_{s|X_n}\simeq\E}. Let \m{s_0\in S} be a closed point. 
Then there exists a neighbourhood \m{S_0\subset S} of \m{s_0} and a morphism \ 
\m{f:S_0\to V_\E} \ such that \ \m{f^\sharp(\ke_\E)\simeq\ke_{|{X_{n+1}\times 
S_0}}}.
\end{subsub}
\begin{proof} Let \ \m{\W=\HHom(p_X^*(\E),\ke_{|X_n\times S})}. The sheaf 
\m{p_{S*}(\W)} is a vector bundle on $S$. For every \m{s\in S} we have \ 
\m{p_{S*}(\W)_s=\Hom(\E,\ke_{s|X_n})}. There exists an affine  neighbourhood 
\m{S_0} of \m{s_0} and \ \m{\psi\in H^0(S_0,p_{S*}(\W))} \ such that for every 
\m{s\in S_0}, \m{\psi(s)} is an isomorphism. This section induces an 
isomorphism \
\m{\bm{\psi}:p_X^*(\E)_{|X_n\times S_0}\to\ke_{|X_n\times S_0}}.

Let \ \m{W=\HHom(p_X^*(E),\ke_{|X\times S})}; \m{p_{S*}(W)} is a vector bundle 
on $S$, and for every \m{s\in S}, \m{p_{S*}(W)_s=\Hom(E,\ke_{s|X})}. Let \ 
\m{\ov{\psi}\in H^0(S_0,p_{S*}(W))} \ induced by $\psi$. This section induces 
an isomorphism \ 
\m{\ov{\bm{\psi}}:p_X^*(E\ot L^n)_{|X\times S_0}\to\ke_{|X_\times 
S_0}\ot p_X^*(L^n)_{|X_\times S_0}}.

From the canonical exact sequence
\xmat{0\ar[r] & \ke_{|X\times S}\ot p_X^*(L^n)\ar[r]^-i & \ke\ar[r]^-p &
\ke_{|X_n\times S}\ar[r] & 0}
we get
\xmat{0\ar[r] & p_X^*(E\ot L^n)_{|X\times S_0}\ar[rr]^-{i\circ\ov{\bm{\psi}}} &
& \ke_{|X_{n+1}\times S_0}\ar[rr]^-{\ov{\bm{\psi}}^{-1}\circ p} & &
p_X^*(\E)_{|X_{n+1}\times S_0}\ar[r] & 0 \ ,}
associated to \ \m{\sigma\in\Ext^1_{\ko_{X_{n+1}\times 
S_0}}(p_X^*(\E)_{|X_{n+1}\times S_0}, p_X^*(E\ot L^n)_{|X\times S_0})}.
Then \ \m{f=\gamma_{S_0}(\sigma)} \ is a morphism \m{S_0\to V_\E} \ such that \ 
\m{f^\sharp(\ke_\E)\simeq\ke_{|{X_{n+1}\times S_0}}}.
\end{proof}

\end{sub}

\sepsub

\Ssect{Families and their extension to higher multiplicity}{FAM}

Let $\E$ be a vector bundle on \m{X_n\times S}, and \ \m{E=\E_{|X\times S}}. 
Let \m{p_X:X_{n+1}\times S\to X_{n+1}}, \m{p_S:X_{n+1}\times S\to S} be the 
projections.

\sepprop

\begin{subsub}\label{lem3}{\bf Lemma: } Let \m{s\in S} be a closed point. Then 
we have \ \m{\Tor^1_{\ko_{X_{n+1}\times S}}(\E,\ko_{X_{n+1}\times\{s\}})=0}.
\end{subsub}
\begin{proof} Let \m{x\in X} and \m{W\subset X_{n+1}\times S} a neighbourhood 
of \m{(x,s)} such that $\E$ is free on $W$. It suffices to prove that
\[\Tor^1_{\ko_W}(\ko_{W\cap(X_n\times S)},\ko_{W\cap(X_{n+1}\times\{s\})})
\ = \ 0 . \]
For this it suffices to show that
\[\Tor^1_{\ko_{X_{n+1}\times S}}(\ko_{X_n\times S},\ko_{X_{n+1}\times\{s\}})
\ = \ 0 \ . \]
Let \m{U_{n+1}} a neighbourhood of $x$ in \m{X_{n+1}}, \m{U_n} the 
corresponding neighbourhood in \m{X_n}, such that \m{\ki_{X,X_{n+1}|U_{n+1}}} 
is globally generated by \m{t\in\ko_{X_{n+1},x}(U_{n+1})}. The result follows 
from the locally free resolution of \m{\ko_{U_n\times S}}:
\xmat{\cdots\ko_{U_{n+1}\times S}\ar[rr]^-{\times t^n} & & 
\ko_{U_{n+1}\times S}\ar[rr]^-{\times t} & & \ko_{U_n\times S}\ar[r] & 0  \ . }
\end{proof}

\sepprop

It follows that an extension \ \m{0\to E\ot p_X^*(L^n)\to\ke\to\E\to 0} \
restricts on \ \m{X_{n+1}\times\{s\}} \ to an extension \ \m{0\to E_s\ot 
L^n\to\ke_s\to\E_s\to 0}.

Let $\E$ be a vector bundle on \m{X_n\times S}, and \ \m{E=\E_{|X\times S}}. We 
have
\[\nabla_0(\E) \ \in \ \Ext^1_{\ko_{X_n\times 
S}}(\E,\E\ot p_{X_n}^*(\Omega_{X_n}))\oplus\Ext^1_{\ko_{X_n\times 
S}}(\E,\E\ot p_S^*(\Omega_S)) \ . \]
We have \ \m{\Delta(\E)=\ov{\sigma}_{X_{n+1}\times S}\nabla_0(\E)}. The 
component of \m{\ov{\sigma}_{X_{n+1}\times S}} in \Nligne
\m{\Ext^1_{\Omega_{X_n}\times S}(\E\ot p_{S}^*(\Omega_S),\E\ot p_{X_n}^*(L^n))} 
vanishes. Let \m{\nabla_1} be the component of \m{\nabla_0(\E)} in \Nligne
\m{\Ext^1_{\ko_{X_n\times S}}(\E,\E\ot p_{X_n}^*(\Omega_{X_n}))}. We have then
\begin{equation}\label{equ4}
\Delta(\E) \ = \ \ov{\sigma}_{X_{n+1}\times S}\nabla_1 \ .
\end{equation}

\sepprop

\begin{subsub}\label{hyp1} Obstruction to the extension of the family -- \rm 
Let \m{p_{S,n}}, \m{p_{S,1}} denote the projections \ \m{X_n\times S\to S}, 
\m{X\times S\to S} \ respectively.
We suppose that for \m{s\in S},
\begin{enumerate}
\item[(i)] $\dim(\Hom(\E_s,\E_s\ot\Omega_{X_n}))$ \ is independent of $s$.
\item[(ii)] $\dim(\End(E_s))$ is independent of $s$,
\item[(iii)] $\dim(\Ext^1_{\ko_X}(E_s,E_s\ot L^n))$ and 
$\dim(\Ext^1_{\ko_{X_n}}(\E_s,\E_s\ot \Omega_{X_n}))$ are independent of $s$,
\item[(iv)] $\dim(\Ext^2_{\ko_X}(E_s,E_s\ot L^n))$ is independent of $s$.
\end{enumerate}
It follows that \ \m{\B=R^1p_{S,n*}(\E^*\ot\E\ot p_{X_n}^*(\Omega_{X_n}))} \ and
\ \m{\D=R^2p_{S,1*}(E^*\ot E\ot L^n)} \ are vector bundles on $S$, with
\[\B_s=\Ext^1_{\ko_{X_n}}(\E_s,\E_s\ot\Omega_{X_n}) \ , \
\D_s=\Ext^2_{\ko_X}(E_s,E_s\ot L^n) \ . \]
\end{subsub}

\sepprop

Let \m{s\in S} be a closed point and
\[R_s:\Ext^2_{\ko_{X\times S}}(E,E\ot p_X^*(L^n))\lra
\Ext^2_{\ko_X}(E_s,E_s\ot L^n)\]
be the canonical map.

\sepprop

\begin{subsub}\label{prop2}{\bf Proposition: } We have \ \m{R_s(\Delta(\E))=
\Delta(\E_s)}.
\end{subsub}
\begin{proof}
Let \ \m{{\bf B}=\Ext^1_{\ko_{X_n\times S}}(\E,\E\ot p_{X_n}^*(\Omega_{X_n}))}, 
\m{{\bf D}=\Ext^2_{\ko_{X\times S}}(E,E\ot p_X^*(L^n))}.
We have a canonical map \ \m{\phi:{\bf B}\to{\bf D}} (multiplication by 
\m{\sigma_{X_{n+1}\times S}}), and \m{\psi:\B\to\D}, where for every 
\m{s\in S}, \m{\psi_s} is the canonical map
\[\Ext^1_{\ko_{X_n}}(\E_s,\E_s\ot\Omega_{X_n})\lra\Ext^2_{\ko_{X_n\times 
S}}(\E_s,\E_s\ot L^n)\]
(multiplication by \m{\sigma_{X_{n+1}}}), and a canonical diagram
\xmat{{\bf B}\ar[rr]^-\Phi\ar[d] & & {\bf D}\ar[d]\\
H^0(S,\B)\ar[rr] & & H^0(S,\D)}
where the vertical arrows come from Leray's spectral sequence. Evaluation at 
$s$ gives the diagram
\xmat{{\bf B}\ar[rr]^-\Phi\ar[d] & & {\bf D}\ar[d]\\
\Ext^1_{\ko_{X_n}}(\E_s,\E_s\ot\Omega_{X_n})\ar[rr] & &
\Ext^2_{\ko_X}(E_s,E_s\ot L^n) \ .}
The result follows from $(\ref{equ4})$, $(\ref{equ2})$ and $(\ref{equ3})$.
\end{proof}

\sepprop

\begin{subsub}\label{coro1}{\bf Corollary: } Suppose that $S$ is affine, and 
that for every closed point \m{s\in S} we have \ \m{\Delta(\E_s)=0}. Then we 
have \ \m{\Delta(\E)=0}.
\end{subsub}
\begin{proof} We have, since $S$ is affine, by Leray's spectral sequence
\[\Ext^2_{\ko_{X\times S}}(E,E\ot p_X^*(L^n)) \ = \ H^0(S,\D) \ , \]
and the result follows from Proposition \ref{prop2}.
\end{proof}

\sepprop

In other words, if $S$ is affine, and if for every \m{s\in S}, \m{\E_s} can be 
extended to a vector bundle on \m{X_{n+1}}, then $\E$ can be extended to a 
vector bundle on \m{X_{n+1}\times S}.

\sepprop

\begin{subsub}\label{ext_B4} The universal family of extensions -- \rm We 
suppose that the hypotheses of \ref{hyp1} are satisfied, that $S$ is affine, 
and that for every \m{s\in S}, $\Delta(\E_s)=0$. From the exact sequence of 
\ref{ext_4} it follows that 
\m{\dim(\Ext^1_{\ko_{X_{n+1}}}(\E_s,E_s\ot L^n))} is independent of $s$.
\end{subsub}

The relative Ext-sheaf (cf. \cite{la}) \m{\V=\EExt_{p_S}^1(\E,E\ot p_X^*(L^n))} 
\ is locally free, with, for every closed point \m{s\in S}, 
\m{\V_s=\Ext^1_{\ko_{X_{n+1}}}(\E_s,E_s\ot L^n))}. Let \m{\kp=\P(\V)}. Let 
\m{p:X_{n+1}\times\kp\to X_{n+1}\times S}, \m{q:X_{n+1}\times\kp\to\kp} be the 
projections. If
\[0\lra p^*(E\ot p_X^*(L^n))\ot q^*(\ko_\kp(1))\lra\ke\lra p^*(\E)\lra 0\]
is an exact sequence of sheaves on \m{X_{n+1}\times\kp}, the induced sequence 
on \m{X_{n+1}\times S}
\xmat{0\ar[r] & E\ot L^n\ot \V^*\fleq[d]\ar[r] & p_*(\ke)\ar[r] & 
\E\fleq[d]\ar[r] & 0\\
& p_*(p^*(E\ot p_X^*(L^n))\ot q^*(\ko_\kp(1))) & & p_*(p^*(\E))}
is also exact. This defines a linear map
\[\bm{\Phi}:\Ext^1_{\ko_{X_{n+1}\times\kp}}(p^*(\E),p^*(E\ot p_X^*(L^n))\ot 
q^*(\ko_\kp(1)))\lra H^0(S,\V\ot\V^*) \ . \]
The proof of the following result is similar to that of Proposition \ref{lem6}:

\sepprop

\begin{subsub}\label{lem7}{\bf Proposition: } $\bm{\Phi}$ is an isomorphism.
\end{subsub}

\sepprop

Let
\[\bm{\sigma}\in\Ext^1_{\ko_{X_{n+1}\times\kp}}(p^*(\E),p^*(E\ot p_X^*(L^n))\ot 
q^*(\ko_\kp(1))) \ \simeq \ H^0(\V\ot\V^*) \ , \]
corresponding to \m{I_\V}, and
\[0\lra p^*(E\ot p_X^*(L^n))\ot q^*(\ko_\kp(1))\lra\bm{\ke}\lra p^*(\E)\lra 0\]
the corresponding extension on \ \m{X_{n+1}\times\kp}. The sheaf $\bm{\ke}$ is 
flat on $\kp$. For every \m{s\in S}, and for every \m{y\in\kp_s}, the 
restriction of this exact sequence to \ \m{X_{n+1}\times\{y\}}
\[0\lra E_s\ot L^n\ot y^*\lra\bm{\ke}_y\lra \E_s\lra 0\]
is exact (by Lemma \ref{lem3}), associated to the inclusion 
\m{i:y\hookrightarrow\V_s},\Nligne \m{i\in y^*\ot 
\V_s=\Ext^1_{\ko_{X_{n+1}}}(\E_s,E_s\ot L^n\ot y^*)}.

The sheaf $\V$ is associated to the presheaf $\kv$ defined by
\[\kv(U) \ = \ \Ext^1_{\ko_{X_{n+1}\times U}}(p^*(\E)_{|X_{n+1}\times U},
p^*(E\ot p_X^*(L^n))_{|X_{n+1}\times U})\]
for every open subset \m{U\subset S}. We have a canonical morphism
\[\Ext^1_{\ko_{X_{n+1}\times U}}(\E_{|X_{n+1}\times U},
(E\ot p_X^*(L^n))_{|X_{n+1}\times U})\lra\qquad\qquad\qquad\qquad \]
\[\qquad\qquad\qquad\to H^0(\EExt^1_{\ko_{X_{n+1}\times 
U}}(\E_{|X_{n+1}\times U},(E\ot p_X^*(L^n))_{|X_{n+1}\times 
U}))=\End(E_{|X_{n+1}\times U})\]
(coming from the Ext spectral sequence) inducing a morphism \ 
\m{\bm{\delta}:\V\to p_*(\EEnd(E))}. For every \m{s\in S}, \m{\bm{\delta}_s} is 
the canonical map
\[\Ext^1_{\ko_{X_{n+1}}}(\E_s,E_s\ot L^n)\lra 
H^0(\EExt^1_{\ko_{X_{n+1}}}(\E_s,E_s\ot L^n))=\End(E_s)\]
We have a canonical section $\mu$ of \m{p_*(\EEnd(E))}, such that 
\m{\mu(s)=I_{E_s}} for every \m{s\in S}. Let
\[\V_\E \ = \ \bm{\delta}^{-1}(\imm(\mu)) \ , \]
and \m{\kp_\E} the image of \m{\V_\E} in $\kp$, which is a smooth subvariety of 
$\kp$. The projection \m{\V_\E\to\kp_\E} is an isomorphism.
The sheaf \ \m{\bm{\ke}_\E=\ke_{|X_{n+1}\times\kp_\E}} \ is a vector 
bundle on \ \m{X_{n+1}\times\kp_\E} (and \m{X_{n+1}\times\V_\E}). We have a 
canonical exact sequence on \m{X_{n+1}\times\V_\E} :
\[0\lra\pi^\sharp(E\ot p_X^*(L^n))\lra\ke_\E\lra\pi^\sharp(\E)\lra 0\]
(where $\pi$ is the projection \m{\V_\E\to S}).

\sepprop

\begin{subsub}\label{537} The case of non affine $S$ -- \rm We don't suppose as 
in Corollary \ref{coro1} that $S$ is affine. It is easy to see that the 
schemes \m{\V_{\E_{|U}}}, \m{\kp_{\E_{|U}}} on affine open subsets \m{U\subset 
S}, can be glued in an obvious way to define the schemes \m{\V_\E}, \m{\kp_\E} 
over $S$. Similarly the vector bundles \m{\ke_{\E_{|\kp_{\E_{|U}}}}} define the 
vector bundle \m{\ke_\E} on \m{X_{n+1}\times\kp_\E}.
\end{subsub}

\sepprop

\begin{subsub}\label{prop4}{\bf Proposition: } Let \m{f:T\to S} be a morphism 
of schemes and $\kf$ a vector bundle on \m{X_{n+1}\times T} such that \ 
\m{\kf_{|X_n\times T}\simeq f^\sharp(\E)}. Then there exists a morphism
\ \m{\delta:T\to\kp_\E} \ over $S$ such that for every closed point \m{t\in T} 
there exists an open neighbourhood \m{T_0} ot $t$ such that \
\m{\delta^\sharp(\ke_\E)_{|X_{n+1}\times T_0}\simeq\kf_{|X_{n+1}\times T_0}}.
\end{subsub}
\begin{proof} Let \m{\V_f=\EExt^1_f(f^\sharp(\E),f^\sharp(E\ot p_X^*(L^n))} be 
the relative Ext-sheaf. By \cite{a-k}, Theorem 1.9, for every affine open subset
\m{U\subset S} containing \m{f(t)}, the canonical morphism
\[\EExt^1_{\ko_{X_{n+1}\times U}}(\E_{|X_{n+1}\times U},E_{|X_{n+1}\times U}\ot 
p_X^*(L^n))\ot_{\ko_U}\ko_{f^{-1}(U)}\lra\qquad\qquad\qquad\qquad\]
\[\qquad\qquad\qquad\qquad\EExt^1_{\ko_{X_{n+1}\times 
f^{-1}(U)}}(f^\sharp(\E_{|X_{n+1}\times U}),f^\sharp(E_{|X_{n+1}\times U}))
\]
is an isomorphism. It follows that we have a canonical isomorphism
\begin{equation}\label{equ6}\V_f \ \simeq \ f^\sharp(\V) \ .
\end{equation}
We define in an obvious 
way the schemes \m{\V_{f^\sharp(\E)}\simeq\kp_{f^\sharp(\E)}} analogous to 
\m{\V_\E\simeq\kp_\E} respectively, which are closed subschemes of \m{\V_f}.
We have a canonical morphism of schemes \ \m{\phi:\V_f\to\V}, coming from  
$(\ref{equ6})$, and we have \ 
\m{\phi(\V_{f^\sharp(\E)})\subset V_\E}. In the same way, we have a 
canonical exact sequence
\begin{equation}\label{equ5}
0\lra{\pi'}^\sharp(E\ot 
p_X^*(L^n))\lra\ke_{f^\sharp(\E)}\lra{\pi'}^\sharp(\E)\lra 0
\end{equation}
on \m{X_{n+1}\times\V_{f^\sharp(\E)}} (where \m{\pi'} is the projection 
\m{\V_{f^\sharp(\E)}\to T}), and \ 
\m{\ke_{f^\sharp(\E)}\simeq\phi^\sharp(\ke_\E)}. The canonical exact sequence
\xmat{0\ar[r] & f^\sharp(E\ot p_X^*(L^n))\ar[r] & \kf\ar[r] & 
f^\sharp(\E)\fleq[d]\ar[r] & 0\\ & & & \kf_{|X_n\times T}} defines a section \ 
\m{\sigma:T\to\V_{f^\sharp(\E)}} \ of \ \m{\V_{f^\sharp(\E)}\to T}. Let
\[\delta=\phi\circ\sigma:T\lra\kp_\E=\V_\E \ . \]
The inverse image of $(\ref{equ5})$ by $\phi$ is an exact sequence
\xmat{0\ar[r] & f^\sharp(E\ot p_X^*(L^n))\ar[r] & 
\sigma^\sharp(\ke_{f^\sharp(\E)})\ar[r]\fleq[d] &  f^\sharp(\E)\ar[r] & 0 \ .\\
& & \delta^\sharp(\ke_\E)}
Over \m{f^{-1}(U)} this exact sequence corresponds to the same element of 
\Nligne
\m{\Ext^1_{\ko_{X_{n+1}\times f^{-1}(U)}}(f^\sharp(\E)_{|X_{n+1}\times 
f^{-1}(U)},
f^\sharp(E\ot p_X^*(L^n))_{|X_{n+1}\times f^{-1}(U)})}, and this implies that 
there exists a commutative diagram
\xmat{0\ar[r] & f^\sharp(E\ot p_X^*(L^n))_{|X_{n+1}\times 
f^{-1}(U)}\fleq[d]\ar[r] & \delta^\sharp(\ke_\E))_{|X_{n+1}\times 
f^{-1}(U)}\ar[d]^{\psi_U}\ar[r] & f^\sharp(\E)_{|X_{n+1}\times 
f^{-1}(U)}\fleq[d]\ar[r] & 0\\
0\ar[r] & f^\sharp(E\ot p_X^*(L^n))_{|X_{n+1}\times f^{-1}(U)}\ar[r] & 
\kf_{|X_{n+1}\times f^{-1}(U)}\ar[r] & f^\sharp(\E)_{|X_{n+1}\times 
f^{-1}(U)}\ar[r] & 0}
where \m{\psi_U} is an isomorphism. We take \ \m{T_0=f^{-1}(U)}.
\end{proof}

\sepprop

\begin{subsub}\label{affin} The affine bundle structure of \m{\kp_\E} -- \rm
Let \ \m{\A=R^1p_{S*}\big(\HHom(\E,E\ot p_X^*(L^n))\big)}. It is a vector 
bundle on $S$ (by \ref{hyp1}, (ii)), and for every closed point \m{s\in S} we 
have \Nligne \m{\A_s=\Ext^1_{\ko_X}(E_s,E_s\ot L^n)}. For every open subset 
\m{U\subset S}, we have a canonical injective map
\xmat{\A(U)\flinc[r] & \Ext^1_{\ko_{X_{n+1}\times 
U}}(p^*(\E)_{|X_{n+1}\times U},p^*(E\ot p_X^*(L^n))_{|X_{n+1}\times U})=\kv(U) 
\ , }
which defines an injective morphism of vector bundles \ \m{\A\to\V}. It follows 
that \m{\kp_\E} has a natural structure of affine bundle over $S$, with 
associated vector bundle $\A$. By Corollary \ref{coro1}, for every affine open 
subset \m{U\subset S}, \m{\kp_{\E|U}} has a section, hence is actually a vector 
bundle isomorphic to \m{\A_{|U}}.
\end{subsub}
Suppose that \ \m{\Hom(E_s,E_s\ot L^n)=\nsp} \ for every \m{s\in S}.
Since \ \m{\Delta(\E_s)=0} \ for every \m{s\in S}, we have \ 
\m{\Delta(\E)\in\imm(\lambda)}, where $\lambda$ is the canonical map
\[H^1(S,R^1p_{S*}(\HHom(\E,E\ot p_X^*(L^n))))\lra H^2(X_n\times S,\HHom(\E,E\ot 
p_X^*(L^n)))\]
(cf. \ref{Ex_Ler}).

Recall that \m{\eta(\kp_\E)} is the element of \m{H^1(S,\A)} associated to 
\m{\kp_\E} (cf. \ref{aff}).
The following is a consequence of Proposition \ref{prop5}:

\sepprop

\begin{subsub}\label{theo3}{\bf Theorem: } We have \ 
\m{\eta(\kp_\E)=\C.\Delta(\E)}.
\end{subsub}

If we don't assume that \ \m{\Hom(E_s,E_s\ot L^n)=\nsp} \ for every \m{s\in S}, 
then we get still an analogous more complicated result, using 3.2.7.

\end{sub}

\sepsub

\Ssect{Families of simple vector bundles}{FAM_s}

We keep the notations of \ref{ext_B} and \ref{FAM}. We don't suppose that $S$ 
is affine. To the hypotheses of \ref{hyp1} we add the following: for every 
\m{s\in S}, \m{\E_s} is simple, and for every \m{m\in\kp_\E}, \m{\ke_{\E,m}} is 
simple. 

We have obvious subbundles of the bundle $\A$ of \ref{affin}:
\[\ko_S\ot\imm(\delta^0) \ \subset \ \ko_S\ot H^1(X,L^n) \ \subset \ \A 
\]
(cf. \ref{isom_cl}).
Let \ \m{\kq_\E=\kp_\E/(\ko_S\ot\imm(\delta^0))} (cf. \ref{aff}), and \ 
\m{\omega:\kp_\E\to\kq_\E} the projection, which is an affine bundle with 
associated vector bundle \ \m{\ko_{\kq_\E}\ot\imm(\delta^0)}.

\sepsubsub

\begin{subsub}\label{glu_q} Invariance of \m{\kq_\E} -- \rm Let \m{\E'} be a 
vector bundle on \ \m{X_n\times S}, \m{E'=\E'_{|X\times S}}, and \m{f:\E\to\E'} 
an isomorphism, inducing obvious isomorphisms \ \m{p_f:\kp_\E\to\kp_{\E'}}, 
\m{q_f:\kq_\E\to\kq_{\E'}}. Let \m{\A'} be the vector bundle on $S$ analogous 
to $\A$, corresponding to \m{\E'}, i.e. for every \m{s\in S}, 
\m{\A'_s=\Ext^1_{\ko_X}(E'_s,E'_s\ot L^n)}. Then there are canonical obvious 
isomorphisms \ \m{\A'\simeq\A}, \m{\kq_\E'\simeq\kq_\E}, an they do not depend 
on $f$.
\end{subsub}

\sepsubsub

\begin{subsub}\rm
Let \m{s\in S}. Then there exists an open neighbourhood \m{S_0} of $s$ such 
that:
\begin{enumerate}
\item[--] $\A_{|S_0}$ is trivial and we have a subbundle 
\m{\D\subset\A_{|S_0}} such that \ \m{\A_{|S_0}=(\ko_{S_0}\ot\imm(\delta^0))
\oplus\D}, 
\item[--] there exists a section \ \m{\nu:S_0\to\kp_{\E|S_0}} \ of \ 
\m{\pi:\kp_{\E|S_0}\to S_0}.
\end{enumerate}
It follows that the affine bundle \ \m{\omega:\kp_{\E|S_0}\to\kq_{\E|S_0}} \ 
has a section $\tau$, so that we can view \m{\kq_{\E|S_0}} as a closed 
subvariety of \m{\kp_{\E|S_0}}. Let \ \m{\eta:\kq_{\E|S_0}\to S_0} \ be the 
projection.
\end{subsub}

\sepprop

\begin{subsub}\label{prop16}{\bf Proposition:} Let \m{m_0\in\pi^{-1}(s)}.
There exists an open affine neighbourhood \Nligne 
\m{M(m_0)\subset\pi^{-1}(S_0)} 
\ of $m$ such that
\m{\ke_{\E|X_{n+1}\times M(m_0)}\simeq 
\omega^*(\ke_{\E|X_{n+1}\times\kq_{\E|S_0}})_{|X_{n+1}\times M(m_0)}} .
\end{subsub}
\begin{proof} Let
\[\U \ = \ \ke_{\E|X_{n+1}\times\pi^{-1}(S_0)} \ , \quad
\V \ = \ \omega^\sharp(\ke_{\E|X_{n+1}\times\kq_{\E|S_0}}) \ . \]
For every \m{m\in \pi^{-1}(S_0)} we have \m{\U_m\simeq\V_m} (by Proposition 
\ref{prop15}), and \Nligne \m{\Hom(\U_m,\V_m)\simeq H^0(\ko_{X_{n+1}})}.

Let \ \m{p_0:X_{n+1}\times\pi^{-1}(S_0)\to\pi^{-1}(S_0)} \ be the projection.
Then \ \m{p_{0*}(\HHom(\U,\V))} \ is a vector bundle. It follows that there 
exists an open affine neighbourhood \ \m{M(m_0)\subset\pi^{-1}(S_0)} \ of $m$ 
such that there exists a section \ \m{\sigma\in 
H^0(M(m_0),p_{0*}(\HHom(\U,\V)))} \ 
such that for every \m{m\in M(m_0)}, \m{\sigma(m)\in\Hom(\U_m,\V_m)} is an 
isomorphism. We obtain in this way an isomorphism \ \m{\U_{|X_{n+1}\times 
M(m_0)}\simeq\V_{|X_{n+1}\times M(m_0)}}.
\end{proof}

\sepprop

By Proposition \ref{prop16} we can find finite open covers \m{(S_i)_{1\leq 
i\leq k}} of $S$, \m{(P_i)_{1\leq i\leq k}} of \m{\kp_\E}, such that 
\m{P_i\subset\pi^{-1}(S_i)}, and vector bundles \m{\F_i} on 
\m{X_{n+1}\times\eta^{-1}(S_i)} such that \ \m{\ke_{\E|X_{n+1}\times 
P_i}\simeq\omega^\sharp(\F_i)_{|X_{n+1}\times P_i}}.

\sepprop

\begin{subsub}\label{theo10}{ \bf Theorem: } Let \m{f:T\to S} be a morphism 
of schemes and $\kf$ a vector bundle on \m{X_{n+1}\times T} such that \ 
\m{\kf_{|X_n\times T}\simeq f^\sharp(\E)}. Then there exists a unique morphism
\ \m{\gamma_\kf:T\to\kq_\E} \ such that for every closed point \m{t\in T}, 
there exists $i$, \m{1\leq i\leq k}, such that \m{f(t)\in S_i}, a neighbourhood 
\m{U\subset T} of $t$ such that \ \m{\gamma_\kf(U)\subset\eta^{-1}(S_i)} \ and \
\m{\kf_{|X_{n+1}\times U}\simeq\gamma_\kf^\sharp(\F_i)_{|X_{n+1}\times U}}.
\end{subsub}
\begin{proof}
This follows from Propositions \ref{prop4} and \ref{prop16}. The unicity of 
\m{\gamma_\kf} comes from the fact that if \m{1\leq i,j\leq k}, and 
\m{q_i\in\eta^{-1}(S_i)}, \m{q_j\in\eta^{-1}(S_j)} are distinct, then by 
Proposition \ref{prop15}, the vector bundles \m{\F_{i,q_i}} and \m{\F_{j,q_j}} 
on \m{X_{n+1}} are not isomorphic.
\end{proof}

\sepprop

Since points in \m{\kq_{\E_s}} correspond exactly to extensions of \m{\E_s} to 
\m{X_{n+1}}, we have

\sepprop

\begin{subsub}\label{prop17}{ \bf Proposition: }
Let \m{f:T\to S} be a morphism of schemes and \m{\kf_1}, \m{\kf_2} vector 
bundles on \m{X_{n+1}\times T}, and \ \m{\theta_i:\kf_{i|X_n\times T}\simeq 
f^\sharp(\E)} \ for \m{i=1,2}. Then the induced isomorphism \ 
\m{\theta=\theta_2^{-1}\theta_1:\kf_{1|X_n\times T}\to\kf_{2|X_n\times T}} \
induces and isomorphism of affine bundles over $T$ \ 
\m{\kq_{\kf_2}\simeq\kq_{\kf_1}} \ which is independent of \m{\theta_1}, 
\m{\theta_2}.
\end{subsub}

\end{sub}

\sepsub

\Ssect{Moduli spaces of vector bundles}{mod}

We use the notations of \ref{FAM} and \ref{FAM_s}. 

Let \m{\chi_n} be a set of isomorphism classes of vector bundles on \m{X_n}. 
Suppose that there is a fine moduli space for \m{\chi_n}, defined by a smooth 
irreducible variety \m{M_n}, an open cover \m{(\km_i)_{i\in I}} of \m{M_n}, and 
for every \m{i\in I},a vector bundle \m{\E_i} on \m{X_n\times\km_i} (cf. 
\ref{Fine2}). 
We suppose that for every \m{\ke\in\chi_n}, if \m{E=\ke_{|X}},
\begin{enumerate}
\item[(i)] $\ke$ and $E$ are simple,
\item[(ii)] $\dim(\Ext^1_{\ko_X}(E,E\ot L^n))$ and 
$\dim(\Ext^1_{\ko_{X_n}}(\ke,\ke\ot\Omega_{X_n}))$ are independent of $\ke$,
\item[(iii)] $\dim(\Ext^2_{\ko_X}(E,E\ot L^n))$ is independent of $\ke$.
\item[(iv)] $\Delta(\ke)=0$ .
\item[(v)] every vector bundle $\kf$ on \m{X_{n+1}} such that \ 
\m{\kf_{|X_n}\simeq\ke} \ is simple. 
\end{enumerate}
In particular the conditions of \ref{hyp1} are verified (for the vector bundles 
\m{\E_i}). Hence the smooth variety \m{\kp_{\E_i}} and the vector bundle 
\m{\ke_{\E_i}} are defined. Let \m{\chi_{n+1}} the set of isomorphism 
classes of vector bundles $\kf$ on \m{X_{n+1}} such that 
\m{\kf_{|X_n}\in\chi_n}. 

The condition (v) is satisfied if the hypotheses of Lemma \ref{endo0} are 
satisfied by every \m{\ke\in\chi_n}. This is the case if $L$ is a non trivial 
ideal sheaf. Condition (v) is also satisfied if \m{\chi_n} consists of line 
bundles.

For every \m{i\in I}, we will apply Theorem \ref{theo10} to \m{\E_i}:
we obtain an affine bundle \m{\eta_i:\kq_{\E_i}\to\km_i}, an open cover 
\m{(\km_{i,k})_{k\in J_i}} of \m{\km_i}, and for every \m{k\in J_i} a vector 
bundle \m{\F_{i,k}} on \ \m{X_{n+1}\times\eta_i^{-1}(\km_{i,k})}.

Finally we get a finer open cover \m{(\kn_j)_{j\in J}} of \m{M_n}, and for 
every \m{j\in J}, a vector bundle \m{\G_j} on \m{X_n\times\kn_j} satisfying the 
conditions of \ref{Fine2}). In particular, for every \m{j,k\in J}, we can find 
an open cover of \m{\kn_j\cap\kn_k} such that for every $U$ in this cover, 
we have \ \m{\G_{j|\kn_i\cap\kn_k\cap U}\simeq \G_{k|\kn_i\cap\kn_k\cap U}}. 
Moreover, for every \m{j\in J}, we have a vector bundle \m{\A_j} on \m{\kn_j}, 
an affine bundle \ \m{\eta_j:\kq_{\G_j}\to\kn_j} with associated vector bundle 
\m{\A_j}, and a vector bundle \m{\H_j} on \ \m{X_{n+1}\times\kq_{\G_i}}.

From Proposition \ref{prop17}, all the vector bundles \m{\A_j}, the affine 
bundles \m{\kq_{\G_j}} can be glued to define respectively a vector bundle $\bf 
A$ on \m{M_n} and an affine bundle \m{\kq} on \m{M_n} with 
associated vector bundle $\bf A$.

\sepprop

\begin{subsub}\label{lem12}{\bf Lemma:} The scheme $\kq$ is separated.
\end{subsub}
\begin{proof} According to \cite{kriz}, Proposition 2.6.1, we need to prove the 
following statement: let $Z$ be a smooth 
variety, $U$, $V$ affine open subsets of $Z$ such that \m{U\cap V} is affine, 
$E$ (resp. $F$) a vector bundle on $U$ (resp. $V$), and \ \m{\phi:E_{|U\cap 
V}\to F_{|U\cap V}} \ an isomorphism. Let $Y$ be the scheme obtained by gluing 
$E$ and $F$ using $\phi$. Then $Y$ is separated. But this is obvious.
\end{proof}

\sepprop

Let \m{\chi_{n+1}} denote the set of isomorphism classes of vector bundles on 
\m{X_{n+1}} whose restriction to \m{X_n} belongs to \m{\chi_n}. From Theorem 
\ref{theo10} follows easily

\sepprop

\begin{subsub}\label{theo11}{\bf Theorem:} The smooth variety $\kq$, the open 
cover \m{(\kq_{\G_j})} of $\kq$ and the vector bundles \m{\H_j} on \ 
\m{X_{n+1}\times\kq_{\G_j}} \ define a fine moduli space for \m{\chi_{n+1}}.
\end{subsub}

\sepprop

This proves Theorem \ref{theo_0_1} if \m{N_n=M_n}. The proof in the general 
case is similar, using the appropriate definition of \m{N_n} given in 
\ref{extens}.

\end{sub}

\sepsub

\Ssect{The subvariety of extensible vector bundles}{extens}

We use the notations of \ref{mod}, except hypothesis (iv). We will indicate 
briefly how to define the subvariety \m{N_n\subset M_n} of bundles that can be 
extended to \m{X_{n+1}}. If this variety is smooth then one can easily modify 
the preceding proofs to obtain a fine moduli space for the vector bundles on 
\m{X_{n+1}} whose restriction to \m{X_n} belongs to $\chi$, and this moduli 
space is an affine bundle over \m{N_n}.

Let \ \m{p_X:X\times\km_i\to X}, \m{p_M:X\times\km_i\to\km_i} \ be the 
projections, and \ \m{E_i=\E_{i|X\times\km_i}}. Then
\[\W_i \ = \ R^1p_{M*}(E_i^*\ot E_i\ot p_X^*(L^n))\]
is a vector bundle on \m{\km_i}, and since all the vector bundles in $\chi$ are 
simple, these vector bundles can be glued naturally to define a vector bundle 
$\W$ on \m{M_n}. We have a canonical section of each \m{\W_i}
\[\xymatrix@R=5pt{\sigma_i:\km_i\ar[r] & \W_i\\ m\fmaps[r] & \Delta(\E_{i,m}) 
.}\]
These sections can be glued to define \ \m{\sigma\in H^0(M_n,\W)}. The variety 
\m{N_n} is the zero subscheme scheme of $\sigma$.

\end{sub}

\sepsec

\section{Picard groups}\label{pic}

We use the notations of \ref{PMS}. Let \m{Y=X_n} be a primitive multiple scheme 
of multiplicity \m{n\geq 2}, with underlying smooth projective irreducible 
variety $X$, and associated line bundle $L$ on $X$. 

\sepsub

\Ssect{Inductive definition}{Pic_univ}

If \ \m{{\bf P}=\Pic^d(X)} \ is an irreducible component of \m{\Pic(X)}, let 
\m{\kl_{n,\bf P}} be the set of line bundles $\kb$ on \m{X_n} such that \ 
\m{\kb_{|X}\in{\bf P}}. Let \m{{\bf P}_0} be the component that contains 
\m{\ko_X}. Then \ \m{\kl_{n,{\bf P}_0}\not=\emptyset} (it contains 
\m{\ko_{X_n}}). Suppose that for every $\bf P$ such that \m{\kl_{n,\bf P}} is 
not empty, there exists a fine moduli space \m{\Pic^{\bf P}(X_n)} for 
\m{\kl_{n,\bf P}}, which is a smooth variety, defined by an open cover 
\m{(P_i)_{i\in I}} of \m{\Pic^{\bf P}(X_n)}, and for every \m{i\in I}, a line 
bundle \m{\L_i} on \m{X_n\times P_i}.
 
The component \m{\Pic^{{\bf P}_0}(X_n)} is an algebraic group. The 
multiplication with an  element of \m{B\in\kl_{n,\bf P}} defines an isomorphism 
\ \m{\eta_B:\Pic^{{\bf P}_0}(X_n)\to\Pic^{\bf P}(X_n)}.

Suppose that \m{X_n} can be extended to a primitive multiple scheme \m{X_{n+1}} 
of multiplicity \m{n+1}. From \ref{ext_B1}, we have a map \ 
\m{\Delta:\Pic(X_n)\to H^2(X,L^n)} \ such that \ \m{\kb\in\Pic(X_n)} \ can be 
extended to a line bundle on \m{X_{n+1}} if and only if \ \m{\Delta(\kb)=0}.
Suppose that some \m{\kb\in\kl_{n,\bf P}} can be extended to \m{X_{n+1}}. The 
set of line bundles in \m{\kl_{n,\bf P}} that can be extended to \m{X_{n+1}} is 
a smooth closed subvariety
\[\Gamma^{\bf P}(X_{n+1}) \ \subset \ \Pic^{\bf P}(X_n) \ .\]
We have \ \m{\Gamma^{\bf P}(X_{n+1})=\eta_\kb\big(\Gamma^{{\bf 
P}_0}(X_{n+1})\big)}, and \m{\Gamma^{{\bf P}_0}(X_{n+1})} is a subgroup of
\m{\Pic^{{\bf P}_0}(X_n)}: from \ref{ext_B1} and 
\ref{CCVB}, $(\ref{equ21})$, \m{\Gamma^{\ko_X}(X_{n+1})} is the kernel of the 
morphism of groups
\[\Delta:\Pic^{{\bf P}_0}(X_n)\to H^2(X,L^n) \ , \]
(the {\em obstruction morphism}).

Let $Z$ be a scheme and \m{\L} a line bundle on \m{X_{n+1}\times Z}. Let 
\m{z\in Z} be a closed point. Then there is an open neighbourhood $U$ of $z$, 
\m{i\in I}, and a morphism \ \m{\phi:U\to P_i}, such that \ \m{\L_{|X_n\times 
U}\simeq\phi^\sharp(\L_i)}.

\sepprop

\begin{subsub}\label{lem11}{\bf Lemma:} We have \ 
\m{\imm(\phi)\subset\Gamma^{\bf P}(X_{n+1})} .
\end{subsub}
\begin{proof} Let \m{p_X:X\times U\to X}, \m{p_U:X\times U\to U} be the 
projections. As in \cite{dr10}, 7.1, we can associate to every line bundle 
\m{\kl_0} on \m{X_n\times U} an element \m{\Delta(\kl_0)} of \m{H^2(X\times 
U,p_X^*(L^n))} such that \m{\kl_0} can be extended to a line bundle on \ 
\m{X_{n+1}\times U} \ if and only if \ \m{\Delta(\kl_0)=0}. So we have \ 
\m{\delta(\L_{|X_{n+1}\times U})=0}. Consider the 
canonical morphism
\[\rho:H^2(X\times U,p_X^*(L^n))\lra H^0(U,R^2p_{U*}(p_X^*(L^n))) \ = \ 
H^0(\ko_U)\ot H^2(X,L^n) \ . \]
As in Proposition \ref{prop2}, \m{\rho(\Delta(\kl_0))=\Delta\circ\phi}. Hence 
\m{\Delta\circ\phi=0}, i.e. \m{\imm(\phi)\subset\Gamma^{\bf P}(X_{n+1})} .
\end{proof}

\sepprop

From Theorem \ref{theo11} we have

\sepprop

\begin{subsub}\label{theo8}{\bf Theorem: } Suppose that \m{\kl_{n+1,{\bf P}}} 
is nonempty. Then there exists a fine moduli space \m{\Pic^{\bf P}(X_{n+1})} 
for \m{\kl_{n+1,{\bf P}}}. It is a smooth irreducible variety, and an affine 
bundle over \m{\Gamma^{\bf P}(X_{n+1})} with associated vector bundle \ 
\m{\ko_{\Gamma^{\bf P}(X_{n+1})}\ot\big(H^1(X,L^n)/\imm(\delta^0)\big)} .
\end{subsub}

\end{sub}

\sepsub

\Ssect{The kernel of the obstruction morphism}{Kern_ob_m}

As in \ref{PMS} we assume that we have an affine cover 
\m{(U_i)_{i\in I}} of $X$ and trivializations \Nligne \m{\delta_i:U_i^{(n)}\to 
U_i\times{\bf Z}_n}. We can write
\[(\delta_{ij}^*)_{|\ko_X(U_{ij})} \ = \ 
I_{|\ko_X(U_{ij})}+tD_{ij}+t^2\Phi_{ij} \ , \]
where \m{D_{ij}} is a derivation of \m{\ko_X(U_{ij})}. 

Suppose that \m{X_n} can be extended to a primitive multiple 
scheme \m{X_{n+1}} of multiplicity \m{n+1}. Let \m{(\delta'_{ij})} be a family 
which defines \m{X_{n+1}}, \m{\delta'_{ij}} being an automorphism of 
\m{U_{ij}\times{\bf Z}_{n+1}} inducing \m{\delta_{ij}}. We have then
\[\delta'_{ij}(t) \ = \ \alpha_{ij}t \ , \]
for some \m{\alpha_{ij}\in[\ko_X(U_{ij})[t]/(t^n)]^*}.

Let $\D$ be a line bundle on 
\m{X_n} that can be extended to a line bundle on \m{X_{n+1}}. According to 
\ref{ext_B1} this is equivalent to  \m{\Delta(\D)=0} (cf. \ref{ext_B1}). Now 
let \m{\D'} be a line bundle on \m{X_n} such that \ 
\m{\D'_{|X_{n-1}}\simeq\D_{|X_{n-1}}}. Then by formula $(\ref{equ3})$, \m{\D'} 
can be extended to a line bundle on \m{X_{n+1}} if and only if
\[\ov{\sigma}_{X_{n+1}}(\nabla_0(\D')-\nabla_0(\D)) \ = \ 0 \ . \]
Suppose that $\D$ is defined by a cocycle \m{(\theta_{ij})}, with \ 
\m{\theta_{ij}\in[\ko_X(U_{ij})(t)/(t^n)]^*}, which is the set of elements \ 
\m{a\in\ko_X(U_{ij})(t)/(t^n)} such that \ \m{a^{(0)}\in\ko_X(U_{ij})} \ 
is invertible. The cocycle relation is
\[\theta_{ik} \ = \ \theta_{ij}\delta_{ij}^*(\theta_{jk})\]
(cf. \ref{const_sh}).
We can assume that \m{\D'} is defined by a cocycle 
\m{(\theta'_{ij})}, with \m{\theta'_{ij}} of the form
\[\theta'_{ij} \ = \ \theta_{ij}+\rho_{ij}t^{n-1} \ , \]
with \ \m{\rho_{ij}\in\ko_X(U_{ij})}. Let \ 
\m{\dsp\beta_{ij}=\frac{\rho_{ij}}{\theta^{(0)}_{ij}}} . The cocycle relation 
for \m{(\theta'_{ij})} is equivalent to
\begin{equation}\label{equ17}\beta_{ik} \ = \ 
\beta_{ij}+(\alpha^{(0)}_{ij})^{n-1}\beta_{jk} \ , \end{equation}
i.e. \m{(\beta_{ij})} represents an element $\beta$ of \m{H^1(X,L^{n-1})}.

The line bundles $\D$, \m{\D'} are extensions to \m{X_n} of a line bundle 
\m{\D_{n-1}} on \m{X_{n-1}}. Let \m{D=\D_{n-1|X}}. We have canonical exact 
sequences
\[0\lra D\ot L^{n-1}\lra\D\lra\D_{n-1}\lra 0 \ , \quad
0\lra D\ot L^{n-1}\lra\D'\lra\D_{n-1}\lra 0 \ , \]
corresponding to \ \m{\sigma,\sigma'\in\Ext^1_{\ko_{X_n}}(\D_{n-1},D\ot 
L^{n-1})} \ respectively. From \ref{ext_4} we have an exact sequence
\xmat{0\ar[r] & H^1(X,L^{n-1})\ar[r]^-\Psi & \Ext^1_{\ko_{X_n}}(\D_{n-1},D\ot 
L^{n-1})\ar[r] & \C \ ,}
and \ \m{\sigma'-\sigma=\Psi(\beta)}.

Now \ \m{\nabla_0(\D')-\nabla_0(\D)} \ is represented by the cocycle 
\m{\dsp\bigg(\frac{d\theta'_{ij}}{\theta'_{ij}}-\frac{d\theta_{ij}}{\theta_{ij}}
\bigg)} (cf. \ref{CCVB}). We have
\[\frac{d\theta'_{ij}}{\theta'_{ij}}-\frac{d\theta_{ij}}{\theta_{ij}} \ = \
(n-1)\beta_{ij}t^{n-2}dt+t^{n-1}d\beta_{ij} \ . \]
Let \ \m{\L=\ki_{X,X_{n+1}|X_2}}, which is a line bundle on \m{X_2}.
Recall that \m{\delta^1_{\L^{n-1}}} denotes \Nligne the connecting morphism \ 
\m{H^1(X,L^{n-1})\to H^2(X,L^n)} \ coming from the exact sequence \ \m{0\to 
L^n\to\L^{n-1}\to L^{n-1}\to 0} \ on \m{X_2}.

\sepprop

\begin{subsub}\label{theo5}{\bf Theorem:} We have \
\m{\ov{\sigma}_{X_{n+1}}(\nabla_0(\D')-\nabla_0(\D))=n\delta^1_{\L^{n-1}}(\beta)}
 . \end{subsub}
\begin{proof}
The line bundle \m{\L^{n-1}} is defined by the family 
\m{\big((\alpha_{ij}^{(0)}+\alpha_{ij}^{(1)}t)^{n-1}\big)}. We have \
\[(\alpha_{ij}^{(0)}+\alpha_{ij}^{(1)}t)^{n-1} \ = \
\big(\alpha_{ij}^{(0)}\big)^{n-1}+(n-1)
\big(\alpha_{ij}^{(0)}\big)^{n-2}\alpha_{ij}^{(1)}t \ . \]
Hence from Lemma \ref{lem10}, \m{n\delta^1_{\L^{n-1}}(\beta)} is represented by 
\m{(\mu_{ijk})}, \m{\mu_{ijk}\in\ko_X(U_{ijk})[t]/(t^n)}, with
\[\mu_{ijk} \ = \ 
n\big(\alpha^{(0)}_{ij}\big)^{n-2}\big((n-1)\alpha^{(1)}_{ij}\beta_{jk}+
\alpha^{(0)}_{ij}D_{ij}(\beta_{jk})\big) \ . \]
We consider the exact sequence
\[0\lra L^n\lra\Omega_{X_{n+1}|X_n}\lra\Omega_{X_n}\lra 0 \ , \]
and the associated map
\[\delta:H^1(X_n,\Omega_{X_n})\lra H^2(X,L^n) \ . \]
We construct \m{\Omega_{X_{n+1}}} as in \ref{PMS_Omeg}. Let 
\[A_{ij} \ = \ (n-1)\beta_{ij}t^{n-2}dt+t^{n-1}d\beta_{ij} \ . \]
Then we have
\[\ov{\sigma}_{X_{n+1}}(\nabla_0(\D')-\nabla_0(\D)) \ = \ \delta(\nabla_0(\D')
-\nabla_0(\D)) \ , \]
and it is represented by the family \m{(\nu_{ijk|X_n})}, \m{\nu_{ijk}\in 
H^0(\Omega_{U_{ijk}\times{\bf Z}_{n+1}})}, where
\[\nu_{ijk} \ = \ A_{ij}-A_{ik}+{\bf d}_{ij}(A_{jk}) \ . \]
A long computation, using the equalities $(\ref{equ17})$, 
\m{D_{ik}=D_{ij}+\alpha_{ij}^{(0)}D_{jk}} \ and
\[\alpha_{ik}^{(1)} \ = \ \alpha_{ij}^{(0)}D_{ij}(\alpha_{jk}^{(0)})+
\big(\alpha_{ij}^{(0)}\big)^2\alpha_{jk}^{(1)}+\alpha_{jk}^{(0)}\alpha_{ij}^{(1)}
\ , \]
shows that \ 
\m{\nu_{ijk}=n\big(\alpha^{(0)}_{ij}\big)^{n-2}\big((n-1)\alpha^{(1)}_{ij}
\beta_{jk}+\alpha^{(0)}_{ij}D_{ij}(\beta_{jk})\big)t^{n-1}dt}, which, by Lemma 
\ref{lem10} and \ref{cech2e}, implies Theorem \ref{theo5}.
\end{proof}

\sepprop

Let $\bf P$ be an irreducible component of \m{\Pic(X)}. Recall that 
\m{\Gamma^{\bf P}(X_{n+1})} is the smooth subvariety of \m{\Pic^{\bf P}(X_n)} 
of line bundles which can be extended to \m{X_{n+1}}. The following result is 
an easy consequence of Theorem \ref{theo5}:

\sepprop

\begin{subsub}\label{coro5}{\bf Proposition:} Let $Z$ be the image of he 
restriction morphism \ \m{\Gamma^{\bf P}(X_{n+1})\to\Gamma^{\bf P}(X_n)}. Then 
$Z$ is smooth and \m{\Gamma^{\bf P}(X_{n+1})} is an affine bundle over $Z$ with 
associated vector bundle \ \m{\ko_Z\ot\ker(\delta^1_{\L^{n-1}})} .
\end{subsub}

\end{sub}

\sepsec

\section{Products of curves}\label{pro_cur}

Let $C$ and $D$ be smooth projective irreducible curves, \m{X=C\times D} and
\ \m{\pi_C:X\to C}, \m{\pi_D:X\to D} \ the projections. Let \m{g_C} (resp. 
\m{g_D}) be the genus of $C$ (resp. $D$). Suppose that \m{g_C\geq 2}, 
\m{g_D\geq 2}.

\sepsub

\Ssect{Basic properties}{bas_pro}

\begin{subsub}\label{bas_pro1} -- \rm
Let \m{L_C} (resp. \m{L_D}) be a line bundle on $C$ (resp. $D$) and \ 
\m{L=\pi_C^*(L_C)\ot\pi_D^*(L_D)}. Then we have canonical isomorphisms
\begin{eqnarray*}
{\rm (i)} \ \ \ H^0(X,L) & \simeq & H^0(C,L_C)\ot H^0(D,L_D) ,\\
{\rm (ii)} \ \ H^1(X,L) & \simeq & \big(H^0(X,L_C)\ot H^1(D,L_D)\big)\oplus
\big(H^1(C,L_C)\ot H^0(D,L_D)\big) ,\\
{\rm (iii)} \ H^2(X,L) & \simeq & H^1(C,L_C)\ot H^1(D,L_D) .
\end{eqnarray*}
\end{subsub}

\sepsubsub

\begin{subsub}\label{bas_pro2} -- \rm
Let \m{M_C} (resp. \m{M_D}) be a line bundle on $C$ (resp. $D$) and \ 
\m{M=\pi_C^*(M_C)\ot\pi_D^*(M_D)}. 

Using \ref{cech2e}, \ref{Ex_Ler}, and the direct sum decompositions of (ii) and 
(iii), it is easy to see that the only components of the canonical map \ 
\m{H^1(X,L)\ot H^1(X,M)\to H^2(X,L\ot M)} \ that could not vanish are the two 
canonical maps
\[\big(H^1(L_C)\ot H^0(L_D)\big)\ot\big(H^0(M_C)\ot H^1(M_D)\big)\lra
H^1(L_C\ot M_C)\ot H^1(L_D\ot M_D) \ , \]
\[\big(H^0(L_C)\ot H^1(L_D)\big)\ot\big(H^1(M_C)\ot H^0(M_D)\big)\lra
H^1(L_C\ot M_C)\ot H^1(L_D\ot M_D) \ , \]
induced by
\[H^1(L_C)\ot H^0(M_C)\to H^1(L_C\ot M_C) \ , \quad
H^0(L_D)\ot H^1(M_D)\to H^1(L_D\ot M_D) \ , \]
and
\[H^0(L_C)\ot H^1(M_C)\to H^1(L_C\ot M_C) \ , \quad
H^1(L_D)\ot H^0(M_D)\to H^1(L_D\ot M_D)\]
respectively.
\end{subsub}

\sepsubsub

\begin{subsub}\label{bas_pro3} -- \rm Similarly, the canonical map
\[H^1(X,\pi_C^*(L_C\ot M_C)\ot
\pi_D^*(L_D))\ot H^1(C,\pi_C^*(M_C^*))\to H^2(X,L)\]
comes from the two canonical maps
\[\big(H^0(L_C\ot M_C)\ot H^1(L_D)\big)\ot H^1(M_C^*)\lra H^1(L_C)\ot 
H^1(L_D) \ , \]
\[\big(H^1(L_C\ot M_C)\ot H^0(L_D)\big)\ot\big(H^0(M_C^*)\ot 
H^1(\ko_D)\big)\lra H^1(L_C)\ot H^1(L_D)\]
induced by \ \m{H^0(C,L_C\ot M_C)\ot H^1(C,M_C^*)\to H^1(C,L_C)} \ and
\[H^1(L_C\ot M_C)\ot H^0(M_C^*)\lra H^1(L_C) \ , \quad
H^0(L_D)\ot H^1(\ko_D)\lra H^1(L_D)\]
respectively.
\end{subsub}

\sepsubsub

\begin{subsub}\label{bas_pro4} -- \rm We have \ 
\m{\Omega_X\simeq\pi_C^*(\omega_C)\oplus\pi_D^*(\omega_D)}, hence \
\m{\omega_X\simeq\pi_C^*(\omega_C)\ot\pi_D^*(\omega_D)} and
\[H^1(X,\Omega_X) \ = \ \big(H^0(\omega_C)\ot H^1(\ko_D)\big)\oplus
H^1(\omega_C)\oplus\big(H^1(\ko_C)\ot H^0(\omega_D)\big)\oplus
H^1(\omega_D) \ . \]
With respect to this decomposition, the canonical class of $L$ is \ 
\m{\nabla_0(L)=\nabla_0(L_C)+\nabla_0(L_D)} .

We have \ \m{T_X\ot L=(\pi_C^*(T_C)\ot L)\oplus(\pi_D^*(T_D)\ot L)}, and
\begin{equation}\label{equ19}H^1(X,T_X\ot L) \ = \ \big(H^0(T_C\ot L_C)\ot 
H^1(L_D)\big)\oplus
\big(H^1(T_C\ot L_C)\ot H^0(L_D)\big)\oplus\end{equation}
\[ \qquad\qquad\qquad\qquad\quad
\big(H^0(L_C)\ot H^1(T_D\ot L_D)\big)\oplus
\big(H^1(L_C)\ot H^0(T_D\ot L_D)\big) \ . \]

\end{subsub}

\end{sub}

\sepsub

\Ssect{Primitive double schemes}{prod_2}

\begin{subsub}\label{ext_L} Line bundles on \m{X_2} -- \rm
We use the notations of \ref{bas_pro}. The non trivial primitive double schemes 
with associated smooth variety $X$ and associated line bundle $L$ are 
parametrized by \m{\P(H^1(X,T_X\ot L))}. Let \ \m{\eta\in H^1(X,T_X\ot L)} \ 
and \m{X_2} the associated double scheme. Using the decomposition 
$(\ref{equ19})$, we see that $\eta$ has four components: 
\m{\eta=(\eta_1,\eta_2,\eta_3,\eta_4)} .
\end{subsub}

We denote by \m{\phi_{L_C}}, \m{\phi_{L_D}} the canonical maps \ \m{H^0(T_C\ot 
L_C)\ot H^1(\omega_C)\to H^1(L_C)} \ and \Nligne \m{H^0(T_D\ot L_D)\ot 
H^1(\omega_D)\to H^1(L_D)} \ respectively. 

Let \ \m{\Psi:H^1(X,T_X\ot L)\ot H^1(X,\Omega_X)\to H^2(X,L)} \ be the 
canonical morphism. Recall that $M$ can be extended to a line bundle on \m{X_2} 
if and only if \ \m{\Delta(M)=0}, where \ 
\m{\Delta(M)=\Psi(\eta\ot\nabla_0(M))} (cf. \cite{dr10}, Theorem 7.1.2). From 
\ref{bas_pro4}, we have
\begin{equation}\label{equ23}\Delta(M) \ = \ (\phi_{L_C}\ot 
I_{H^1(L_D)})(\eta_1\ot\nabla_0(M_C))+(I_{H^1(L_C)}\ot
\phi_{L_D})(\eta_4\ot\nabla_0(M_D)) \ . \end{equation}
The canonical class of a line bundle on a smooth projective curve is in fact an 
integer. For every line bundle $\kl$ on $C$ we have \ 
\m{\nabla_0(\kl)=\deg(\kl)c}, where \ \m{c=\nabla_0(\ko_C(P))}, for any \m{P\in 
C}.

\end{sub}

\sepsub

\Ssect{The case \m{\omega_{X_2}\simeq\ko_{X_2}}}{K3}

Let \m{X_2} be a primitive double scheme such that \m{(X_2)_{red}=X}, with 
associated line bundle $L$ on $X$. According to Corollary \ref{coro4}, we have 
\ \m{\omega_{X_2}\simeq\ko_{X_2}} \ if and only if \ \m{L\simeq\omega_X}. 

Suppose that \ \m{L\simeq\omega_X}. We have an isomorphism \ 
\m{H^1(\omega_C)\simeq\C} such that if $\kl$ is a line bundle of degree $d$ on 
$C$, we have \m{\nabla_0(\kl)=d} (and similarly on $D$). Let
\[\eta \ = \ (\eta_1,\eta_2,\eta_3,\eta_4) \ \in H^1(X,T_X\ot\omega_X)\]
(using the decomposition $(\ref{equ19})$). Here we have
\[\eta_1\in\C \ , \quad \eta_2\in H^1(C,\ko_C)\ot H^0(D,\omega_D) \ , \quad
\eta_3\in H^0(C,\omega_C)\ot H^1(D,\ko_D) \ , \quad \eta_4\in\C \ . \]
If \ \m{M=\pi_C^*(M_C)\ot\pi_D^*(M_D)}, 
\[\Delta(M) \ = \ \eta_1\nabla_0(M_C)+\eta_4\nabla_0(M_D) \ , \]
hence $M$ can be extended to a line bundle on \m{X_2} if and only if \ 
\m{\eta_1\deg(M_C)+\eta_4\deg(M_D)=0}.

\sepprop

\begin{subsub}\label{coro6}{\bf Corollary: } The scheme \m{X_2} is projective 
if and only if \ \m{\eta_1=\eta_4=0} \ or \m{\eta_1\eta_4<0} \ and 
\m{\dsp\frac{\eta_1}{\eta_4}} is rational.
\end{subsub}

\sepprop

We have by \ref{bas_pro1}
\[H^2(T_X\ot L^2) \ \simeq \ \big(H^1(C,\omega_C)\ot H^1(D,\omega_D^2)\big)
\oplus\big(H^1(C,\omega_C^2)\ot H^1(D,\omega_D)\big) \ = \ 0 \ . \]
Hence by \ref{ext-I}, \m{X_2} can be extended to a primitive multiple scheme of 
multiplicity 3 if and only \m{\omega_X} can be extended to a line bundle on 
\m{X_2}. This is the case if and only \ \m{(g_C-1)\eta_1+(g_D-1)\eta_4=0}, and 
then \m{X_2} is projective.

\end{sub}

\sepsub

\Ssect{An example of extensions of double primitive schemes}{ext_db}

Let $\Theta$ be a theta characteristic on $C$ such that \m{h^0(C,\Theta)>0}. We 
take \ \m{L_C=\Theta}, \m{L_D=\omega_D}. So we have
\[H^1(L) \ = \ \big(H^1(\Theta)\ot H^0(\omega_D)\big)\oplus\big(H^0(\Theta)\ot 
H^1(\omega_D)\big) \ , \quad H^2(L)=H^1(\Theta)\ot H^1(\omega_D) \ . \]
In the decomposition $(\ref{equ19})$ the first summand is 0, and
\[H^1(T_X\ot L) \ = \ \big(H^1(\Theta^{-1})\ot 
H^0(\omega_D)\big)\oplus\big(H^0(\Theta)\ot
H^1(\ko_D)\big)\oplus H^1(\Theta) \ . \]
Let \m{X_2} be a non trivial primitive double scheme, corresponding to \Nligne 
\m{\eta=(\eta_1,\eta_2\eta_3,\eta_4)\in\P(H^1(T_X\times L))}. We have 
\m{\eta_1=0}, and we will also suppose that \m{\eta_4=0}. Such double schemes 
are parametrized by \ 
\m{\P((H^1(\Theta^{-1})\ot H^0(\omega_D))\oplus(H^0(\Theta)\ot H^1(\ko_D)))}.

We have, from $(\ref{equ23})$, \m{\Delta(M)=0}. Hence every line bundle on $X$ 
can be extended to a line bundle on \m{X_2}, and \m{X_2} is projective. Such 
extensions of $M$ are parametrized by an affine space, with underlying vector 
space \m{H^1(L)}. 

\sepprop

\begin{subsub}\label{prop19}{\bf Proposition: } Suppose that \ \m{\eta_4=0} \ 
and \ \m{\eta_3\not=0}. Then there exists an extension of \m{X_2} to a 
primitive multiple scheme \m{X_3} of multiplicity 3, and \m{X_3} is projective.
\end{subsub}
\begin{proof} By Theorem \ref{theo6}, to prove the existence of \m{X_3} it 
suffices to show that if
\[\chi:H^1(T_X\ot L)\ot H^1(L)\lra H^2(T_X\ot L^2)\]
is the canonical map, then \ \m{\chi(\{\eta\}\ot H^1(L))=H^2(T_X\ot L^2)} . We 
have
\begin{eqnarray*} H^2(T_X\ot L^2) & = & H^2(p_D^*(\omega_D^2))\oplus 
H^2(p_C^*(\omega_C)\ot
p_D^*(\omega_D)) \ = \  H^2(p_C^*(\omega_C)\ot p_D^*(\omega_D))\\
& = & H^1(\omega_C)\ot H^1(\omega_D) \ \simeq \ \C \ ,
\end{eqnarray*}
and $\chi$ is equivalent to the canonical map
\[\chi':\big(H^0(\Theta)\ot H^1(\ko_D)\big)\ot\big(H^1(\Theta)\ot 
H^0(\omega_D)\big)\lra H^1(\omega_C)\ot H^1(\omega_D) \ , \]
and we have to prove that \ \m{\chi'(\{\eta_3\}\ot(H^1(\Theta)\ot 
H^0(\omega_D))=H^1(\omega_C)\ot H^1(\omega_D)} .
Let \m{(\sigma_i)_{1\leq i\leq g_D}} be a basis of \m{H^1(\ko_D)}, and 
\m{(\sigma_i^*)_{1\leq i\leq g_D}} the dual basis of \m{H^0(\omega_D)}. We can 
write \ \m{\eta_3=\sigg_{i=1}^g\mu_i\ot\sigma_i}, with \ \m{\mu_i\in 
H^0(\Theta)}. If \m{1\leq j\leq g_D} and \m{\rho\in H^1(\Theta)}, we have
\ \m{\chi'(\eta_3\ot(\rho\ot\sigma_j^*))=\rho\mu_j}. Since \m{\eta_3\not=0}, we 
can choose $j$ such that \m{\mu_j\not=0}. We have an exact sequence on $C$
\xmat{0\ar[r] & \Theta\ar[rr]^-{\times\mu_j} & & \Theta^2=\omega_C\ar[r] & 
\T\ar[r] & 0 \ ,}
where $\T$ is a torsion sheaf. It follows that \ 
\m{H^1(\times\mu_j):H^1(\Theta)\to H^1(\omega_C)} \ is surjective. Since \ 
\m{H^2(L^2)=\nsp}, every line bundle on \m{X_2} can be extended to 
\m{X_3}. This proves Proposition \ref{prop19}.
\end{proof}

\end{sub}

\sepsub

\Ssect{Picard groups}{pic_b}

We use the notations of \ref{ext_L}.
Let \m{M_C\in\Pic(C)}, \m{M_D\in\Pic(D)}, \m{M=\pi_C^*(M_C)\ot\pi_D^*(M_D)} be 
such that \ \m{\Delta(M)=0}, so that $M$ can be 
extended to a line 
bundle on \m{X_2}. Let $\bf P$ be the irreducible component of \m{\Pic(C\times 
D)} that contains \m{\pi_C^*(M_C)\ot\pi_D^*(M_D)}. Let \ \m{\Pic^{\bf 
P}(X_2)\subset\Pic(X_2)} \ be the irreducible component of line bundles $\kl$ 
such that \m{\kl_{|X}\in{\bf P}}. From \ref{pic}, it is an affine bundle 
over $\bf P$, with associated vector bundle \ \m{\ko_{\bf P}\ot H^1(X,L)}.

Let \m{C'=C}, \m{X=C\times D}, and \m{p_{C}:X\times C'\to 
C}, \m{p_{C'}:X\times C'\to C'}, \m{p_D:X\times C'\to 
D},\Nligne \m{p_{X}:X\times C'\to X}, \m{p_{C\times C}:X\times 
C'\to C\times C'} the projections.

Let \ \m{n=\deg(M_C)} \ and \ \m{\Delta_C\subset C\times C} \ the diagonal. 
Then \ \m{\L=p_{C\times C}^*(\ko(n\Delta_C))\ot p_D^*(M_D)} \ is a family of 
lines bundles on $X$ parametrized by \m{C'}. These line bundles 
belong to $\bf P$, and $\L$ defines an immersion \m{C'\to\Pic^{\bf P}(X)}. 

We have
\[\Delta(\L) \ = \ \Delta(p_{C\times C}^*(\ko(n\Delta_C)))+
\Delta(p_D^*(M_D)) \ , \]
Let \ \m{q_{C\times C}:p_{C\times C}^*(\Omega_{C\times C})\to\Omega_{X\times 
C'}}, \m{q_D:p_D^*(\omega_D)\to \Omega_{X\times C'}} \ 
be the injective canonical morphisms. We have
\[\nabla_0(p_{C\times C}^*(\ko(n\Delta_C))) \ = \ H^1(q_{C\times 
C})(\nabla_0(\ko(n\Delta_C))) 
\ , \quad \nabla_0(p_D^*(M_D)) \ = \ H^1(q_D)(\nabla_0(M_D)) \ . \]
We have \ \m{H^1(X\times C',\Omega_{X\times 
C'})=A_{C}\oplus A_{C'}\oplus A_D}, with
\begin{eqnarray*}A_{C} & = & 
H^1(\omega_{C})\oplus\big(H^0(\omega_{C})\ot(H^1(\ko_D)\oplus 
H^1(\ko_{C'}))\big) \ ,\\
A_{C'} & = & 
H^1(\omega_{C'})\oplus\big(H^0(\omega_{C'})\ot(H^1(\ko_D)\oplus 
H^1(\ko_{C}))\big) \ ,\\
A_{D} & = & 
H^1(\omega_{D})\oplus\big(H^0(\omega_{D})\ot(H^1(\ko_{C})\oplus 
H^1(\ko_{C'}))\big) \ ,\\
\end{eqnarray*}
and \ \m{\imm(H^1(q_D)) \ = \ H^1(\omega_D)},
\[\imm(H^1(q_{C\times C})) \ = 
H^1(\omega_{C})\oplus(H^0(\omega_{C})\ot H^1(\ko_{C'}))\oplus 
H^1(\omega_{C'})\oplus(H^0(\omega_{C'})\ot H^1(\ko_{C})) \ . \]
We have \ \m{\nabla_0(p_D^*(M_D))=\nabla_0(M_D)}, and by \ref{th4_pro},
 \Nligne \m{q_{C\times 
C}(\nabla_0(\ko(n\Delta_C)))=n(1,I_{H^0(\omega_C)},-1,-I_{H^0(\omega_C)})}.

We see \ \m{X_2\times C'} as a primitive multiple scheme, with underlying 
smooth scheme \m{X\times C'} and associated line bundle \m{p_X^*(L)}. We have
\[H^1(T_{X\times C'}\ot p_X^*(L)) \ \simeq \ H^1(T_X\ot L)\oplus
(H^0(T_X\ot L)\ot H^1(\ko_{C'}))\oplus B \ , \]
with
\[B \ = \ (H^1(L_{C})\ot H^0(L_D)\ot H^0(\omega_{C'}))\oplus\qquad 
\qquad\qquad\qquad\quad \]
\[\qquad\qquad\qquad\qquad (H^0(L_{C})\ot H^1(L_D)\ot H^0(\omega_{C'}))\oplus
(H^0(L_{C})\ot H^0(L_D)\ot H^1(\omega_{C'})) ,\]
The element \m{\eta'} of \m{H^1(T_{X\times C'}\ot p_X^*(L))} associated to 
\m{X_2\times C'} is $\eta$ (in the factor \m{H^1(T_X\ot L)}).

We have
\[H^2(p_X^*(L)) \ = \ (H^1(L_C)\ot 
H^1(L_D))\oplus\qquad\qquad\qquad\qquad\qquad\qquad\qquad\qquad\qquad\]
\[\qquad\qquad\qquad(H^1(L_C)\ot H^0(L_D)\ot H^1(\ko_{C'}))\oplus
(H^0(L_C)\ot H^1(L_D)\ot H^1(\ko_{C'})) \ . \]

Let \ \m{\beta_{L_C}:H^0(\omega_C)\ot H^1(T_C\ot L_C)\to H^1(L_C)} \ be the 
canonical map. As in \ref{ext_L} we can compute \m{\Delta(\L)} and find that it 
has the component \ \m{\Delta(M)=0} \ in the factor \ \m{H^1(L_C)\ot 
H^1(L_D)},\break and the component \ \m{(\beta_{L_C}\ot I_{H^1(\omega_C)\ot 
H^0(L_D)})(\eta_2\ot I_{H^0(\omega_C)})} \ in the factor \Nligne \m{H^1(L_C)\ot 
H^0(L_D)\ot H^1(\ko_{C'})}. It follows that

\sepprop

\begin{subsub}\label{theo16}{\bf Theorem: } If \ \m{n(\beta_{L_C}\ot 
I_{H^1(\omega_C)\ot H^0(L_D)})(\eta_2\ot I_{H^0(\omega_C)})\not=0}, then the 
affine bundle \m{\Pic^{\bf P}(X_2)} is not a vector bundle over $\bf P$.
\end{subsub}

\sepprop

\begin{subsub}\label{K3_2} The case \m{L=\omega_{C\times D}} -- \rm We have 
\m{\eta_1,\eta_4\in\C}, 
\m{\eta_2\in H^1(C,\ko_C)\ot H^0(D,\omega_D)}, and \ \m{\Delta(M)=0} \ is 
equivalent to \m{\eta_1\deg(M_C)+\eta_4\deg(M_D)=0},. Theorem \ref{theo16} 
implies 
that if \ \m{\deg(M_C)\not=0} \ and \m{\eta_2\not=0}, then \m{\Pic^{\bf 
P}(X_2)} is not a vector bundle over $\bf P$.
\end{subsub}

\end{sub}

\sepsec

\section{Moduli spaces of vector bundles on primitive multiple 
curves}\label{pic_rib0}

Let $C$ be a smooth irreducible projective curve of genus \m{g\geq 2}, and 
\ \m{Y=C_2} \ a {\em ribbon} over $C$, i.e. a primitive multiple curve of 
multiplicity 2 such that \m{Y_{red}=C}. Let \ \m{L=\ki_{C,Y}} \ be the 
associated line bundle on $C$.

\sepsub

\Ssect{Construction of the moduli spaces}{rib}

We suppose that \ \m{\deg(L)<0}.
Let $r$, $d$ be coprime integers such that \m{r>0}. Let \m{M_C(r,d)} be 
the moduli space of stable sheaves on $C$ of rank $r$ and degree $d$. It is a 
smooth projective variety of dimension \ \m{r^2(g-1)+1}. On \ \m{M_C(r,d)\times 
C} \ there is an {universal bundle} $E$. For every closed point \m{m\in 
M_C(r,d)}, \m{E_m} is the stable vector bundle corresponding to $m$.

Let \ \m{p_M:M_C(r,d)\times C\to M_C(r,d)}, \m{p_C:M_C(r,d)\times C\to C} \ be 
the projections. The vector bundle of \ref{affin},
\[\A \ = \ R_1p_{M*}\big(\HHom(E,E\ot p_C^*(L))\big)\]
is of rank \ \m{r^2(\deg(L)+g-1)}. For every \ \m{m\in M_C(r,d)}, 
\m{\A_m=\Ext^1_{\ko_C}(E_m,E_m\ot L)}.

The conditions of \ref{mod} are satisfied. Hence we obtain a fine moduli space 
\m{{\bf M}_Y(r,d)} for the vector bundles $\E$ on $Y$ such that \ \m{\E_{|C}\in 
M_C(r,d)}. We have a canonical morphism
\[\tau_{r,d}:{\bf M}_Y(r,d)\lra M(r,d)\]
associating \m{\E_{|C}} to $\E$, and \m{{\bf M}_Y(r,d)} is an affine bundle on 
\m{M(r,d)} with associated vector bundle $\A$.

The same kind of construction can be applied in higher multiplicity to obtain 
Theorems \ref{theo_0_3} and \ref{theo_0_4}.

\end{sub}

\sepsub

\Ssect{Picard groups}{pic_rib}

We don't suppose here that \m{\deg(L)<0}. We will use the results of \ref{pic}.
For \m{r=1}, we have \ \m{M(r,d)=\Pic^{d}(C)}, the variety of 
line bundles of degree $d$ on $C$. As on $C$ we have a decomposition of the 
Picard group of $Y$:
\[\Pic(Y) \ = \ \bigcup_{d\in\Z}\Pic^{d}(Y) \ , \]
where \ \m{\Pic^{d}(Y)={\bf M}_Y(1,d)} \ is the variety of line bundles on 
$Y$ whose restriction to $C$ has degree $d$. With the notations of \ref{pic}, 
\m{\Pic^{d}(Y)=\Pic^{\bf P}(Y)}, where $\bf P$ is the set of line bundles of 
degree $d$ on $C$. Consider the restriction morphism
\[\tau_{1,d}:\Pic^{d}(Y)\lra \Pic^{d}(C)\]
which is an affine bundle over \m{\Pic^{d}(C)}, with associated line bundle \ 
\m{\ko_{\Pic^{d}(C)}\ot H^1(C,L)} .

\sepprop

Let \ \m{0\to L\to\Omega_{Y|C}\to\omega_C\to 0} \ be the canonical exact 
sequence, associated to \Nligne \m{\sigma_Y\in H^1(C,T_C\ot L)}.
Let \ \m{p_Y:Y\times C\to Y} \ be the projection. Let
\[\Sigma_{Y\times C}:\qquad
0\lra p_Y^*(L)\lra\Omega_{Y\times C|C\times C}\lra\Omega_{C\times C}\lra 0 \ , 
\]
be the exact sequence $(\ref{equ1})$ corresponding to \m{Y\times C},
associated with
\[\sigma_{Y\times C} \ \in \ \Ext^1_{\ko_{C\times 
C}}(p_1^*(\omega_{C}),p_1^*(L))\oplus\Ext^1_{\ko_{C\times 
C}}(p_2^*(\omega_C),p_1^*(L)) \ . \]
Then \m{\Sigma_{Y\times C}} is \m{p_1^*(\Sigma_Y)\oplus p_2^*(\Sigma)}, 
where $\Sigma$ is the exact sequence
\xmat{0\ar[r] & 0\ar[r] & \omega_C\fleq[r] & \omega_S\ar[r] & 0 \ .}
It follows that the component of \m{\sigma_{Y\times C}} in 
\m{\Ext^1_{\ko_{C\times C}}(p_{2}^*(\omega_C),p_1^*(L))} vanishes, and that 
\[\sigma_{Y\times C} \ = \ \sigma_Y \in H^1(C,T_C\ot 
L)\subset\Ext^1_{\ko_{C\times C}}(p_1^*(\omega_{C}),p_1^*(L))\]
(cf. $(\ref{equ10})$).

According to Theorem \ref{theo3} it is equivalent to say that \m{\Pic^{d}(Y)} 
is not a vector bundle over \m{\Pic^{d}(C)}, and that a Poincar\'e 
bundle $\D$ on \ \m{C\times\Pic^{d}(C)} \ cannot be extended to a line 
bundle of \ \m{Y\times\Pic^{d}(C)}. To see this we can suppose that \m{d=-1}. 
We have then an embedding
\[\xymatrix@R=5pt{C\flinc[r] & \Pic^{-1}(C)\\ P\fmaps[r] & \ko_C(-P) \ .} \]
We will show that with suitable hypotheses, \m{\D_{|C\times C}} cannot be 
extended to \m{Y\times C}.
Let \ \m{p_1:C\times C\to C} (resp. \m{p_2:C\times C\to C}) be the 
first (resp. second) projection. It suffices to prove that for some \ 
\m{N\in\Pic(C)}, \m{\D_{|C\times C}\ot p_2^*(N)} cannot be extended to 
\m{Y\times C}. So we can assume that \ \m{\D_{|C\times C}=\ko_{C\times 
C}(-\Gamma)}, where \ \m{\Gamma\subset C\times C} \ is the diagonal.

We see \m{Y\times C} as a primitive multiple scheme, with associated line 
bundle \m{p_1^*(L)} on \m{C\times C}. We will show that \ 
\m{\Delta(\ko_{C\times C}(-\Gamma))\in H^2(C\times C,p_1^*(L))} \ is non zero. 
This implies that \m{\ko_{C\times C}(-\Gamma)} cannot be extended to \m{Y\times 
C} (cf. \ref{ext_B1}).

We have
\[\Delta(\ko_{C\times C}(-\Gamma)) \ \in \ H^2(C\times 
C,p_1^*(L)) \ = \ L(H^0(C,\omega_C),H^1(C,L)) \ , \]
and \ \m{\Delta(\ko_{C\times C}(-\Gamma))=\sigma_{Y\times 
C}.\nabla_0(\ko_{C\times C}(-\Gamma))} (from \cite{dr10}, Theorem 7.1.2).
It follows easily from \ref{bas_pro2} that \m{\Delta(\ko_{C\times C}(-\Gamma))} 
is the image of \m{\sigma_Y} by the canonical map
\[\eta:H^1(C,T_C\ot L)\lra L(H^0(C,\omega_C),H^1(C,L)) \ . \]
It follows that if \ \m{\eta(\sigma_Y)\not=0}, then \m{\Pic^{d}(Y)} is not a 
vector bundle over \m{\Pic^{d}(C)}.

We will prove

\sepprop

\begin{subsub}{\bf Theorem :}\label{theo4} 
\m{\Pic^{d}(Y)} is not a vector bundle over \m{\Pic^d(C)} if $Y$ is not trivial 
and either $C$  is not hyperelliptic and \ \m{\deg(L)\leq 2-2g}, or 
\m{L=\omega_C}.
\end{subsub}
\end{sub}

\sepsub

\Ssect{The duality morphism}{coh_prod}

We will use the results of \ref{bas_pro}.

Let $Z$ be smooth irreducible projective curve. Let \ \m{p_Z:Z\times C\to Z} \ 
and \ \m{p_C:Z\times C\to C} \ be the projections. Let $N$ be a vector bundle 
on $C$. Then we have (from \ref{bas_pro})
\begin{equation}\label{equ10}
H^1(Z\times C, p_C^*(N)) \ \simeq \ H^1(C,N)\oplus\big(H^0(C,N)\ot 
H^1(Z,\ko_Z)\big) \ . 
\end{equation}
We have also \ \m{H^2(Z\times C,p_C^*(N))\simeq H^1(C,N)\ot H^1(Z,\ko_Z)} .

\sepprop

\begin{subsub}\label{can_b}Canonical bundles -- \rm We have \ 
\m{\Omega_{Z\times C}=p_Z^*(\omega_Z)\oplus p_C^*(\omega_C)}, hence
\begin{equation}\label{equ11}H^1(Z\times C,\Omega_{Z\times C}) \ \simeq \ 
H^1(Z,\omega_Z)\oplus\big(H^0(Z,\omega_Z)\ot 
H^1(C,\ko_C)\big)\oplus\quad\qquad\qquad\qquad\quad\quad\end{equation}
\[ \qquad H^1(C,\omega_C)\oplus\big(H^0(Z,\ko_Z)\ot 
H^1(C,\omega_C)\big)\]
\[ \quad\quad\qquad\qquad\quad\qquad = \ \C\oplus 
L(H^0(C,\omega_C),H^0(Z,\omega_Z))
\oplus\C\oplus L(H^0(Z,\omega_Z),H^0(C,\omega_C)) \ . \]
We have \ \m{H^2(Z\times C,\Omega_{Z\times C})\simeq H^1(Z,\omega_Z)\ot
H^1(C,\omega_C)=\C} . 
\end{subsub}

\sepprop

\begin{subsub}\label{15}{\bf Lemma:} The duality morphism
\[H^1(Z\times C,\Omega_{Z\times C})\times H^1(Z\times C,\Omega_{Z\times C})\lra
H^2(Z\times C,\omega_{Z\times C})=\C\]
is
\xmat{\big((\alpha_1,\phi_1,\alpha_2,\phi_2),(a_1,f_1,a_2,f_2)\big)\fmaps[rr]
& & \alpha_1a_2-\alpha_2a_1+tr(f_2\phi_1)-tr(f_1\phi_2) \ . }
\end{subsub}
\begin{proof}
This follows from \ref{bas_pro}. For the minus signs, see \ref{cech2e}.
\end{proof}

\end{sub}

\sepsub

\Ssect{Proof of Theorem \ref{theo4}}{th4_pro}

\begin{subsub}The canonical class of \m{\ko_{C\times C}(-\Gamma)} -- \rm We 
have \ \m{\nabla_0(\ko_{C\times C}(-\Gamma))=-{\rm cl}_\Gamma}, where \m{{\rm 
cl}_\Gamma} is the cohomology class of $\Gamma$. It is the element of 
\m{H^1(C\times C,\Omega_{C\times C})} corresponding to the canonical linear form
\[\tau:H^1(C\times C,\Omega_{C\times C})\lra 
H^1(\Gamma,\omega_\Gamma)=H^1(C,\omega_C)=\C\]
using the duality described in \ref{can_b}. For every
\[(a_1,f_1,a_2,f_2) \ \in \
\C\oplus\End(H^0(C,\omega_C))\oplus\C\oplus\End(H^0(C,\omega_C)) \ . \]
we have
\[\tau(a_1,f_1,a_2,f_2) \ = \ a_1+a_2+tr(f_1)+tr(f_2) \ . \]
It follows that
\[{\rm cl}_\Gamma \ = \ (1,I_{H^0(C,\omega_C)},-1,-I_{H^0(C,\omega_C)}) \ . \]
\end{subsub}

\sepprop

\begin{subsub} Proof of Theorem \ref{theo4} -- \rm The case \ \m{L=\omega_C} \ 
is obvious. So we suppose that $C$ is not hyperelliptic and \ \m{\deg(L)\leq 
2-2g}. The transpose of the evaluation morphism \
\m{\ko_C\ot H^0(C,\omega_C)\to\omega_C} \ is injective and its cokernel $F$ is 
stable (from \cite{pa_ra}, Corollary 3.5). From the exact sequence
\[0\lra T_C\ot L\lra L\ot H^0(C,\omega_C)^*\lra F\ot L\lra 0 \ , \]
we have an exact sequence
\xmat{0\ar[r] & H^0(C,F\ot L)\ar[r] & H^1(C,T_C\ot L)\ar[r]^-\eta &
L(H^0(C,\omega_C),H^1(C,L)) \ . }
Since \ \m{\deg(F\ot L)\leq 0} \ and \m{F\ot L} is stable, we have \ 
\m{H^0(C,F^*\ot L)=\nsp}, and $\eta$ is injective. Since \ \m{\sigma_Y\not=0} 
($Y$ is not trivial), we have \ \m{\Delta(\ko_{C\times C}(-\Gamma))\not=0}. 
This proves Theorem \ref{theo4}.
\end{subsub}

\end{sub}

\sepsub

\Ssect{Moduli spaces of vector bundles of rank \m{r\geq 2} and degree 
$-1$}{mod_rg2}

\begin{subsub}{\bf Theorem :}\label{theo15} Suppose that $C$ is not 
hyperelliptic, $Y$ not trivial and \ \m{\deg(L)\leq 2-2g}.
Then \m{{\bf M}_Y(r,-1)} is not a vector bundle over \m{M(r,-1)}.
\end{subsub}
\begin{proof}
Let $E$ be a stable vector bundle on $C$ of rank $r$ and degree $0$ and \m{x\in 
C}. Then for every non zero linear form \ \m{\varphi:E_x\to\C}, the kernel of 
the induced morphism \m{\tilde{\varphi}:E\to\ko_x} is a stable vector bundle: 
let \m{F\subset\ker(\tilde{\varphi})} be a subbundle such that \ 
\m{0<\rk(F)<r}; then since \m{F\subset E} and $E$ is stable, we have 
\m{\deg(F)<0}. Hence \ \m{\dsp\frac{\deg(F)}{\rk(F)}<-\frac{1}{r}}.

Let \m{\P=\P(E^*)} and \m{\pi:\P\to C}, \m{p_\P:\P\times C\to C}, 
\m{p_C:\P\times C\to C} the projections. Let \Nligne \m{\Delta\subset C\times 
C} \ be the diagonal. We have a canonical obvious surjective morphism on 
sheaves on \m{\P\times C}
\[{\bf p}:p_C^*(E)\ot p_\P^*(\ko_\P(-1))\lra\ko_{(\pi\times p_C)^{-1}(\Delta)} 
\ . \]
For every \m{x\in C} and non zero \m{\varphi\in E_x^*}, \m{{\bf p}_{\C\varphi}}
is the morphism \m{E\ot\C\varphi\to\ko_x} induced by $\varphi$. The kernel of 
$\bf p$ is a vector bundle $\E$ on \m{\P\times C} that is a family of stable 
vector bundles on $C$, of rank $r$ and degree \m{-1}, parametrized by $\P$. 

Let $\ke$ be a universal bundle on \m{M(r,-1)\times C}.
Since \m{r\geq 2}, \m{E^*} has a rank 1 subbundle $D$, inducing a section $s$ 
of $\pi$. Then \m{\F=s^*(\E)} is a family of stable vector bundles on $C$, of 
rank $r$ and degree \m{-1}, parametrized by $C$. Hence there exists a morphism 
\ \m{\tau:C\to M(r,-1)} and a line bundle $H$ on $C$ such that \ 
\m{\F\simeq(\tau\times I_C)^*(\ke)\ot q^*(H)} (where \m{q:C\times C\to C} is 
the first projection). 

We have a canonical immersion \ \m{i:C\to\Pic^{-1}(C)}, sending $x$ to
\m{\ko_C(-x)\ot\det(E)}. Let \ \m{{\bf det}:M(r,-1)\to\Pic^{-1}(C)} \ be 
the determinant morphism. Then \ \m{{\bf det}\circ\tau} \ is the identity 
morphism (from $C$ to \m{i(C)=C}). It follows that $\tau$ is an immersion.

Now suppose that \m{{\bf M}_Y(r,-1)} is a vector bundle over \m{M(r,-1)}. By 
Theorem \ref{theo3}, $\ke$ can be extended to a vector bundle on \ \m{{\bf 
M}_Y(r,-1)\times Y}, hence \m{\ke_{|\tau(C)\times C}} can be extended to a 
vector bundle on \m{\tau(C)\times Y}. Hence \m{\det(\ke_{|\tau(C)\times C})}
can be extended to a line bundle on \m{i(C)\times C}. But this is impossible by 
the proof of Theorem \ref{theo4}.
\end{proof}

\end{sub}

\sepsub

\Ssect{Picard groups of primitive multiple schemes of multiplicity 3}{geq0}

Let \m{C_3} be an extension of \m{C_2} of multiplicity 3. We use the notations 
of \ref{pic}. Let \ \m{{\bf P}=\Pic^d(C)}. Every line bundle on \m{C_2} can be 
extended to \m{C_3}, since the obstructions lie in $H^2$ of vector bundles on 
$C$. Hence we have \ \m{\Gamma^{\bf P}(C_3)=\Pic^{d}(C_2)}.

Let \m{\Pic^{d}(C_3)} denote the component of \m{\Pic(C_3)} of line bundles 
whose restriction to $C$ is of degree $d$. From \ref{pic} we see that 
\m{\Pic^{d}(C_3)} is a smooth irreducible variety, and an affine 
bundle over \m{\Pic^{d}(C_2)} with associated vector bundle \ 
\m{\ko_{\Pic^{d}(C_2)}\ot\big(H^1(C,L^2)/\imm(\delta^0)\big)} , where \ 
\m{\delta^0:H^0(\ko_{C_2})\to H^1(C,L^2)} \ is the map coming from the exact 
sequence \Nligne \m{0\to L^2\to\ko_{C_3}\to\ko_{C_2}\to 0} .

We have an exact sequence
\[0\lra H^0(C,L)\lra H^0(\ko_{C_2})\lra H^0(\ko_C)=\C\lra 0 \ , \]
so that \m{\imm(\delta^0)} is the same as the image of the restriction of 
\m{\delta^0} to \m{H^0(C,L)}. 

Now suppose that we have another extension \m{C'_3} of \m{C_2} to multiplicity 
$3$. Then \m{\ki_{C,C_3}} and \m{\ki_{C,C'_3}} are two extensions of $L$ to a 
line bundle on \m{C_2}, corresponding to \ 
\m{\sigma,\sigma'\in\Ext^1_{\ko_{C_2}}(\ko_C,L)} \ respectively. We have also 
two maps \ \m{\delta^0,{\delta'}^0:H^0(C,L)\to H^1(C,L^2)}.

According to 
\ref{ext_4}, we have \ \m{\sigma'-\sigma=\Psi(\mu)}, for some \ \m{\mu\in 
H^1(C,L)}. Let \Nligne \m{\chi:H^0(C,L)\ot H^1(C,L)\to H^1(C,L^2)} \ be the 
canonical product.

From Lemma \ref{lem17} follows easily

\sepprop

\begin{subsub}\label{prop18}{\bf Proposition: } We have, for every \ \m{\eta\in 
H^0(C,L)}, \m{({\delta'}^0-\delta^0)(\eta)=\chi(\eta\ot\mu)}.
\end{subsub}

\sepprop

This implies that it is possible that \m{\imm(\delta^0)\not=0}. If \m{C_2} is 
the trivial ribbon (cf. \ref{PMS_def}), then \m{\delta^0=0}, so we obtain a 
formula for \m{{\delta'}^0}.

\end{sub}

\sepsec

{\bf Acknowledgements.} I am grateful to the anonymous referee for giving some 
useful remarks.

\vskip 4cm

\end{document}